\documentclass[12pt, twoside, leqno, letterpaper]{huthesis}

\usepackage{amssymb, amsthm, amsmath, graphicx, makeidx, enumerate}
\usepackage{color} \usepackage[all]{xy}

\makeindex

\def\dotfil{\leaders\hbox to.6em{\hss .\hss}\hfil}%

\newcommand{\dfi}[1]{{\em #1}}
\newcommand{\df}[2]{{\index{#1}{\em #2}}}
\newcommand{\comma}{,}

\newcommand{\propwidth}{350pt}

\DeclareMathOperator{\st}{st} \DeclareMathOperator{\bst}{bst}
 \DeclareMathOperator{\dom}{dom}
\DeclareMathOperator{\rel}{rel} \DeclareMathOperator{\Int}{Int}
\DeclareMathOperator{\Bd}{Bd} \DeclareMathOperator{\Cl}{Cl}
\DeclareMathOperator{\im}{im} \DeclareMathOperator{\jn}{jn}
 \DeclareMathOperator{\mesh}{mesh}

\newcommand{\nhom}{{\ \stackrel{n}{\sim}\ }}

\newcommand{\E}[1]{{\protect\mathcal{E}_{#1}}}
\newcommand{\F}[1]{{\protect\mathcal{F}_{#1}}}

\newcommand{\G}[1]{{\protect\mathcal{G}_{#1}}}
\newcommand{\N}[1]{{\protect\mathcal{N}_{#1}}}
\newcommand{\U}[1]{{\protect\mathcal{U}_{#1}}}
\newcommand{\V}[1]{{\protect\mathcal{V}_{#1}}}
\newcommand{\W}[1]{{\protect\mathcal{W}_{#1}}}

\newcommand{\BF}{{\protect\mathcal{B}_\F{}}}
\newcommand{\BFn}{{\protect\mathcal{B}^{(n)}_\F{}}}

\theoremstyle{plain}
\newtheorem{theorem}{Theorem}[chapter]
\newtheorem{corollary}[theorem]{Corollary}
\newtheorem{lemma}[theorem]{Lemma}
\newtheorem{proposition}[theorem]{Proposition}

\newtheorem*{characterization theorem}{Characterization theorem}
\newtheorem*{twierdzenie o charakteryzacji}{Twierdzenie o charakteryzacji}
\newtheorem*{rigidity theorem}{Topological rigidity theorem}
\newtheorem*{twierdzenie o sztywnosci}{Twierdzenie o sztywno\'sci}
\newtheorem*{z-set unknotting theorem}{$Z$-Set unknotting theorem}
\newtheorem*{open embedding theorem}{Open embedding theorem}
\newtheorem*{local z-set unknotting theorem}{Local $Z$-set unknotting theorem}
\newtheorem*{sum theorem}{Sum theorem}
\newtheorem*{main theorem}{Main theorem}
\newtheorem*{whitehead's characterization}{Whitehead's characterization}
\newtheorem*{corollary*}{Corollary}
\newtheorem*{nerve theorem}{Nerve theorem}
\newtheorem*{carrier theorem}{Carrier theorem}
\newtheorem*{theorem*}{Theorem}
\newtheorem*{proposition*}{Proposition}
\newtheorem*{pump up theorem}{Theorem~\ref{thm:pump up the regularity}}
\newtheorem*{Torunczyk's lemma}{Toru\'nczyk's lemma}

\theoremstyle{definition}
\newtheorem{definition}[theorem]{Definition}
\newtheorem*{definition*}{Definition}
\newtheorem*{construction*}{Construction}

\newtheorem{example}[theorem]{Example}

\newtheorem{remark}[theorem]{Remark}

\begin{document}

\title{Characterization and topological rigidity of N\"obeling manifolds}
\author{Andrzej Nag\'orko}
\degreemonth{January}
\degreeyear{2006}
\degree{Doctor of Philosophy}
\field{Mathematics}
\department{Mathematics}
\advisor{Henryk Toru\a'nczyk}

\maketitle
\copyrightpage

\begin{sspabstract}

  We develop a theory of N\"obeling manifolds similar to the theory of
  Hilbert space manifolds. We show that it reflects the theory of
  Menger manifolds developed by M.~Bestvina~\cite{bestvina1988} and is
  its counterpart in the realm of complete spaces. In particular we
  prove the N\"obeling manifold characterization conjecture.
  
  We define the $n$-dimensional universal N\"obeling space $\nu^n$ to
  be the subset of $R^{2n+1}$ consisting of all points with at most
  $n$ rational coordinates.  To enable comparison with the infinite
  dimensional case we let $\nu^\infty$ denote the Hilbert space. We
  define an $n$-dimensional N\"obeling manifold to be a Polish space
  locally homeomorphic to $\nu^n$. The following theorem for $n =
  \infty$ is the characterization theorem of
  H.~Toru\'nczyk~\cite{torunczyk1981}. We establish it for $n <
  \infty$, where it was known as the N\"obeling manifold
  characterization conjecture.

  \vspace{1mm}
  \noindent {\bf Characterization theorem } { \em
  An $n$-dimensional Polish $ANE(n)$-space is a N\"o\-be\-ling
  manifold if and only if it is strongly universal in dimension
  $n$.
  }
  \vspace{1mm}
  
  The following theorem was proved by D. W. Henderson and R.
  Schori~\cite{hendersonschori1970} for $n = \infty$. We establish it in the
  finite dimensional case.

  \vspace{1mm}
  \noindent {\bf Topological rigidity theorem.} { \em
  Two $n$-dimensional N\"obeling manifolds are homeomorphic if and
  only if they are $n$-homotopy equivalent.
  }
  \vspace{1mm}
  
  We also establish the open embedding theorem, the $Z$-set unknotting
  theorem, the local $Z$-set unknotting theorem and the sum theorem for
  N\"obeling manifolds.

\end{sspabstract}
\newpage\thispagestyle{empty}\ 
\vspace{ 0.7cm}
{ \center  \bf \large Key words and phrases.

}

\vspace{ 5mm} N\"obeling manifold classification, N\"obeling space
characterization, $Z$-set unknotting theorem, open embedding theorem, carrier,
nerve theorem, regular cover, semiregular cover

\vspace{ 10mm}

{ \center \bf \large 2000 Mathematics Subject Classification. 

} 

\vspace{ 5mm}
Primary 55M10, 54F45; Secondary 54C20.

\begin{plsspabstract}
  
  Celem niniejszej rozprawy jest rozwini\c{e}cie sko\'nczenie wymiarowego
  odpowiednika klasycznego dzia{\l}u topologii niesko\'nczenie wymiarowej:
  teorii rozmaito\'sci modelowanych na przestrzeni Hilberta. W
  szczeg\'olno\'sci dowodzimy sko\'nczenie wymiarowej wersji twier\-dze\-nia
  Toru\'nczyka o charakteryzacji.
  \begin{twierdzenie o charakteryzacji}
    Polska $n$-wymiarowa $ANE(n)$-przestrze\'n jest roz\-ma\-i\-to\-\'sci\c{a}
    N\"obelinga jedynie je\'sli jest mocno uniwersalna w wymiarze $n$.
  \end{twierdzenie o charakteryzacji}

  Dowodzimy te\.z sko\'nczenie wymiarowej wersji twierdzenia Hendersona
  klasyfikuj\c{a}ego rozmaito\'sci modelowane na przestrzeni Hilberta.

  \begin{twierdzenie o sztywnosci}
    Dwie $n$-wymiarowe rozmaito\'sci N\"obelinga s\c{a} homeomorficzne jedynie
    je\'sli s\c{a} $n$-homotopijnie r\'ownowa\.zne.
  \end{twierdzenie o sztywnosci}
  
  Poza tym dowodzimy twierdzenia o jednakowym po{\l}o\.zeniu $Z$-zbior\'ow,
  jego lokalnej wersji, twierdzenia o otwartym zanurzeniu i twierdzenia o
  sumie dla rozmaito\'sci N\"obelinga.
\end{plsspabstract}
\newpage\thispagestyle{empty}\ 

\newpage
\addcontentsline{toc}{section}{Table of Contents}
\tableofcontents

\begin{acknowledgments}\hsp
  Completing this thesis was an overwhelming experience and would be
  impossible withoth the help of others. Foremost I would like to thank
  my family for support and encouragement during these years.

  I am grateful to my thesis supervisor, professor H. Toru\'nczyk, who
  helped with his enormous knowledge and was always available to talk
  and review my writings.  I am grateful to reviewers of the first
  version of this thesis, professors W. Marciszewski and S. Spie\.z,
  who helped to improve this text a lot.

  I would like to thank the organizers of the Topology seminar at the
  Institute of Mathematics of Polish Academy of Sciences - professors
  J. Krasinkiewicz and S. Spie\.z - who gave me an opportunity to give
  a long series of talks on the subject title. I'm grateful to
  professors T. Dobrowolski and M. Sobolewski who also attended these
  talks.

  I'm grateful to D.~Michalik, R.~G\'orak, M.~Sawicki and professors
  E.~Pol and A.~Chigogidze for helping me out with various questions I
  have had.
  
  The starting point and the original inspiration for this work
  was a paper by M. Levin, K. Kawamura and E. D. Tymchatyn, who proved
  the characterization theorem for one-dimensional N\"obeling
  spaces~\cite{kawamuralevintymchatyn1997}. I had the pleasure of
  meeting professors M. Levin and K. Kawamura at the International
  Conference and Workshops on Geometric Topology in June 2005, where
  the results of this thesis were presented.

\end{acknowledgments}

\begin{citations}

Contents of chapter~\ref{ch:carrier and nerve theorems} were published as a
part of
\begin{quote}
  ``Carrier and nerve theorems in the extension theory'', A.~Nag\'orko,
  Proc. Amer. Math. Soc. 135(2):551-558, 2007.
\end{quote}

\end{citations}
\newpage\thispagestyle{empty}\ 

\dedication

%
%

\ \vspace{ 7.7in}
{ \flushright \em \Large To my father \\ }

\newpage\thispagestyle{empty}\ 

\newpage
\ssp
\startarabicpagination

\part{Introduction and preliminaries}
\chapter{Introduction}

{ \em \raggedleft \small ,,Perhaps this [characterization conjecture
for N\"obeling manifolds]\\ is one of the most exciting
problems of general topology''. \\
\vspace{1mm} A. N. Dranishnikov~\cite{dranishnikov1997} \\ }
\vspace{ 5 mm}

The aim of this paper is to develop a theory of N\"obeling manifolds
similar to the theory of Hilbert space manifolds. In particular,
the N\"obeling manifold characterization conjecture is proven.

We define the $n$-dimensional universal \df{N\"obeling
  space}{N\"obeling space $\nu^n$} to be the set of all points of
$\mathbb{R}^{2n+1}$ with at most $n$ rational coordinates. To enable
comparison with the infinite dimensional case we let $\nu^\infty$
denote the Hilbert space~$\mathbb{R}^\infty$. We define an
\df{N\"obeling manifold}{$n$-dimensional N\"obeling manifold} to be a
Polish space locally homeomorphic to $\nu^n$. The following theorem
for $n = \infty$ is the characterization theorem of
H.~Toru\'nczyk~\cite{torunczyk1981}. We establish it for $n < \infty$,
where it was known as the N\"obeling manifold characterization
conjecture~\cite{chigogidze1996, dobrowolski1990, kawamura2002,
  west1990}.

\begin{characterization theorem}
  \index{theorem!characterization} An $n$-dimensional Polish space is
  a N\"obeling manifold if and only if it is an absolute neighborhood
  extensor in dimension~$n$ that is strongly universal in dimension
  $n$.
\end{characterization theorem}

D. W. Henderson and R. Schori proved~\cite{hendersonschori1970} the
following theorem for $n = \infty$. We give a proof of the finite
dimensional case.

\begin{rigidity theorem}
  \index{theorem!topological rigidity}
  Two $n$-dimensional N\"obeling manifolds are homeomorphic if and
  only if they are $n$-homotopy equivalent.
\end{rigidity theorem}

Topological rigidity and characterization theorems imply a
characterization theorem for N\"obeling spaces: every
$n$-dimensional strongly universal Polish absolute extensor in
dimension $n$ is homeomorphic to the $n$-dimensional N\"obeling
space. Another consequence of topological rigidity theorem is the
open embedding theorem. For $n = \infty$ it was proved by
Henderson~\cite{henderson1970}. We prove it in the finite
dimensional case.

\begin{open embedding theorem}
  \index{theorem!open embedding}
  Every $n$-dimensional N\"obeling manifold is homeomorphic to
  an open subset of the $n$-dimensional N\"obeling space.
\end{open embedding theorem}

The topological rigidity theorem has the following generalization.
Again, it was known for $n = \infty$ by results of R.~D.~Anderson
and J.~D.~McCharen~\cite{andersonmccharen1970} and we prove it in
the finite dimensional case.

\begin{z-set unknotting theorem}
  \index{theorem!Z-set unknotting@$Z$-set unknotting} A homeomorphism
  between $Z$-sets in $n$-dimensional N\"obeling manifolds extends to
  an ambient homeomorphism if and only if it extends to an
  $n$-homotopy equivalence.
\end{z-set unknotting theorem}

We also prove a local version of the $Z$-set unknotting theorem, for
$n < \infty$. It was proved for $n = \infty$ independently by W.
Barit~\cite{barit1969} and by Cz.~Bessaga and
A.~Pe{\l}czy\'nski~\cite{bessagapelczynski1970}.

\begin{local z-set unknotting theorem}
  \index{theorem!local $Z$-set unknotting} For every open cover $\U{}$
  of a N\"obeling manifold there exists an open cover $\V{}$ such that
  every homeomorphism between $Z$-subsets of the manifold that is
  $\V{}$-close to the inclusion extends to a homeomorphism of the
  entire manifold that is $\U{}$-close to the identity.
\end{local z-set unknotting theorem}

By an example of R.~B.~Sher a space that is an union of two Hilbert
cubes that meet in a Hilbert cube needs not to be a Hilbert
cube~\cite{sher1977}. We show that no such example exists in the
theory of N\"obeling spaces.

\begin{sum theorem}
  \index{theorem!sum}
  If a space $X$ is an union of two closed $n$-dimensional N\"obeling manifolds
  whose intersection is also an $n$-dimensional N\"obeling manifold, then $X$
  is an $n$-dimensional N\"obeling manifold.
\end{sum theorem}

As we shall see in chapter~\ref{ch:proof of the main theorem} the
proof of the above mentioned theorems can rather easily be reduced to
a proof of the main theorem, stated below. We say that a space is an
\df{abstract N\"obeling manifold}{abstract N\"obeling manifold} if it
is an $n$-dimensional Polish absolute neighborhood extensor in
dimension~$n$ that is strongly universal in dimension~$n$, i.e. if it
satisfies the conditions of the characterization theorem.

\begin{main theorem}
  \index{theorem!main} A homeomorphism between $Z$-sets in abstract
  N\"obeling manifolds extends to an ambient homeomorphism if and only
  if it extends to an $n$-homotopy equivalence.
\end{main theorem}

The theory developed in this paper reflects the theory of Menger
manifolds developed by M.~Bestvina~\cite{bestvina1988}.  He proved
finite dimensional analogues of the characterization of the Hilbert
cube given by Toru\'nczyk~\cite{torunczyk1980} and of the
topological rigidity of Hilbert cube manifolds given by
J.~West~\cite{west1970} and T.~A.~Chapman~\cite{chapman1974}. Proofs
given in the present paper does not admit direct generalizations neither 
to the compact case nor to the complete, infinite dimensional case.

Special cases of the characterization and the $Z$-set unknotting
theorems, for the case of N\"obeling spaces rather than that of
N\"obeling manifolds, were established by several authors. The
one-dimensional cases were proved by K.~Kawamura, M.~Levin and
E.~D.~Tymchatyn~\cite{kawamuralevintymchatyn1997}. S.~Ageev
announced a proof of the characterization theorem and of a special
case of the $Z$-set unknotting theorem~\cite{ageev2000} but recently
he declared it contained a gap~\cite{ageev2005}. The Z-set
unknotting theorem for a certain model of a N\"obeling space was
proved by Levin~\cite{levin2005}. Levin's result and those of this
thesis were presented at the International Conference and Workshops
on Geometric Topology, 3-10 July 2005. Very recently, Levin posted a
preprint~\cite{levin2006} containing a proof of the characterization
theorem for N\"obeling spaces of all dimensions.

\subsection*{How the paper is organized}

The text is divided into three parts. In the first part we recall the
general results used. In the second part we reduce the proof of the
theorems menitoned above to the proofs of two theorems (\ref{thm:pump
  up the regularity} and~\ref{thm:retraction onto a complex}) that
give sufficient conditions for the existence of covers of abstract
N\"obeling manifolds that resemble triangulations of simplicial
complexes. In the third part these two theorems are proved.



\chapter{Preliminaries}
\label{ch:preliminaries}

To facilitate reading of this paper we recall the general
results and notions used throughout~it. It makes the paper
self-contained in a sense that, with few exceptions, in the next
chapters there are no references to external sources.

The image and the domain of a function $f$ are denoted by \df{im
  f@$\im f$}{$\im f$} and by \df{dom f@$\dom f$}{$\dom f$}
respectively. We write $f(A)$ for $f(A \cap \dom f)$ for each map $f$
and each set $A$, i.e. we do not require $A$ to be a subset of the
domain of $f$.  The unit interval is denoted by~$[0,1]$ and the
symbol~$I$ is used exclusively to denote a set of indices of a
collection of sets.  If $\F{} = \{ F_i \}_{i \in I}$ is an indexed
collection of sets and $J$ is a non-empty subset of the indexing set
$I$, then we let \df{FJ@$F_J$}{$F_J = \bigcap_{j \in J} F_j$}.  We use
the following symbols to denote standard topological operations:
$\Cl_X A = \Cl A$ denotes the \df{Cl A@$\Cl_X A$}{closure} of a subset
$A$ of a space $X$, $\Int_X A = \Int A$ denotes the \df{Int A@$\Int_X
  A$}{interior} of a subset~$A$ of a space $X$ and $\partial X$
denotes the \df{d X@$\partial X$}{geometrical boundary} of an
euclidean manifold~$X$.

We say that a space is \df{Polish space}{Polish} is it is separable
and completely metrizable. A \df{map}{map} is a continuous function
between two topological spaces. We say that sets \df{intersecting
  sets}{intersect} if their intersection is non-empty.

\section{Covers and interior covers}

We say that a collection $\F{}$ of sets is \df{indexed
  collection}{indexed}, whenever there is a chosen \df{indexing
  set}{indexing set} $I$ and an \df{indexing function}{indexing
  function} that assigns to each element $i$ of $I$ an element~$F_i$
of~$\F{}$. We shall write $\F{} = \{ F_i \}_{i \in I}$ for short. Note
that the indexing function is onto, but needs not to be one-to-one.
Nonetheless, we will distinguish between elements of $\F{}$ that are
equal, but have different indices. Unless otherwise stated, we shall
assume that every collection of sets is indexed.

\begin{definition*}
  A \df{cover}{cover} of a topological space $X$ is a collection of
  subsets of $X$ whose union equals to $X$. An
  \df{cover!interior}{interior cover} of $X$ is a collection of
  subsets of $X$ whose interiors cover~$X$.
\end{definition*}

Symbols $\F{}, \G{}$ will usually denote closed covers and symbols
$\U{}, \V{}$ will usually denote open covers.

\begin{definition*}
  Elements of two indexed collections of sets are \df{corresponding
    elements}{corresponding} if they have equal indices. Covers with
  equal indexing sets are \df{cover!isomorphic}{isomorphic} if for
  each set of intersecting elements in one cover the corresponding
  elements in the other cover intersect.
\end{definition*}

In more formal terms, a cover $\F{} = \{ F_i \}_{i \in I}$ is
isomorphic to a cover $\G{} = \{ G_i \}_{i \in I}$ if for each
non-empty $J \subset I$ the intersection $F_J = \bigcap_{j \in J} F_j$
is non-empty if and only if the intersection $G_J = \bigcap_{j \in J}
G_j$ is non-empty. For every $i \in I$ elements~$F_i$ of~$\F{}$
and~$G_i$ of~$\G{}$ are corresponding.

\begin{definition}\label{def:restricted cover}
  Let $\F{} = \{ F_i \}_{i \in I}$ and $\G{} = \{ G_i \}_{i \in I}$ be
  collections of subsets of a space $X$, indexed by the same indexing
  set $I$.
  \begin{enumerate}[(a)]
  \item We say that \df{swelling}{$\G{}$ is a swelling of $\F{}$}, if
    it is isomorphic to~$\F{}$ and if every element of~$\F{}$ is a
    subset of the corresponding element of $\G{}$.
  \item We say that \df{shrinking}{$\F{}$ is a shrinking of $\G{}$},
    if $\G{}$ is a swelling of $\F{}$.
  \item For each subset $A$ of~$X$, we let $\{ F_i \cap A \}_{i \in
      I}$ be \df{restriction of a collection}{the restriction of
      $\F{}$ to $A$} and denote it by \df{F/A}{$\F{} / A$}.
  \item For each subset $A$ of $X$, we say that \df{cover!equal on a
      set}{$\F{}$ is equal to $\G{}$ on $A$}, if $F_i \cap A = G_i
    \cap A$ for each $i$ in $I$, i.e. if the restrictions of $\F{}$
    and $\G{}$ to $A$ are equal.
  \item We say that \df{refinement}{$\F{}$ refines a collection
      $\mathcal{H}$} of subsets of $X$, if every element of $\F{}$ is
    a subset of an element of $\mathcal{H}$.
  \item For each point $x$ in $X$, we let \df{st F x@$\st_\F{}
      x$}{$\st_\F{} x$} be the union of all elements of $\F{}$ that
    contain~$x$.  We say that $\st_\F{} x$ is a \df{star of a
      point}{star of $x$ in $\F{}$}.
  \item For each subset $A$ of $X$, we let \df{st F A@$\st_\F{}
      A$}{$\st_\F{} A$} be the union of $A$ with all elements of
    $\F{}$ that intersect $A$.  We say that $\st_\F{} A$ is a \df{star
      of a set}{star of $A$ in $\F{}$}.
  \item For each collection $\mathcal{H}$ of subsets of $X$, we let
    \df{st F G@$\st_\mathcal{H} \F{}$}{$\st_\mathcal{H} \F{} = \{
      \st_\mathcal{H} F_i \}_{i \in I}$}. We say that $\st_\mathcal{H}
    \F{}$ is a \df{star of a collection}{star of $\F{}$ in
      $\mathcal{H}$}.
  \item For each collection $\mathcal{H}$ of subsets of $X$ and each
    $m > 0$, we let $\st^m_\mathcal{H} \F{} = \st_\mathcal{H}
    \st^{m-1}_\mathcal{H} \F{}$ and $\st^0_\mathcal{H} \F{} = \F{}$.
    We say that $\st^m \F{}$ is \df{mth star@$m$th star}{the $m$th
      star of $\F{}$ in $\mathcal{H}$}.
  \item For each $m > 0$, we let $\st^m \F{} = \st_\F{} \st^{m-1}
    \F{}$ and $\st^0 \F{} = \F{}$. We say that $\st^m \F{}$ is
    \dfi{the $m$th star of $\F{}$}.
  \end{enumerate}
\end{definition}

Note that $\F{}$ and $\st_\mathcal{H} \F{}$ have the same indexing
set, hence it makes sense to talk about corresponding elements of
$\F{}$ and $\st_\mathcal{H} \F{}$.

\begin{lemma}
  \label{lem:star of a star}
  If $k$ and $m$ are nonnegative integers and $\F{}$ is a collection
  of sets, then
  \[ 
  \st^k_\F{} \st^m \F{} = \st^{k + m} \F{} \text{ and } \st^k \st^m
  \F{} = \st^{2km + m + k} \F{}.
  \]
\end{lemma}

\section{Absolute extensors}

Many characterizations of absolute extensors for various classes of
topological spaces made it into quite a mess in notation. We are
going to settle for the following definitions.

\begin{definition*}
  We say that a space $X$ is an \df{absolute extensor in
    dimension~$n$}{absolute extensor in dimension~$n$} if it is a
  metric space and if every map into $X$ from a closed subset of an
  $n$-dimensional metric space extends over the entire space. The
  class of absolute extensors in dimension $n$ is denoted by $AE(n)$
  and its elements are called \df{AE(n)@$AE(n)$}{$AE(n)$-spaces}.
\end{definition*}

\begin{definition*}
  We say that a metric space $X$ is an \df{absolute neighborhood
    extensor in dimension~$n$}{absolute neighborhood extensor in
    dimension $n$} if it is a metric space and if every map into $X$
  from a closed subset $A$ of an $n$-dimensional metric space extends
  over an open neighborhood of $A$. The class of absolute neighborhood
  extensors in dimension $n$ is denoted by $ANE(n)$ and its elements
  are called \df{ANE(n)@$ANE(n)$}{$ANE(n)$-spaces}.
\end{definition*}

The following characterization of the class of absolute neighborhood
extensors in dimension~$n$ is proven in~\cite{hu1965}, with references
to original papers by K.~Borsuk and J.~Dugundji.

\begin{definition}\label{def:u-close}
  We say that a map $g$ is \df{U-@$\U{}$-!close}{$\U{}$-close} to a map~$f$,
  or that~$g$ is a \df{U-@$\U{}$-!approximation}{$\U{}$-approximation} of~$f$,
  if~$\U{}$ is an open cover of the common codomain of~$f$ and~$g$ and the
  collection $\{ f^{-1}(U) \cap g^{-1}(U) \}_{U \in \U{}}$ covers the common
  domain of $f$ and $g$.
\end{definition}

\begin{theorem}
 \label{thm:ane characterization}
 The following conditions are equivalent for every metric space~$X$ and every
 integer $n < \infty$:
  \begin{enumerate}
    \item $X$ is an $ANE(n)$-space.
    \item \label{item:lcn} $X$ is locally $k$-connected for every $k < n$. In
      other words, for each point $x$ in $X$ and every neighborhood $U$ of
      $x$ there exists a smaller neighborhood $V$ of $x$ such that every map
      from the boundary of a $(k+1)$-dimensional ball into $V$ extends to a
      map of the entire ball into~$U$.
    \item \label{item:v-close u-homotopic} for each open cover $\U{}$ of $X$
      there is an open cover $\V{}$ of $X$ such that any two $\V{}$-close maps
      into $X$, defined on an at most $(n-1)$-dimensional Polish space, are
      $\U{}$-homotopic. Recall that two maps are
      \df{U-@$\U{}$-!homotopic}{$\U{}$-homotopic} if they are homotopic by a
      homotopy whose paths refine $\U{}$.
  \end{enumerate}
\end{theorem}

The third condition of theorem~\ref{thm:ane characterization} implies
the following proposition, as the space of maps (taken with the $\sup$
topology) from a $k$-dimensional sphere into a separable metric space
is separable.
\begin{proposition}
  \label{pro:countable homotopy groups}
  If a space is a separable absolute neighborhood extensor in
  dimension $n$, then its homotopy groups of dimensions less than~$n$
  are countable.
\end{proposition}

A characterization of the class of absolute extensors in dimension~$n$
was given by J.~Dugundji~\cite{dugundji1958} in a form of the
following theorem.

\begin{theorem}
  \label{thm:dugundji characterization}
  A space is an absolute extensor in dimension~$n$ if and only if it
  is an absolute neighborhood extensor in dimension $n$ and its
  homotopy groups of dimensions less than $n$ vanish.
\end{theorem}

Theorem~\ref{thm:sum theorem for ane-spaces} is a closed sum theorem
for absolute neighborhood extensors in dimension~$n$.

\begin{theorem}[{\cite[p. 49]{hu1965}}] 
  \label{thm:sum theorem for ane-spaces}
  A space is an absolute neighborhood extensor in dimension $n$
  whenever it can be represented as an union of two closed
  $ANE(n)$-subspaces that meet in an $ANE(n)$-subspace.
\end{theorem}

We cite the following theorem after~\cite[Proposition 4.1.7, p.
131]{chigogidze1996}, with $(*)_{n-1}$ renamed to $(*)_n$.

\begin{theorem}
  \label{thm:n-homotopy extension theorem}

  Every open cover $\U{}$ of an $ANE(n)$-space~$Y$ has an open
  refinement~$\V{}$ such that the following condition holds:
  \begin{enumerate}
  \item[$(*)_n$] For any at most $n$-dimensional metric space $B$, any
    closed subspace $A$ of $B$, and any two $\V{}$-close maps $f, g
    \colon A \to Y$, if $f$ has an extension $\tilde f \colon B \to
    Y$, then~$g$ has an extension $\tilde g \colon B \to Y$ that is
    $\U{}$-close to $\tilde f$.
  \end{enumerate}
\end{theorem}

\begin{lemma}
  \label{lem:extension property}
  For every open neighborhood $U$ of a closed $AE(n)$-subset of an at
  most $n$-dimensional $ANE(n)$-space there exists an open
  neighborhood $V$ such that every map from a closed subset of an at
  most $n$-dimensional metric space into $V$ extends over the entire
  space to a map into $U$.
\end{lemma}
\begin{proof}
  By theorem~\ref{thm:ane characterization}, $U$ is an absolute
  neighborhood extensor in dimension $n$.  By
  theorem~\ref{thm:n-homotopy extension theorem}, there exists an open
  cover $\V{}$ of $U$, for which $(*)_n$ is satisfied with $\U{} = \{
  U \}$.  Let $A$ be a closed subset of $U$ that is an $AE(n)$-space.
  Since $U$ is at most $n$-dimensional, there exists a retraction $r
  \colon U \to A$. Let
  \[
    V = \bigcup_{ W \in \V{} } r^{-1}(W) \cap W.
  \]
  It is an open neighborhood of $A$. We shall prove that it
  satisfies the assertion.

  Let $f$ be a map into $V$ from a closed subset $B$ of an at most
  $n$-dimensional metric space $X$. By the definition of $V$, maps $r
  \circ f$ and $f$ are $\V{}$-close. Since $A$ is an absolute extensor
  in dimension $n$, the map $r \circ f$ extends over $X$ to a map into
  $A$. Hence, by~$(*)_n$, $f$ extends over $X$ to a map into $U$. We
  are done.
\end{proof}

\section{Nerves of covers and barycentric stars}
\label{sec:simplicial complexes}

We define simplicial complexes in the usual way~\cite[p.  99]{hu1965}.
We do not allow infinite-dimensional simplices.  Unless otherwise
stated, simplicial complexes are endowed with the metric
topology~\cite[p.~100]{hu1965}. We do not assume that they are
countable or locally finite.

\begin{theorem}[{\cite[p.~107]{hu1965}}] \label{thm:complete
    complexes} Every simplicial complex is topologically complete.
  Every countable simplicial complex is separable.
\end{theorem}

\begin{definition*}
  Let $K$ be a simplicial complex. A \df{bst v@$\bst v$}{barycentric
    star~$\bst v$} of a vertex~$v$ of~$K$ is the union of all
  simplices of the first barycentric subdivision of~$K$ that
  contain~$v$.  We let \df{bst L@$\bst L$}{$\bst L$} be the union of
  barycentric stars of vertices of a subcomplex $L$ of~$K$.  We treat
  it as a subcomplex of the barycentric subdivision of~$K$ and define
  \df{bstm L@$\bst^m L$}{the $m$th barycentric star $\bst^m L$ of~$L$}
  recursively by the formula $\bst^m L = \bst(\bst^{m-1} L)$.
  {%
\begin{figure}[ht]
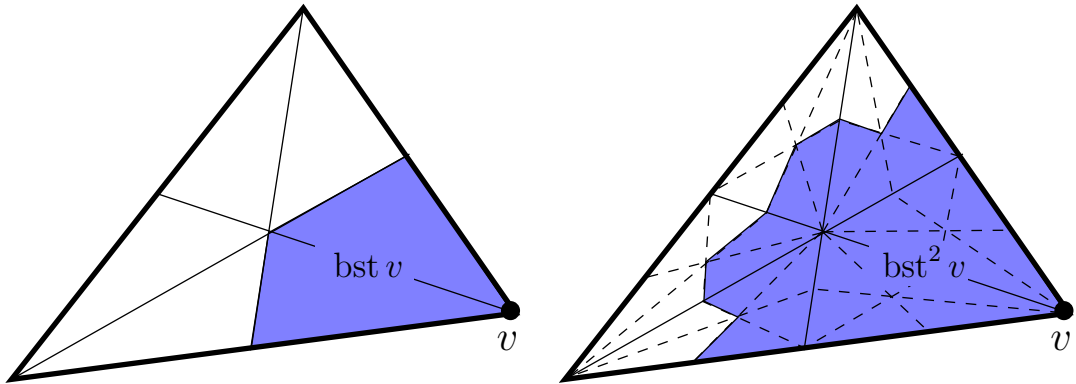

\begin{center}
$\begin{array}{c@{\hspace{5mm}}c}
\resizebox{0.45\textwidth}{!}{\input{psfigures/barycentric1.pstex_t}} &
\resizebox{0.45\textwidth}{!}{\input{psfigures/barycentric2.pstex_t}}
\end{array}$
\end{center}
\caption{{The first and the second
    barycentric star of a vertex~$v$ of a triangle.}}
\label{barycentric star}
\end{figure}
}
\end{definition*}

\begin{definition*}
  The nerve of an point-finite indexed collection~$\F{}$ of sets is a 
  simplicial
  complex whose set of vertices is in one-to-one correspondence with
  non-empty elements of~$\F{}$ and whose vertices span a simplex if
  and only if the corresponding elements in~$\F{}$ intersect.
\end{definition*}

\begin{lemma}
  \label{lem:bst isomorphism} A collection of barycentric stars of
  vertices of a simplicial complex is its closed locally finite
  cover. The nerve of a cover of a simplicial complex~$K$ by
  barycentric stars of its vertices is isomorphic to~$K$. The cover of
  the nerve of a cover~$\F{}$ by barycentric stars of its vertices is
  isomorphic to~$\F{}$.
\end{lemma}
\begin{proof}
  Let $K$ be a simplicial complex. Observe that
  \em\begin{enumerate}
  \item a barycentric star of a vertex of $K$ intersects a simplex in
    $K$ if and only if the vertex lies in the simplex,
  \item vertices of $K$ span a simplex if and only if their
    barycentric stars intersect.
  \end{enumerate}\em

  Choose a point in $K$ and let $\delta$ be a simplex that contains
  it. By (1), the union of barycentric stars of vertices contained in
  $K \setminus \delta$ is disjoint from $\delta$. By~\cite[p.~99]{hu1965}, every subcomplex of a simplicial complex is closed, hence
  the complement of this union is an open neighborhood of $\delta$
  (and of the chosen point). By the definition, it meets only
  barycentric stars of vertices contained in $\delta$, hence the cover
  of $K$ by barycentric stars of its vertices is locally finite. 

  Let $\mathcal{B}$ be a cover of $K$ by barycentric stars of its
  vertices. By~(2), a set of vertices of $K$ spans a simplex if and
  only if their barycentric stars intersect. By the definition, the
  elements of $\mathcal{B}$ intersect if and only if the corresponding
  vertices of the nerve of $\mathcal{B}$ span a simplex. Hence the
  second assertion of the lemma.

  By the definition, a set of elements of a collection $\F{}$ of sets
  has non-empty intersection if and only if the corresponding vertices
  in the nerve of $\F{}$ span a simplex. By~(2), these vertices span a
  simplex if and only if their barycentric stars intersect. Hence the
  third assertion of the lemma.
\end{proof}

\begin{lemma}
  \label{lem:contractible stars} A cover of a simplicial complex by
  barycentric stars of its vertices has the property that every
  intersection of its elements is either empty or contractible.
\end{lemma}
\begin{proof}
  Let $K$ be a simplicial complex and let $\mathcal{B}$ be a cover of
  $K$ by barycentric stars of its vertices.  By lemma~\ref{lem:bst
    isomorphism}, the nerve of $\mathcal{B}$ is isomorphic to $K$.
  Hence if elements of~$\mathcal{B}$ intersect, then the collection
  $\{ v_k \}$ of corresponding vertices of $K$ span a simplex~$\delta$
  in~$K$. By the definition, a barycentric star $\bst v_k$ is a full
  subcomplex of the first barycentric subdivision of $K$ that contains
  all barycenters of faces that contain $v_k$ (with a zero dimensional
  face $\{ v_k \}$ included). Hence the intersection $\bigcap \bst
  v_k$ is a full subcomplex of the first barycentric subdivision of
  $K$ that contains all barycenters of faces that contain $\delta$.
  Observe that this subcomplex is a cone with apex equal to the the
  barycenter of $\delta$. By~\cite[p. 116]{spanier1989}, it is
  contractible. We are done.
\end{proof}

\section{Strong universality}

In the present paper, we shall endow all function spaces with the
limitation topology, which embraces a concept of closeness introduced
in definition~\ref{def:u-close}.

\begin{definition*}
  We let $C(X, Y)$\index{C(X{\comma}Y)@$C(X\comma Y)$} denote the set
  of all maps from a space $X$ into a space~$Y$. For each map $f
  \colon X \to Y$ and for each open cover $\U{}$ of $Y$, we let $B(f,
  \U{})$\index{B(f{\comma}U)@$B(f{\comma}\U{})$} denote the set of all
  maps that are $\U{}$-close to $f$. The \df{limitation
    topology}{limitation topology} on $C(X, Y)$ is defined by a basis
  of neighborhoods $\{ B(f, \U{}) \}$, where $f$ runs over all
  elements of $C(X, Y)$ and $\U{}$ runs over all open covers of $Y$.
\end{definition*}

It should be emphasized that this topology needs not to be metrizable
and the neighborhoods $B(f, \U{})$ defined above need not to be open
in it~\cite{bowers1989}.

\begin{lemma}[{\cite{torunczyk1981}}]
  \label{lem:torunczyk's lemma}
  If $X$ and $Y$ are completely metrizable, then $C(X, Y)$ endowed
  with the limitation topology has the Baire property. Moreover, every
  closed subset of $C(X, Y)$ has the Baire property.
\end{lemma}

See~\cite[p. 429]{bowers1989} for a detailed discussion.

\begin{definition*}
  We say that \df{map!approximable}{a map $f$ is approximable by
    embeddings} to express that ``for every open cover $\U{}$ of the
  codomain of $f$ there exists an embedding, which is $\U{}$-close to
  $f$''. Likewise for closed embeddings, $Z$-embeddings, etc.
\end{definition*}

In other words, a map $f$ from $X$ into $Y$ is approximable by
elements of a subset of~$C(X, Y)$ if it lies in its closure.

\begin{definition*}
  We say that a Polish space $X$ is \df{strong universality in
    dimension $n$}{strongly universal in dimension $n$} if every map
  from an at most $n$-dimensional Polish space into $X$ is
  approximable by closed embeddings.
\end{definition*}

In other words, a Polish space $Y$ is strongly universal in dimension
$n$, if for each at most $n$-dimensional Polish space $X$, the set of
closed embeddings is a dense subset of~$C(X, Y)$.

\begin{proposition}[{\cite[p. 34]{chigogidze1996}}]
  \label{pro:closed embeddings are g-delta}
  If $X$ is a Polish space, then the set of closed embeddings of $X$
  into a space $Y$ is a $G_\delta$-set in $C(X, Y)$ endowed with the
  limitation topology.
\end{proposition}

The following two theorems imply that every $n$-dimensional N\"obeling
manifold is strongly universal in dimension~$n$, which is the ``easy''
part of the characterization theorem.

\begin{theorem}
  \label{thm:nobeling space is strongly universal}
  The $n$-dimensional N\"obeling space is strongly universal in
  dimension~$n$.
\end{theorem}

Following the spirit of the literature, we omit the proof. Hint: by
lemma~1.4(b) from~\cite{torunczyk1981}, it suffices to show that every
map from a countable simplicial complex into~$\nu^n$ is approximable
by closed embeddings.

%

\begin{theorem}
  \label{thm:strong universality is a local property}
  If a Polish $ANE(n)$-space $X$ is strongly universal in dimension
  $n$, then every open subset of $X$ is strongly universal in
  dimension $n$. If every point of a Polish $ANE(n)$-space $X$ has a
  neighborhood that is strongly universal in dimension~$n$, then $X$
  is strongly universal in dimension $n$.
\end{theorem}

Theorem~\ref{thm:strong universality is a local property} appears to
be a folk theorem, well known to the experts, but unpublished.

\begin{proof}
  We base the proof on the observation that a map $g$ is a closed
  embedding whenever there is a closed interior cover $\G{}$ of the
  codomain such that for every $G \in \G{}$ the restriction
  $g_{|g^{-1}(G)}$ is a closed embedding.

  Consider the following condition for each collection~$\F{}$ of
  subsets of a space~$Y$.

  \vspace{ 1mm}
  \begin{tabular}{rl}
    ($*_{\F{}, Y}$) &
    \begin{tabular}{p{\propwidth}}\noindent\em
      for each $F \in \F{}$,  every map from at most $n$-dimensional Polish
      space into~$F$ is approximable by closed embeddings into~$Y$ (with a control 
      by covers of~$Y$).
    \end{tabular}
  \end{tabular}

  \vspace{ 1mm}
  \noindent We will show that if ($*_{\F{}, Y}$) is satisfied for some
  interior cover~$\F{}$ of an $ANE(n)$-space~$Y$, then~$Y$ is strongly
  universal in dimension $n$.

  Let $\F{}$ be an interior cover of an $ANE(n)$-space $Y$. Let $\G{}$
  denote a closed interior cover of $Y$ whose second star
  refines~$\F{}$.  Let $X$ denote an at most $n$-dimensional Polish
  space and let $f$ be a map from $X$ into $Y$. Let $B$ denote the
  interior, with respect to the limitation topology, of the set of
  maps from $X$ into $Y$ that are $\G{}$-close to $f$. It contains $f$
  (see~\cite{bowers1989}). It has the Baire property by
  lemma~\ref{lem:torunczyk's lemma}. Fix $G \in \G{}$ and let $B_G$
  denote a subset of $B$ that contains maps whose restrictions to
  $f^{-1}(\st_\G{} G)$ are closed embeddings. By
  proposition~\ref{pro:closed embeddings are g-delta}, it is a
  $G_\delta$-set in $B$.  We'll prove that it is dense in $B$.  Let
  $g$ denote an element of $B$ and let $\U{}$ be an open cover of~$Y$.
  We will show that $g$ has a $\U{}$-approximation by a map that
  belongs to~$B_G$. Since $B$ is open, we may assume that every map
  $\U{}$-close to $g$ belongs to $B$.  Let $\V{}$ denote a cover that
  satisfies condition $(*)_n$ of theorem~\ref{thm:n-homotopy extension
    theorem} for $\U{}$. Since $g$ is $\G{}$-close to $f$ and the
  second star of $\G{}$ refines $\F{}$, there is $F \in \F{}$ such
  that $g(f^{-1}(\st_\G{} G)) \subset F$.  By ($*_{\F{}, Y}$) there is
  a $\V{}$-approximation of $g_{|f^{-1}(\st_\G{} G)}$ by a closed
  embedding into $Y$ and by the definition of~$\V{}$ it extends over
  $X$ to a map that is $\U{}$-close to $g$. Therefore $B_G$ is dense
  in $B$.

  By the Baire property of $B$ the intersection $\bigcap_{G \in \G{}}
  B_G$ is dense in $B$. Therefore~$f$ is approximable by maps that are
  $\G{}$-close to $f$ and whose restrictions to $f^{-1}(\st_\G{} G)$
  are closed embeddings for every $G \in \G{}$. Let $g$ denote such a
  map.  Observe that $g^{-1}(G) \subset f^{-1}(\st_\G{} G)$, since $g$
  is $\G{}$-close to $f$. Therefore $g$ is a closed embedding on
  $g^{-1}(G)$ for every $G \in \G{}$. Hence it is a closed embedding
  of $X$ into $Y$.  Therefore we proved that if ($*_{\F{}, Y}$) is
  satisfied for an interior cover~$\F{}$ of an $ANE(n)$-space $Y$,
  then~$Y$ is strongly universal in dimension~$n$.

  To finish the proof it suffices to show that if we assume hypthesis
  of either of the implications of the theorem, then there exists an
  interior cover $\F{}$ of $X$ that satisfies condition~($*_{\F{},
    X}$).

  Let $X$ be a Polish $ANE(n)$-space that is strongly universal in
  dimension $n$. Let $U$ be an open subset of $X$. For each $x \in U$
  let $F_x$ be a closed (in $X$) neighborhood of $x$ such that $F_x
  \subset U$. We will show that the interior cover $\F{} = \{ F_x
  \}_{x \in U}$ of $U$ satisfies ($*_{\F{}, U}$). Pick $F_x \in \F{}$
  and let~$\U{}$ be an open cover of $U$. Let $\V{} = \{ X \setminus
  F_x \} \cup \U{}$. It is a cover of $X$. By strong universality
  of~$X$, every map from at most $n$-dimensional Polish space into
  $F_x$ is $\V{}$-approximable by closed embeddings into $X$. But
  every $\V{}$-approximation of a map into $F_x$ must have image in
  $U$. Hence $\F{}$ satisfies condition ($*_{\F{}, U}$). By
  theorem~\ref{thm:ane characterization}, $U$ is an $ANE(n)$-space.
  Hence $U$ is stringly universal in dimension $n$. This finishes the
  proof of the first implication.

  Let $X$ be a Polish $ANE(n)$-space such that each point $x$ in $X$
  has a neighborhood $N_x$ that is strongly universal in
  dimension~$n$. This neighborhood may be neither open nor closed. Let
  $F_x$ be a closed neighborhood of $x$ and let $U_x$ be an open
  neighborhood of $x$ such that $F_x \subset U_x \subset \Cl_X U_x
  \subset N_x$. We will show that the interior cover $\F{} = \{ F_x
  \}_{x \in X}$ of $X$ satisfies ($*_{\F{}, X}$). Pick $F_x \in \F{}$
  and let $\U{}$ be an open cover of $X$. Let $f$ be a map from at
  most $n$-dimensional Polish space into $F_x$. We have to show that
  $f$ has a $\U{}$-approximation by a closed embedding into $X$.
  Without a loss of generality we may assume that $\U{}$ refines $\{
  U_x, X \setminus F_x \}$. Let $\V{} = \{ U \cap N_x \}_{U \in
    \U{}}$. Since $N_x$ is strongly universal in dimension $n$, $f$
  has a $\V{}$-approximation by a closed embedding into $N_x$. Since
  $\V{}$ refines $\{ U_x, X \setminus F_x \}$ and the image of $f$ is
  contained in $F_x$, the image of this approximation is contained in
  $U_x$. This image is closed in $\Cl_X U_x$, as $\Cl_X U_x$ is a
  subset of $N_x$, so it must be closed in $X$. Hence ($*_{\F{}, X}$)
  is satisfied and $X$ is strongly universal in dimension $n$. This
  finishes the proof of the second implication.
\end{proof}

The characterization theorem states that every $n$-dimensional Polish
strongly universal $ANE(n)$-space is a N\"obeling manifold. It is
convenient to have a shorthand for the class of spaces that possess
these properties. To this end we adapt the notion of an $\N{n}$-space
from~\cite{kawamuralevintymchatyn1997}, with the exception that we do
not require that an $\N{n}$-space has vanishing homotopy groups of
dimensions less than $n$.

\begin{definition}\label{def:n-space}
  We say that a space is an \df{Nn-@$\N{n}$-!space}{$\N{n}$-space} if it is an
  $n$-dimensional Polish $ANE(n)$-space that is strongly universal in
  dimension~$n$. By \df{Nn@$\N{n}$}{$\N{n}$} we denote the class of
  $\N{n}$-spaces.
\end{definition}

In the language used in the introduction, the class of $\N{n}$-spaces
is the class of abstract $n$-dimensional N\"obeling manifolds.

\begin{theorem}
  \label{thm:open subsets of n-spaces}
  Every non-empty open subset of an $\N{n}$-space is an $\N{n}$-space.
  If every point of a separable space $X$ has a neighborhood that is
  an $\N{n}$-space, then $X$ is an $\N{n}$-space.
\end{theorem}
\begin{proof}
  Let $U$ be a non-empty open subset of an $\N{n}$-space $X$. By
  strong universality of $X$, one can embed an $n$-dimensional space
  into $U$, hence $U$ itself is $n$-dimensional. By \cite[theorem
  4.3.23]{engelking1989}, every open subset of a Polish space is
  Polish, hence $U$ is Polish.  By theorem~\ref{thm:ane
    characterization}, $U$ is an absolute extensor in dimension~$n$.
  By theorem~\ref{thm:strong universality is a local property}, $U$ is
  strongly universal in dimension $n$.  Hence $U$ is an $\N{n}$-space.

  If every point in $X$ has an $\N{n}$-neighborhood, then $X$ is
  clearly $n$-dimensional. Complete metrizability is a local property,
  hence $X$ is Polish. By theorem~\ref{thm:ane characterization}, $X$
  is an $ANE(n)$-space. By theorem~\ref{thm:strong universality is a
    local property}, $X$ is strongly universal in dimension $n$. Hence
  $X$ is an $\N{n}$-space.
\end{proof}

We finish this section with a technical lemma that will be used in
chapter~\ref{ch:swelling}.

\begin{lemma}
  \label{lem:limit close to the inclusion}
  For each open cover $\U{}$ of a Polish space $X$ there exists a
  sequence $\U{1}, \U{2}, \ldots$ of open covers of $X$ such that for
  each sequence $f_1, f_2, \ldots$ of maps from $X$ into $X$ such that
  $\im f_{k-1} \subset \dom f_k$ and such that $f_k$ is $\U{k}$-close
  to the inclusion $\dom f_k \subset X$, the pointwise limit $\lim_{k
    \to \infty} f_k$ exists, is continuous and is $\U{}$-close to the
  identity.
\end{lemma}
\begin{proof}
  Let $\alpha$ be a map from $X$ into $(0, \infty)$. We say that a map
  $g \colon X \to X$ is $\alpha$-close to a map $f \colon X \to X$, if
  $d(f(x),g(x)) \leq \alpha(f(x))$ for each $x \in X$. By lemma~2.1.6
  of~\cite{chigogidze1996} there is an open cover $\V{\alpha}$ of $X$
  such that if a map $g$ is $\V{\alpha}$-close to a map $f$, then it
  is $\alpha$-close to it. We define a map $\alpha_\U{} \colon X \to
  (0, \infty)$ by the formula $\alpha_\U{}(x) = (1/2) \sup \{ d(x, X
  \setminus U) \colon U \in \U{} \}$ (the formula is taken from
  page~35 of~\cite{chigogidze1996}). Again by lemma~2.1.6
  of~\cite{chigogidze1996} if $g$ is $\alpha_\U{}$-close to $g$, then
  it is $\U{}$-close to it. We let
  \[
    \U{k} = \V{ (1/2^{k+1}) \alpha_\U{}}.
  \]

  By the triangle inequality for the Hausdorff metric for each $U$ the map $x
  \mapsto d(x, X \setminus U)$ is $1$-Lipschitz. Hence $\alpha_\U{}$ is
  $1$-Lipschitz. Let $f_k \colon X \to X$, $k \geq 0$, be a sequence of maps
  such that $f_k$ is $\U{k}$-close to $f_{k-1}$. A direct computation shows
  that the distance between $f_k(x)$ and $f_0(x)$ is at most
  $(1-1/2^k)\alpha_\U{}(f_0(x))$. Therefore by completeness of $X$ the limit
  $f_\infty = \lim_{k \to \infty} f_k$ exists and the map $f_\infty$ is
  $\alpha_\U{}$-close, hence $\U{}$-close, to $f_0$.

\end{proof}

\section{$n$-Homotopy equivalence}

First we introduce a notion of a weak $n$-homotopy equivalence. Note
that by a characterization theorem of J.~H.~C.~Whitehead (to be given
on page~\pageref{thm:whitehead characterization}), the phrase
``$n$-homotopy equivalences'' in theorems stated in the introduction
can be replaced by ``weak $n$-homotopy equivalences''.

\begin{definition*}
  We say that a map is a \df{weak $n$-homotopy equivalence}{weak
  $n$-homotopy equivalence} if it induces isomorphisms on homotopy
  groups of dimensions less than~$n$, regardless of the choice of the
  base point.
\end{definition*}

\begin{definition*}
  We say that maps $f, g \colon X \to Y$ are
  \df{maps!n-homotopic@$n$-homotopic}{$n$-homotopic} if for every
  map~$\varPhi$ from a complex of dimension less than $n$ into $X$,
  the compositions $f \circ \varPhi$ and $g \circ \varPhi$ are
  homotopic in the usual sense.
\end{definition*}

\begin{lemma}\label{lem:n-homotopic to equivalence is equivalence}
  If a map is $n$-homotopic to a weak $n$-homotopy equivalence, then
  it is a weak $n$-homotopy equivalence.
\end{lemma}
\begin{proof}
  Let $\varphi \colon X \to Y$ ba a weak $n$-homotopy equivalence and
  let $\psi$ be a map that is $n$-homotopic to $\varphi$. Fix a base
  point in $X$ and let $0 \leq k < n$. We have to show that $\psi$
  induces an isomorphism on $k$th homotopy groups of $X$ and $Y$. Let
  $K$ be a wedge of $k$-dimensional spheres and let $f \colon K \to Y$
  be a map that induces an epimorphism on $k$th homotopy groups of $K$
  and $Y$. By the definition of $n$-homotopy, $\varphi \circ f$ and
  $\psi \circ f$ are homotopic. Let $\gamma$ be a patch of the
  base-point under this homotopy. Let $\varphi_*$ be an isomorphism of
  $\pi_k(X)$ and $\pi_k(Y)$ induced by $\varphi$.  The homomorphism
  $\psi_*$ induced by $\psi$ on $\pi_k(X)$ and $\pi_k(Y)$ is equal to
  an isomorphism $\varphi_*$ conjugated with $\gamma$ (see
  \cite[p.~341]{hatcher2002} for precise definitions). Hence $\psi_*$
  is an isomorphism as well.
\end{proof}

\begin{lemma}
  \label{lem:set of weak equivalences is open}
  Every $ANE(n)$-space $Y$ admits an open cover $\U{}$ such that every
  map that is $\U{}$-close to a weak $n$-homotopy equivalence of a
  metric space $X$ and $Y$ is also a weak $n$-homotopy equivalence. 
\end{lemma}

In particular, the set of weak $n$-homotopy equivalences of $X$ and
$Y$ is an open subset of $C(X, Y)$.

\begin{proof}
  By theorem~\ref{thm:ane characterization}, there exists an open
  cover $\U{}$ of $Y$ such that every two $\U{}$-close maps defined on
  a metric space of dimension less than $n$ are homotopic. In
  particular, every two $\U{}$-close maps into $Y$ are $n$-homotopic.
  By lemma~\ref{lem:n-homotopic to equivalence is equivalence}, the
  cover $\U{}$ satisfies our claim.
\end{proof}

\begin{lemma}\label{lem:lifting under n-homotopy}
  If $\varphi \colon X \to Y$ is a weak $n$-homotopy equivalence, $K$
  is a simplicial complex whose dimension is less than $n$ and $f$ is
  a map from $K$ into $Y$, then there exists a map $f' \colon K \to X$
  whose composition with $\varphi$ is homotopic to $f$. If $L$ is a
  subcomplex of $K$ and if we are given a map from $L$ into $X$ whose
  composition with $\varphi$ is equal to $f_{|L}$, then we may require
  that $f'$ extends this map and that the homotopy of $f$ and $\varphi
  \circ f'$ is $\rel L$.  If $K$ is a subcomplex of an at most
  $n$-dimensional simplicial complex $M$ and $f$ extends to a map $g$
  from $M$ into $Y$, then $f'$ extends to a map from $M$ into $X$. In
  particular, there exists a map from $M$ into $X$ whose composition
  with~$\varphi$ is $n$-homotopic to $g$.
\end{lemma}
\begin{proof}
  The proof of the first two implications follows exactly the proof of
  proposition~6.4 from page~77 of~\cite{bessagapelczynski1975}. To
  prove the third implication, consider an extension $g \colon M \to
  X$ of $f$. Without a loss of generality we may assume that there are
  no simplices of dimension less than $n$ in $M \setminus K$ (by
  enlarging $K$ if needed) and that $M$ is endowed with a weak
  topology (by theorem~6.4 on page~81 of~\cite{bessagapelczynski1975},
  $M$ endowed with a weak topology is homotopy equivalent to $M$
  endowed with a metric topology). Hence all we have to do is to show
  how to extend $f'$ over a single $n$-dimensional simplex~$\delta$
  of~$M$. The map $f'$ is defined on $\partial \delta$. By the
  assumptions, $\varphi \circ f'_{|\partial \delta}$ is homotopic to
  $f_{|\partial \delta}$ and $f_{|\partial \delta}$ extends over
  $\delta$ to a map into $X$. Hence $\varphi \circ f'_{|\partial
    \delta}$ represents a trivial element of $\pi_{n-1}(Y)$. Since
  $\varphi$ is a weak $n$-homotopy equivalence, $f'_{|\partial
    \delta}$ must represent a trivial element of $\pi_{n-1}(X)$.
  Hence $f'_{|\partial \delta}$ extends over $\delta$ to a map into
  $X$.  Hence $f'$ admits an extension to a map $g'$ from $M$ into
  $X$.  It follows from the cellular approximation
  theorem~\cite[Corollary~4.12]{hatcher2002}, that $\varphi \circ g'$
  is $n$-homotopic to $g$.
\end{proof}

\begin{definition*}
  We say that a map $\varphi \colon X \to Y$ is an \df{n-homotopy
    equivalence@$n$-homotopy equivalence}{$n$-homotopy equivalence} if there
  exists a map $\psi \colon Y \to X$ such that $\psi \circ \varphi$ is
  $n$-homotopic to $id_X$ and $\varphi \circ \psi$ is $n$-homotopic to $id_Y$.
  We say that $X$ is \df{n-homotopy equivalent@$n$-homotopy
    equivalent}{$n$-homotopy equivalent} to $Y$ if there exists an
  $n$-homotopy equivalence of $X$ and $Y$.
\end{definition*}

\begin{whitehead's characterization}
  \label{thm:whitehead characterization}\index{Whitehead's characterization}
  A map of two $n$-dimensional $ANE(n)$-spaces is a weak $n$-homotopy
  equivalence if and only if it is an $n$-homotopy equivalence.
\end{whitehead's characterization}

The theorem was stated by Whitehead~\cite{whitehead1949} with
slightly stronger assumptions, although his proof works in this case
as well. For completeness, we give a proof based on the argument
from~\cite{bessagapelczynski1975}. A proof of an even more general
case appeared in~\cite{chigogidzekarasev2003}.

\begin{proof}
  The proof is a modification of the proof of theorem~6.1 from page~78
  of~\cite{bessagapelczynski1975}. Let $\varphi$ be a map from an at
  most $n$-dimensional $ANE(n)$-space $X$ into an at most
  $n$-dimensional $ANE(n)$-space $Y$. We shall write $f \nhom g$ to
  denote that $f$ and $g$ are $n$-homotopic.

  Assume that $\varphi$ is an $n$-homotopy equivalence. By the
  definition, there exists a map $\psi \colon Y \to X$ such that $\psi
  \circ \varphi \nhom id_X$ and $\varphi \circ \psi \nhom id_Y$. The
  identity on $X$ is a weak $n$-homotopy equivalence, hence by
  lemma~\ref{lem:n-homotopic to equivalence is equivalence}, $\psi
  \circ \varphi$ is a weak $n$-homotopy equivalence. Hence for each $k
  < n$, $\psi \circ \varphi$ induces an isomorphism on $\pi_k(X)$, so
  $\varphi$ induces a monomorphism of $\pi_k(X)$ and $\pi_k(Y)$. Again
  by lemma~\ref{lem:n-homotopic to equivalence is equivalence},
  $\varphi \circ \psi$ induces an isomorphism on $\pi_k(Y)$, hence
  $\varphi$ induces an epimorphism of $\pi_k(X)$ and $\pi_k(Y)$.
  Hence $\varphi$ is a weak $n$-homotopy equivalence.

  Assume that $\varphi$ is a weak $n$-homotopy equivalence.  By
  theorem~4.1 on page 156 in~\cite{hu1965} and by theorem~\ref{thm:ane
    characterization}, every $n$-dimensional $ANE(n)$-space is
  $n$-homotopy dominated by an at most $n$-dimensional complex. Hence
  we get diagrams $K \stackrel{g_K}{\longrightarrow} X
  \stackrel{f_K}{\longrightarrow} K$ and $L
  \stackrel{g_L}{\longrightarrow} Y \stackrel{f_L}{\longrightarrow} L$
  where $K$, $L$ are at most $n$-dimensional complexes, $g_K \circ f_K
  \nhom id_X$ and $g_L \circ f_L \nhom id_Y$. By
  lemma~\ref{lem:lifting under n-homotopy} there is a map $g'_L \colon
  L \to X$ such that $\varphi \circ g'_L \nhom g_L$.  Let $\psi = g'_L
  \circ f_L$. We'll show that $\varphi$ and $\psi$ are $n$-homotopy
  inverses of each other. Observe that $\varphi \circ \psi = \varphi
  \circ g'_L \circ f_L \nhom f_L \circ g_L \nhom id_Y$.  By
  lemma~\ref{lem:n-homotopic to equivalence is equivalence}, $\psi$ is
  a weak $n$-homotopy equivalence. Hence by lemma~\ref{lem:lifting
    under n-homotopy}, there is a map $g'_K \colon K \to Y$ such that
  $\psi \circ g'_K \nhom g_K$.  Therefore $\varphi \nhom \varphi \circ
  g_K \circ f_K \nhom \varphi \circ \psi \circ g'_K \circ f_K \nhom
  g'_K \circ f_K$.  So $\psi \circ \varphi \nhom \psi \circ g'_K \circ
  f_K \nhom g_K \circ f_K \nhom id_X$ and $\varphi$ is an $n$-homotopy
  equivalence.
\end{proof}

\begin{theorem}
  \label{thm:subspace retraction} If $X$ is an at most $n$-dimensional
  $ANE(n)$-space and $A$ is a closed $ANE(n)$-subset of $X$ such that
  the inclusion $A \subset X$ is a weak $n$-homotopy equivalence, then
  $A$ is a retract of $X$. Moreover, every retraction of $X$ onto $A$
  is a weak $n$-homotopy equivalence.
\end{theorem}
\begin{proof}\footnote{I'm grateful to A.~Chigogidze for showing me this simple proof.}
  By Whitehead's characterization, the inclusion $\imath \colon A
  \subset X$ is an $n$-homotopy equivalence. By the definition of an
  $n$-homotopy equivalence, there exists a map $\psi \colon X \to A$
  whose restriction to $A$ is $n$-homotopic to the identity on $A$. By
  proposition~4.1.8 from~\cite{chigogidze1996}, the identity on $A$
  extends over $X$ to a retraction.
  
  Let $r$ be any retraction of $X$ onto $A$. The composition $r \circ
  \imath$ is equal to $\imath$, which is a weak $n$-homotopy
  equivalence by the assumptions. Hence $r$ must be a weak
  $n$-homotopy equivalence as well.
\end{proof}

\begin{lemma}
  \label{lem:excision}
  If $A_1$, $A_2$ and $A_1 \cap A_2$ are at most $n$-dimensional closed
  $ANE(n)$-subspaces of $A_1 \cup A_2$ and the inclusion $A_1 \cap A_2 \subset
  A_2$ is a weak $n$-homotopy equivalence, then the inclusion $A_1 \subset A_1
  \cup A_2$ is a weak $n$-homotopy equivalence.
\end{lemma}
\begin{proof}
  By theorem~\ref{thm:subspace retraction}, there exists a retraction
  of $A_2$ onto $A_1 \cap A_2$. It extends to a retraction $r \colon
  A_1 \cup A_2 \to A_1$, since both $A_1$ and $A_2$ are closed in $A_1
  \cup A_2$. Therefore the inclusion $A_1 \subset A_1 \cup A_2$
  induces monomorphisms on all homotopy groups.
  
  By theorem~\ref{thm:ane characterization}, there exists an open
  cover $\U{}$ of $A_1 \cup A_2$ such that every two $\U{}$-close maps
  from an at most $(n-1)$-dimensional Polish space are homotopic.
  Let~$\V{}$ be an open cover of $A_1 \cup A_2$ that satisfies
  condition $(*)_n$ of theorem~\ref{thm:n-homotopy extension theorem}.
  Like we argued in the proof of lemma~\ref{lem:extension property},
  the restriction of $r$ to a small closed neigbhorhood of $A_1$ in
  $A_1 \cup A_2$ is $\V{}$-close to the identity. By the choice of
  $\V{}$, this restriction admits an extension to a map $s$ from $A_1
  \cup A_2$ into $A_1 \cup A_2$ that is $\U{}$-close to the identity.

  Fix $k < n$ and consider a map $\varphi \colon S^k \to A_1 \cup
  A_2$. Let $\psi = s \circ \varphi$. By the choice of $\U{}$, $\psi$
  is homotopic to $\varphi$. By the construction of $s$, there exist
  open subsets $U, V \subset S^k$ such that $U \subset
  \psi^{-1}(A_1)$, $V \subset \psi^{-1}(A_2)$ and $U \cup V = S^k$.
  Hence we may fix subcomplexes $C$ and $D$ of $S^k$ such that $C \cup
  D = S^k$ and $\psi(C) \subset A_1$ and $\psi(D) \subset A_2$. By
  Whitehead's characterization, the inclusion $A_1 \cap A_2 \subset
  A_2$ is an $n$-homotopy equivalence. By lemma~\ref{lem:lifting under
    n-homotopy}, the map $\psi_{|D}$ is homotopic $\rel C \cap D$ with
  a map into $A_1 \cap A_2$. Obviously this homotopy extends to a
  homotopy $\rel C$.  Therefore the inclusion $A_1 \subset A_1 \cup
  A_2$ induces epimorphisms on homotopy groups of dimensions less than
  $n$.  We are done.
\end{proof}

\begin{corollary}
  \label{cor:excision}
  If $A_1$, $A_2$ and $A_1 \cap A_2$ are at most $n$-dimensional closed
  $ANE(n)$-subspaces of a space $X$ and inclusions $A_1 \subset X$ and $A_1
  \cap A_2 \subset A_2$ are weak $n$-homotopy equivalences, then the inclusion
  $A_1 \cup A_2 \subset X$ is a weak $n$-homotopy equivalence.
\end{corollary}
\begin{proof}
  By lemma~\ref{lem:excision}, the inclusion $A_1 \subset A_1 \cup A_2$ is a
  weak $n$-homotopy equivalence. So, by Whitehead's characterization, there
  exists a weak $n$-homotopy equivalence $A_1 \cup A_2 \to A_1$. We are done,
  by the assumption that $A_1 \subset X$ is a weak $n$-homotopy equivalence.
\end{proof}

\begin{lemma}
  \label{lem:equivalence with a complex} If $X$ is a separable
  $ANE(n)$-space and $f$ is a map of an at most $n$-dimensional
  locally finite countable simplicial complex $K$ into $X$, then there
  exists an at most $n$-dimensional locally finite countable
  simplicial complex $L$, containing $K$ as a subcomplex, and an
  extension of $f$ to a weak $n$-homotopy equivalence of~$L$ and~$X$.
\end{lemma}
\begin{proof}
  This is a standard argument that follows the argument given
  in~\cite{hatcherCWapprox}. The only difference is that we want $L$
  to be a locally finite countable simplicial complex. Let $K_0 = K$
  and let $f_0 = f$. For $m = 1, 2, \ldots, n$ we shall construct a
  locally finite countable simplicial complex $K_m$ that contains
  $K_{m-1}$ as a subcomplex and a map $f_m \colon K_m \to X$ that is
  an extension of $f_{m-1}$ to a weak $m$-homotopy equivalence (note
  that we allow the dimension of $K_m$ to be greater than $n$).

  Fix $0 \leq m < n$ and assume that we already constructed $f_{m}
  \colon K_{m} \to X$. By the separability of $X$, the space of maps
  from $S^{m}$ into $X$ is separable. By theorem~\ref{thm:ane
    characterization}, the $m$-dimensional homotopy group $\pi_m(X)$
  is countable. Hence we may enlarge $K_{m}$ to a locally finite
  countable simplicial complex $K'_{m}$ and extend $f_{m}$ to a map
  $f'_{m} \colon K'_{m} \to X$ that is a weak $m$-homotopy
  equivalence and induces an epimorphism $\pi_{m}(f'_{m})$ on $m$th
  homotopy groups. Let $\{ \varphi_k \colon S^{m} \to K'_{m} \}_{k
    \geq 1}$ be a set of representants of the kernel of
  $\pi_{m}(f'_{m})$. By the simplicial approximation
  theorem~\cite[Theorem 2C.1]{hatcher2002}, we may assume that
  $\varphi_k$'s are simplicial. Let $K''_{m} = K'_{m} \times [1,
  \infty)$, let $p \colon K''_{m} \to K'_{m}$ be a projection along
  $[1, \infty)$ and let $i_k \colon K'_{m} \to K'_m \times \{ k \}$ be
  a map such that $p \circ i_k$ is the identity for each $k \geq 1$.
  We identify $K'_{m} \times \{ 1 \}$ with $K'_{m}$ so that $K_{m}$ is
  a subcomplex of $K''_{m}$.  Observe that $f''_{m} = f'_{m} \circ p$
  is an extension of $f'_{m}$ over $K''_{m}$ that is a weak
  $m$-homotopy equivalence that induces an epimorphism on $m$th
  homotopy groups.  The set $\{ i_k \circ \varphi_k \}_{k \geq 1}$ is
  a set of representants of the kernel of $\pi_{m}(f''_{m})$. Let
  $\varPhi_k \colon B^{m+1} \to K''_{m}$ be a partial map defined on
  the boundary $S^{m}$ of $B^{m+1}$ to be equal to $\varphi_k$. Let
  $K_{m+1}$ be a mapping cylinder of an disjoint union of
  $\varPhi_k$'s. It is a locally finite countable simplicial complex.
  Since each $\varphi_k$ represents an element of the kernel of
  $\pi_{m}(f''_{m})$, map $f''_{m}$ extends to a map $f_{m+1} \colon
  K_{m+1} \to X$. This map is a weak $(m+1)$-homotopy
  equivalence~\cite{hatcherCWapprox}.

  By the cellular approximation
  theorem~\cite[Corollary~4.12]{hatcher2002}, the restriction of $f_n$
  to the $n$-dimensional skeleton of $K_n$ is a weak $n$-homotopy
  equivalence. It is a map that we were looking for.
\end{proof}

Now we are ready to prove that every $n$-dimensional absolute
neighborhood extensor in dimension $n$ is $n$-homotopy equivalent to
an open subset of $\nu^n$. This theorem together with characterization
and rigidity theorems implies the open embedding theorem for
N\"obeling manifolds.

\begin{proposition}
  \label{pro:intersection of nobeling space with an open set}
  For every open subset $U$ of $\mathbb{R}^{2n+1}$ the inclusion of $U
  \cap \nu^n$ into $U$ is a weak $n$-homotopy equivalence.
\end{proposition}
\begin{proof}
  The following statement can be proved by the standard ``general position''
  argument.

  \vspace{ 1mm}
  \begin{tabular}{rl}
    (*) &
    \begin{tabular}{p{\propwidth}}\noindent\emph{Every map from an at
        most $n$-dimensional finite simplicial complex into
        $\mathbb{R}^{2n+1}$ is approximable by maps into
        $\mathbb{R}^{2n+1} \setminus H$, where~$H$ is an arbitrary
        $n$-dimensional hyperplane. Moreover, if the map restricted to
        a subcomplex $L$ has its image disjoint from~$H$, then such
        approximations exist $\rel L$.  }\end{tabular}
  \end{tabular}
  \vspace{ 1mm}

  Let $K$ be a finite simplicial complex, let $L$ be a subcomplex of
  $K$ and let $\psi$ be a map from $L$ into $\nu^n$. Let $C_\psi$ be a
  set of maps from $K$ into $\mathbb{R}^{2n+1}$ that are equal to
  $\psi$ on $L$. It is a closed subset of $C(X, Y)$. By (*), a subset
  of $C_\psi$ consisting of maps into the complement of an
  $n$-dimensional hyperplane that is disjoint from the image of $\psi$
  is dense in $C_\psi$. By lemma~\ref{lem:torunczyk's lemma}, $C_\psi$
  has the Baire property, hence the set of maps from $K$ into $\nu^n$
  that are equal to $\psi$ on $L$ is dense in $C_\psi$.

  Let $\varphi$ be a map from a $k$-dimensional sphere into $U$,
  for~$k \leq n$. By letting $K = S^k$ and $L = \emptyset$ we see that
  $\varphi$ can be approximated by a map into $U \cap \nu^n$ and by
  theorem~\ref{thm:ane characterization}, if the approximation is
  close enough, then it is homotopic to $\varphi$ in $U$. Therefore
  the inclusion induces epimorphisms on $k$th homotopy groups, for $k
  \leq n$.

  Let $\varphi$ be a map from a $k$-dimensional sphere into $U \cap
  \nu^n$ and assume that $\varphi$ extends over a $(k+1)$-dimensional
  ball to a map $\varPhi$ into $U$. By letting $K = B^{k+1}$ and $L =
  S^k$, we see that $\varPhi$ can be approximated $\rel S^k$ by a map
  into $U \cap \nu^n$. Therefore the inclusion induces monomorphisms
  on $k$th homotopy groups, for $k < n$.  We are done.
\end{proof}

\begin{corollary}
  \label{cor:nobeling space is ae(n)}
  The $n$-dimensional N\"obeling space is an absolute extensor in
  dimension $n$.
\end{corollary}
\begin{proof}
  Every point of $\nu^n$ has arbitrarily small neighborhoods
  homeomorphic to~$\nu^n$. By proposition~\ref{pro:intersection of
    nobeling space with an open set}, the N\"obeling space~$\nu^n$ has
  vanishing homotopy groups of dimensions less than $n$. Therefore
  theorems~\ref{thm:dugundji characterization} and~\ref{thm:ane
    characterization} imply the assertion.
\end{proof}

\begin{proposition}
  \label{pro:equivalence with an open subset}
  For every separable $ANE(n)$-space $X$ there exists an open subset $U$ of
  $\nu^n$ and a weak $n$-homotopy equivalence $U \to X$.
\end{proposition}
\begin{proof}
  By lemma~\ref{lem:equivalence with a complex}, for every separable
  $ANE(n)$-space $X$ there exists an $n$-dimensional countable locally
  finite simplicial complex $K$ and a weak $n$-homotopy equivalence $K
  \to X$. Realize $K$ in $\mathbb{R}^{2n+1}$ in such a way that every
  ball around origin of the space intersects only finitely many
  vertices of $K$. With such realization it is clear that $K$ has an
  open neighborhood $U$ (in $\mathbb{R}^{2n+1}$) that is homotopy
  equivalent to $K$. By proposition~\ref{pro:intersection of nobeling
    space with an open set}, the set $U \cap \nu^n$ is weak
  $n$-homotopy equivalent to~$U$. The proof is finished, since $U$ is
  homotopy equivalent to $K$ and $K$ is weak $n$-homotopy equivalent
  to $X$.
\end{proof}

\section{$Z$-sets}

The definition of a $Z$-set was given by Anderson to describe a
general concept of a negligible set~\cite{anderson1967}. There are
many flavors of the definition to be found in the literature and most
of them are equivalent in the class of N\"obeling manifolds~\cite[p.
191]{chigogidze1996}. We are going to settle for the following one,
proposed by Toru\'nczyk. 

\begin{definition*}
  We say that a subset~$A$ of a space~$X$ is a
  \df{Z-@$Z$-!set}{$Z$-set} in~$X$ if it is closed in $X$ and if the
  identity on $X$ is approximable by maps into $X \setminus A$. A
  countable union of $Z$-sets is said to be a
  \df{Zs-set@$Z_\sigma$-set}{$Z_\sigma$-set}. A map is said to be a
  \df{Z-@$Z$-!embedding}{$Z$-embedding} if it is an embedding and its
  image is a $Z$-set in the codomain.
\end{definition*}

\begin{theorem}[{\cite[Proposition~5.1.2]{chigogidze1996}}]\label{thm:Z characterization}
  A closed subset $A$ of an $\N{n}$-space $X$ is a $Z$-set if and only
  if every map from $I^n$ into $X$ is approximable by maps into the
  complement of $A$.
\end{theorem}

\begin{proposition}[{\cite[Corollary~5.1.6]{chigogidze1996}}]
  \label{pro:compact is Z}
  Every compact subset of an $\N{n}$-space is a $Z$-set.
\end{proposition}

\begin{proposition}
  \label{pro:closed Zsigma is Z} If a closed subset of an
  $\N{n}$-space~$X$ is a countable union of $Z$-sets in $X$, then it
  is a $Z$-set in $X$.
\end{proposition}
\begin{proof}
  Let $A = \bigcup_{n \geq 1} A_n$ be a closed countable union of
  $Z$-sets in an $\N{n}$-space $X$. Let $G_n = \{ f \in C(I^n, X)
  \colon \im f \cap A_n = \emptyset \}$. By theorem~\ref{thm:Z
    characterization}, $G_n$ is a dense subset of $C(I^n, X)$. By the
  assumption that $A_n$ is closed, $G_n$ is an open subset of $C(I^n,
  X)$. By lemma~\ref{lem:torunczyk's lemma}, the intersection
  $\bigcap_{n \geq 1} G_n$ is a dense subset of $C(I^n, X)$. By
  theorem~\ref{thm:Z characterization}, $A$ is a $Z$-set in $X$.
\end{proof}

\begin{proposition}
  \label{pro:Z is a local property} Let $A$ be a closed subset of an
  $\N{n}$-space $X$. If $A$ is a $Z$-set, then for each open subset
  $U$ of $X$, $A \cap U$ is a $Z$-set in~$U$. If each point of $A$ has
  a neighborhood $U$ in $X$ such that $A \cap U$ is a $Z$-set in $U$,
  then $A$ is a $Z$-set in $X$.
\end{proposition}
\begin{proof}
  Let $A$ be a $Z$-set in $X$.  Every open subset $U$ of $X$ is a
  countable union $\bigcup_{n \geq 1} F_n$ of closed subsets of $X$.
  Directly from the definition, the intersection $F_n \cap A$ is a
  $Z$-set in~$X$. Hence the set $C(I^n, X \setminus (F_n \cap A))$ is
  a dense subset of $C(I^n, X)$. The set $C(I^n, U)$ is an open subset
  of $C(I^n, X)$.  Hence $C(I^n, U \setminus (F_n \cap A))$ is a dense
  subset of $C(I^n, U)$.  By theorem~\ref{thm:Z characterization},
  $F_n$ is a $Z$-set in $U$.  By proposition~\ref{pro:closed Zsigma is
    Z}, $A = \bigcup_{n \geq 1} F_n$ is a $Z$-set in~$U$.

  Let $A$ be a closed subset of $X$ such that each point $x$ of $A$
  has a neighborhood $U_x$ such that $A \cap U_x$ is a $Z$-set in
  $U_x$. Every open subset $U$ of $X$ admits an open cover $\U{}$ such
  that if a map $U \to U$ is $\U{}$-close to the identity on $U$, then
  it admits a countinuous extension by the identity on $X \setminus U$
  (we will elaborate on this in the proof of lemma~\ref{lem:local
    approximation}). Hence it follows from the definition that if a
  closed subset of $X$ is a $Z$-set in an open subset $U$ of $X$, then
  it is a $Z$-set in $X$. Since $X$ is separable metric, there exists
  a countable cover $\U{}$ of $X$ that is a subset of $\{ U_x \}_{x
    \in X}$. Let $\F{}$ be a closed shrinking of $\U{}$. By
  proposition~\ref{pro:closed Zsigma is Z}, $A = \bigcup_{F \in \F{}}
  A \cap F$ is a $Z$-set in $X$.
\end{proof}

\begin{remark}\label{rem:closed locally compact is Z}
  Propositions~\ref{pro:compact is Z} and~\ref{pro:Z is a local
    property} imply that every closed locally compact subset of an
  $\N{n}$-space is a $Z$-set.
\end{remark}

\begin{theorem}[{\cite[Proposition~5.1.7]{chigogidze1996}}]
  \label{thm:z-approximation}
  If $X$ is an $\N{n}$-space, then every map $f$ from at most $n$-dimensional
  Polish space into $X$ is approximable by $Z$-embeddings. Moreover, if
  $f$ restricted to a closed subset~$A$ of its domain is a closed
  $Z$-embedding, then such approximations exist $\rel A$.
\end{theorem}

We'll need the following corollaries of
theorem~\ref{thm:z-approximation}.

\begin{corollary}
  \label{cor:inverse equivalence} Let $h$ denote a homeomorphism of a
  $Z$-subset of an $\N{n}$-space $X$ and a $Z$-subset of an $\N{n}$-space $Y$.
  If $h$ extends to an $n$-homotopy equivalence of $X$ and $Y$, then $h^{-1}$
  extends to an $n$-homotopy equivalence of $Y$ and $X$. Moreover, we may
  require the extension of $h^{-1}$ to be a closed embedding.
\end{corollary}
\begin{proof}
  By theorem~\ref{thm:z-approximation} the $n$-homotopy equivalence of $X$
  and $Y$ can be approximated $\rel Z_1$ by a closed embedding $h_1 \colon X
  \to Y$. By theorem~\ref{thm:ane characterization} a sufficiently close
  approximation of an $n$-homotopy equivalence is an $n$-homotopy equivalence,
  therefore we may assume that $h_1$ is an $n$-homotopy equivalence. By
  theorem~\ref{thm:subspace retraction} there exists a retraction $r_1 \colon
  Y \to h_1(X)$. The map $h_1^{-1} \circ r_1$ is an extension of $h^{-1}$ to
  an $n$-homotopy equivalence of $Y$ and $X$.
\end{proof}

\begin{corollary}
  \label{cor:z-approximation}
  If $D$ is a $Z_\sigma$-set in an $\N{n}$-space $Y$, then every map $f$ from
  an $n$-dimensional Polish space into $Y$ is approximable by closed
  embeddings with images disjoint from~$D$. If $f$ restricted to a closed
  subset $A$ of its domain is a $Z$-embedding with image disjoint from
  $D$, then such approximations exist $\rel A$.
\end{corollary}
\begin{proof}
  Let $C(X, Y)$ denote the space of maps from $X$ into $Y$ endowed
  with the limitation topology. It is easy to verify that the set $G$
  of $Z$-embeddings is a $G_\delta$ in $C(X, Y)$.
  Let $D = \bigcup Z_n$ and assume that every
  $Z_n$ is a $Z$-set in $Y$. Let $G_n = \{ f \in G \colon f(X) \cap
  Z_n = \emptyset \}$. Let $F$ denote the set of maps equal to $f$ on
  $A$. By theorem~\ref{thm:z-approximation} the set $F \cap G_n$ is
  dense in $F$.  Therefore by Toru\'nczyk's lemma, the set $F \cap
  \bigcap G_n$ is dense in $F$ and we are done.
\end{proof}

\begin{lemma}\label{lem:local approximation}
  Every open subset $U$ of a metric space $Y$ admits an open cover
  $\U{}$ such that if $V$ is an open neighborhood of $Y \setminus U$,
  $f$ is a map from a metric space into $Y$ whose restriction to
  $f^{-1}(V)$ is a (closed) embedding into $V$ and a map $g$ is a
  $\U{}$-approximation of the restriction of $f$ to $f^{-1}(U)$ by a
  (closed) embedding into $U$, then $f_{|f^{-1}(Y \setminus U)} \cup
  g$ is a (closed) embedding into $Y$.
\end{lemma}
\begin{proof}
  Let $A$ be the complement of $U$ in $Y$. Let $\U{}$ be a cover (we
  allow it to be uncountable) of $U$ by a collection of balls such
  that each ball from~$\U{}$ has a radius equal to one-fifth of the
  distance from its center to $A$.  Let $V$ be an open neighborhood of
  $A$ and let $f$ be a map from a metric space $X$ into $Y$ whose
  restriction to $f^{-1}(V)$ is an embedding into~$V$.  Let $g \colon
  f^{-1}(U) \to U$ be a $\U{}$-approximation of the restriction of~$f$
  to~$f^{-1}(U)$ by an embedding into $U$.  We shall prove that $h =
  f_{|f^{-1}(A)} \cup g$ is an embedding of $X$ into $Y$.

  The image of $f_{|f^{-1}(A)}$ is a subset of $A$, the image of $g$
  lies in the complement of $A$ and both of these maps are one-to-one,
  hence $h$ is one-to-one. Let $x_0, x_1, \ldots$ be a sequence
  in~$X$. If $h(x_0) \in U$, then $x_n \to x_0$ iff $g(x_n) \to
  g(x_0)$ (by which we understand that $g(x_n)$ is defined for almost
  all $n$ and converges to $g(x_0)$) iff $h(x_n) \to h(x_0)$ (because
  $U$ is open).  If $h(x_0) \in A$, then
  \[\tag{*}
  \frac{2}{3} d(f(x_n), f(x_0)) \leq d(h(x_n), h(x_0)) \leq
  \frac{3}{2} d(f(x_n), f(x_0)),
  \]
  by the choice of $\U{}$. Hence if $h(x_0) \in A$, then $x_n \to x_0$
  iff $f(x_n) \to f(x_0)$ (because $f$ is an embedding on $f^{-1}(V)$)
  iff $h(x_n) \to h(x_0)$ (by (*)). Hence
  \[
  \lim_{n \to \infty} x_n = x_0 \Leftrightarrow \lim_{n \to \infty} h(x_n) =
  h(x_0),
  \]
  so $h$ is a homeomorphism.

  Assume that the restriction of $f$ to $f^{-1}(V)$ is a closed
  embedding into $V$ and that $g$ is a closed embedding into $U$. To
  show that $h$ is a closed embedding it suffices to show that $h(X)$
  is closed in $Y$. If $y \in A$, then $h(x_n) \to y$ iff $f(x_n) \to
  y$, by (*). Hence $A \cap \Cl_X h(X) = A \cap \Cl_X f(X) = A \cap
  f(X) = A \cap h(X)$. If $y \in U$, then $h(x_n) \to y$ iff $g(x_n)
  \to y$. Hence $U \cap \Cl_X h(X) = U \cap \Cl_X g(X) = U \cap g(X) =
  U \cap h(X)$. Therefore $h(X)$ is equal to its closure. We are
  done.
\end{proof}

\begin{proposition}
  \label{pro:approximation rel inverse image}
  Assume that $V$ is an open subset of an $\N{n}$-space $Y$ and $f$ is
  a map from an at most $n$-dimensional Polish space $X$ into $Y$. If
  $f$ is a closed embedding into $V$ on $f^{-1}(V)$ and $A$ is a
  closed (in $Y$) subset of $V$, then $f$ is approximable $\rel
  f^{-1}(A)$ by closed embeddings. Moreover, if $Z$ is a subset of $X$
  and $f_{|Z}$ is a $Z$-embedding, then such approximations exist
  $\rel f^{-1}(A) \cup Z$. Even more, if $Z_Y$ is a $Z_\sigma$-set in
  $Y$ and $Z_Y$ is disjoint from $f(Z)$, then we may additionally
  require that the approximations have images disjoint from $Z_Y$.
\end{proposition}
\begin{proof}
  Let $\V{}$ be an open cover of $Y$. Let $U = Y \setminus A$ and let
  $\U{}$ be an open cover of $U$ obtained via lemma~\ref{lem:local
    approximation} applied to $U$ and $Y$. Without a loss of
  generality, we may assume that $\U{}$ refines $\V{}$. By
  proposition~\ref{pro:Z is a local property}, the restriction of $f$
  to $Z \cap f^{-1}(U)$ is a $Z$-embedding into $U$. By
  theorem~\ref{thm:strong universality is a local property}, $U$ is
  strongly universal in dimension~$n$, hence by
  corollary~\ref{cor:z-approximation} applied with $D = Z_Y \cap U$,
  there exists a $\U{}$-approximation $g$ of $f_{|f^{-1}(U)}$ $\rel Z
  \cap f^{-1}(U)$ by a closed embedding into $U$. By
  lemma~\ref{lem:local approximation} and by the choice of~$\U{}$, the
  map $g \cup f_{|f^{-1}(U)}$ is a closed embedding into $Y$ and is
  $\V{}$-close to $f$.
\end{proof}

\begin{definition*}
  A \df{Nn-@$\N{n}$-!neighborhood}{closed $\N{n}$-neighborhood} of a set $A$
  is a closed set that is an $\N{n}$-space and contains $A$ in its interior.
\end{definition*}

\begin{remark}\label{rem:small n-neighborhoods}
  Proposition~\ref{pro:approximation rel inverse image} and
  theorem~\ref{thm:strong universality is a local property} imply that
  every closed subset of an $\N{n}$-space has arbitrarily small closed
  $\N{n}$-neighborhoods.
\end{remark}


\newpage\thispagestyle{empty}

\part{Reducing the proof of the main results to the con\-stru\-ction of
  $n$-regular and $n$-semiregular $\N{n}$-covers}

\chapter{Approximation within an $\N{n}$-cover}
\label{ch:approximation within a cover}

The aim of this section is to prove theorem~\ref{thm:approximation
  within a cover} about approximation within closed $\N{n}$-covers.
Its proof relies on proposition~\ref{pro:compact is Z}, which states
that every compact subset of an $\N{n}$-space is a $Z$-set. This
property and the theorem fail in the case of Hilbert cube
manifolds~\cite{sher1977}.  Corollary~\ref{cor:pasting} is a sum
theorem for abstract N\"obeling manifolds.

Recall that $\N{n}$ is the class of abstract $n$-dimensional
N\"obeling manifolds (definition~\ref{def:n-space}). We define
$\N{n}$-covers in the following way.

\begin{definition*}
  Let~$\mathcal{C}$ be a class of topological spaces. We say that a
  collection is a
  \df{C-@$\protect\mathcal{C}$-collection}{$\mathcal{C}$-collection}
  if the intersection of each non-empty collection of its elements
  either belongs to~$\mathcal{C}$ or is empty. If the collection
  covers the underlying space, then we say that it is a
  \df{cover!C-@$\protect\mathcal{C}$-}{$\mathcal{C}$-cover}.
\end{definition*}

\section{$Z(\F{})$-sets}

In order to effectively work with closed $\N{n}$-covers we need the
following basic generalization of the concept of a $Z$-set.

\begin{definition*}
  Let $\F{}$ be a collection of subsets of a space $X$ and let~$A$ be
  a subset of~$X$. We say that $A$ is a \df{ZF-@$Z(\F{})$-!set in
    $X$}{$Z(\F{})$-set in $X$} if it is a $Z$-set in~$X$ and if for
  each subset $\E{}$ of $\F{}$ the intersection $A \cap \bigcap \E{}$
  is a $Z$-set in $\bigcap \E{}$.  We say that $A$ is a
  \df{ZF-@$Z(\F{})$-!set}{$Z(\F{})$-set} if it is a $Z(\F{})$-set in
  $\bigcup \F{}$.  We say that a map is a \df{ZF-@$Z(\F{})$-!embedding
    into $X$}{$Z(\F{})$-embedding into~$X$} if it is an embedding and
  its image is a $Z(\F{})$-set in $X$. We say that it is a
  \df{ZF-@$Z(\F{})$-!embedding}{$Z(\F{})$-embedding} if it is a
  $Z(\F{})$-embedding into $\bigcup \F{}$.
\end{definition*}

\begin{remark}\label{rem:ZU-sets}
  By proposition~\ref{pro:Z is a local property}, if $\U{}$ is a
  collection of open subsets of a space $X$, then every $Z$-set in $X$ is 
  a $Z(\U{})$-set in $X$.
\end{remark}

The following propositions are direct corollaries of
propositions~\ref{pro:compact is Z}, \ref{pro:closed Zsigma is Z}
and~\ref{pro:Z is a local property}.

\begin{proposition}\label{pro:compact is ZF}
  If $A$ is a compact subset of a space $X$ and $\F{}$ is a closed
  $\N{n}$-cover of $X$, then $A$ is a $Z(\F{})$-set.
\end{proposition}

\begin{proposition}\label{pro:closed ZFsigma is ZF}
  If a closed subset of a space $X$ is a countable union of
  $Z(\F{})$-sets in $X$, then it is a $Z(\F{})$-set in $X$.
\end{proposition}

Recall that by definition~\ref{def:restricted cover}, $\F{} / U = \{
F_i \cap U \}_{i \in I}$ for each collection $\F{} = \{ F_i \}_{i \in
  I}$ and each set $U$.

\begin{proposition}\label{pro:ZF is a local property}
  Let $A$ be a closed subset of a space $X$ and let $\F{}$ be a closed
  cover of $X$. If $A$ is a $Z(\F{})$-set, then for each open subset
  $U$ of $X$, $A \cap U$ is a $Z(\F{}/U)$-set. If each point of $A$
  has a neighborhood $U$ in $X$ such that $A \cap U$ is a
  $Z(\F{}/U)$-set, then $A$ is a $Z(\F{})$-set.
\end{proposition}
%

\begin{remark}\label{rem:closed locally compact is ZF}
  Propositions~\ref{pro:compact is ZF} and~\ref{pro:ZF is a local
    property} imply that if $\F{}$ is a closed $\N{n}$-cover of a
  space $X$, then every closed locally compact subset of $X$ is a
  $Z(\F{})$-set.
\end{remark}

\section{Approximation within a cover}

\begin{definition*}
  Assume that $\U{}$ and $\F{}$ are covers of a space $Y$ and that~$\U{}$
  is open. We say that a map $g \colon X \to Y$ is a \df{approximation
    within a cover}{$\U{}$-approximation within $\F{}$} of a map $f
  \colon X \to Y$ if it is $\U{}$-close to $f$ and $gf^{-1}(F) \subset
  F$ for each $F \in \F{}$.  We say that $f$ is \df{approximable
    within a cover}{approximable within~$\F{}$} by embeddings (closed
  embeddings, $Z$-embeddings, etc.) if for each open cover $\U{}$ it
  admits a $\U{}$-approximation within $\F{}$ that is an embedding
  (closed embedding, $Z$-embedding, etc.).
\end{definition*}

Note the asymmetry: the assumption that $g$ is a
$\U{}$-approximation of $f$ within $\F{}$ doesn't imply that $f$ is
a $\U{}$-approximation of $g$ within~$\F{}$.

\begin{theorem}
  \label{thm:approximation within a cover} Let $\G{}$ be a closed
  countable star-finite locally finite $\N{n}$-collection in a
  space~$Y$. If~$\G{}$ covers $Y$, then every map~$f$ from an
  $n$-dimensional Polish space into~$Y$ is approximable within $\G{}$
  by $Z(\G{})$-embeddings.  Moreover, if~$f$ restricted to a closed
  subset $A$ of its domain is a $Z(\G{})$-embedding, then the
  approximating $Z(\G{})$-embeddings can be taken to coincide with~$f$
  on~$A$.
\end{theorem}

\begin{proof}
  Let $\G{} = \{ G_i \}_{i \in I}$ and let $G_J = \bigcap_{j \in J}
  G_j$ for each $J \subset I$.  On this occasion, we allow $J$ to be
  the empty set and let $G_\emptyset = Y$. If $G_J$ is non-empty, then
  the space $C(I^n, G_J)$ of maps from $I^n$ into $G_J$ is separable
  (cf. proposition~\ref{pro:countable homotopy groups}). By the
  assumption that $f(A)$ is a $Z(\G{})$-set, the intersection of
  $f(A)$ with $G_J$ is a $Z$-set in $G_J$, hence we can pick a
  countable dense subset of $C(I^n, G_J)$ whose elements have images
  disjoint from $f(A)$. Let~$D$ denote the union of images of maps
  belonging to these countable dense subsets, over all $J \subset I$.
  By theorem~\ref{thm:Z characterization},
  \begin{enumerate}
  \item every closed embedding with image disjoint from $D$ is a
    $Z(\G{})$-embedding.
  \end{enumerate}

  Since~$\G{}$ is countable and locally finite, $G_J$ is non-empty
  only for countably many $J$'s, hence $D$ is sigma-compact. By the
  assumptions, $G_J$ is an $\N{n}$-space, hence by
  proposition~\ref{pro:compact is Z}, the intersection of $D$ with
  $G_J$ is a $Z_\sigma$-set in $G_J$. By (1) and by
  corollary~\ref{cor:z-approximation},
  \begin{enumerate}
    \addtocounter{enumi}{1}
  \item for every $J \subset I$ every map from $n$-dimensional Polish
    space into $G_J$ is approximable by closed embeddings with images
    disjoint from $D \cap G_J$.  Moreover, if the map restricted to a
    closed subset $A$ of its domain is a closed embedding with image
    disjoint from $D \cap G_J$, then such approximations exist $\rel
    A$.
  \end{enumerate}
  
  Let $\U{}$ be an open cover of~$Y$. We will show that the given
  map~$f$ can be $\U{}$-approximated within $\G{}$ by a
  $Z(\G{})$-embedding.  Let $\mathcal{I} = \{ J \subset I \colon J
  \neq \emptyset, G_J \neq \emptyset \}$. Order $\mathcal{I}$ into a
  sequence $\{ J_k \}_{k \in \mathbb{N}}$ non increasing in the order
  by inclusion, i.e. such that if $J_l \supsetneq J_k$, then $l < k$.
  Such ordering exists because $\G{}$ is star-finite. By
  theorem~\ref{thm:n-homotopy extension theorem}, there is a sequence
  of pairs of covers $(\V{k}, \U{k})$ such that for each $k$
  \begin{enumerate}\addtocounter{enumi}{2}
    \item $\V{k}$ and $\U{k}$ are open covers of $G_{J_k}$ that
    satisfy condition $(*)_{n}$ of theorem~\ref{thm:n-homotopy
    extension theorem},
    \item $\st \U{k} \prec \U{}$ and for every $l$ if $J_l
    \supsetneq J_k$, then $\st \U{l} \prec \V{k}$.
  \end{enumerate}

  The construction of a sequence $(\U{k}, \V{k})$ starts with $k$'s that
  correspond to sets $J_k$ of cardinality one and continues backwards by a
  recursive application of theorem~\ref{thm:n-homotopy extension theorem}.

  Let $f$ denote the given map from an $n$-dimensional Polish
  space~$X$ into~$Y$. Let $F_i = f^{-1}(G_i)$ for every $i \in I$. We
  let $g_0 = f_{|A}$ and recursively define a sequence of maps $g_k$
  in such a way that $g_k$ extends $g_{k-1}$ over $F_{J_k}$ and
  \begin{enumerate}\addtocounter{enumi}{4}
  \item $g_k(F_{J_k}) \subset G_{J_k}$ and $g_k$ is $\st \U{k}$-close
    to $f$ on $F_{J_k}$,
  \item $g_k$ is a closed embedding with image disjoint from $D$.
  \end{enumerate}

  In $k$th step of the construction we extend $g_{k-1}$ over
  $F_{J_k}$. By (5) and by the order of $J_k$'s, $g_{k-1}(F_{J_k})
  \subset G_{J_k}$. By the construction, the intersection of $F_{J_k}$
  with the domain of $g_{k-1}$ is equal to $A \cup \bigcup_{l < k, J_l
    \supset J_k} F_{J_l}$. By the choice of $g_0$, $g_{k-1}$ is equal
  to $f$ on the set $A$. By (5), $g_{k-1}$ is $\st \U{l}$-close to $f$
  on $F_{J_l}$ for each $l < k$ such that $J_l \supset J_k$. By (4),
  $\st \U{l}$ refines $\V{k}$, hence $g_{k-1}$ is $\V{k}$-close to $f$
  on the intersection of $F_{J_k}$ with the domain of $g_{k-1}$.
  By~(3), the restriction of $g_{k-1}$ to $F_{J_k}$ extends over
  $F_{J_k}$ to a map from $F_{J_k}$ into $G_{J_k}$ that is
  $\U{k}$-close to $f$. Let $\widetilde g_{k-1}$ be an union of any
  such extension with the restriction of $g_{k-1}$ to $B =
  g_{k-1}^{-1}(G_{J_k})$. By (6), the restriction of $\widetilde
  g_{k-1}$ to $B$ has image disjoint from $D$. By (2), $\widetilde
  g_{k-1}$, defined \marginpar[(*)]{(*)} on $B \cup F_{J_k}$, is
  approximable $\rel B$ by a closed embedding into $G_{J_k}$ that is
  $\st \U{k}$-close to $f$ and has image disjoint from $D$. We let
  $g_k$ to be equal to this embedding on $F_{J_k}$ and equal to
  $g_{k-1}$ on $\dom g_{k-1} \setminus F_{J_k}$. The definition of $B$
  guarantees that $g_k$ is a closed embedding.  Conditions (5) and (6)
  are satisfied directly from the construction.

  Let $g = \bigcup_{k \in \mathbb{N}} g_k$. It is a closed embedding
  because $\G{}$ is locally finite.  By (6), it has image disjoint
  from $D$. Hence, by (1), it is a $Z(\G{})$-embedding. By (4) it is
  $\U{}$-close to $f$ and by (5) it is an approximation of $f$ within
  $\G{}$. Since $g_{|A}$ is equal to $f$, we are done.
\end{proof}

We can draw two important corollaries from
theorem~\ref{thm:approximation within a cover}.

\begin{corollary}
  \label{cor:approximation within a cover}
  If $\F{}$ is a closed countable star-finite locally finite
  $\N{n}$-cover of a space~$X$ and~$Z$ is a $Z(\F{})$-set, then every
  map of an $n$-dimensional Polish space into~$X$ is approximable
  within~$\F{}$ by closed embeddings with images disjoint from~$Z$.
\end{corollary}
\begin{proof}
  Let $f$ be the given map. We can assume that the domain of~$f$ is
  separated from~$Z$, so $f \cup id_Z$ is a well-defined map. By
  theorem~\ref{thm:approximation within a cover} it is approximable
  within $\F{}$ and $\rel Z$ by a closed embedding.  The restriction
  of such approximation to the domain of $f$ is an embedding that we
  were looking for.
\end{proof}

The second corollary follows from theorem~\ref{thm:sum theorem for
  ane-spaces} and theorem~\ref{thm:approximation within a cover}
applied to the cover $\G{} = \{ X, Y \}$ of $X \cup Y$.
\begin{corollary}
  \label{cor:pasting}
  If $X$, $Y$ and $X \cap Y$ are $\N{n}$-spaces and $X \cap Y$ is closed both
  in~$X$ and in~$Y$, then $X \cup Y$ is an $\N{n}$-space.
\end{corollary}

\section[$Z$-collections]{$Z$-collections\protect\footnote{we shall
    not use the theorem proved in ths section until the third part of
    the paper.}}

We prove that under stronger assumptions we can require in the
statement of theorem~\ref{thm:approximation within a cover} not only
that the constructed map~$g$ approximates the given map~$f$
within~$\G{}$, but also that~$f$ approximates~$g$ within~$\G{}$. The
stronger assumption is that~$\G{}$ is a $Z$-collection in~$Y$.

\begin{definition*}\label{def:z-collection}
  We say that a collection $\F{}$ of closed subsets of a topological
  space $X$ is a \df{Z-@$Z$-!collection}{$Z$-collection} if every
  element of $\F{}$ is a $Z$-set in $X$ and if for every $\mathcal{A}
  \subsetneq \mathcal{B} \subset \F{}$ the set $\bigcap \mathcal{B}$
  is a $Z$-set in $\bigcap \mathcal{A}$.
\end{definition*}

Observe that a locally finite $Z$-collection cannot, by the
definition, cover the underlying space. Hence, we need the following,
modified statement.

\begin{theorem}
  \label{thm:approximation within a Z-collection}
  If the assumption that $\G{}$ covers $Y$ in theorem~\ref{thm:approximation
    within a cover} is replaced by the assumption that $Y$ is an $\N{n}$-space
  and if we assume that $\G{}$ is a $Z$-collection, then the approximating map
  $g$ can be taken so that $g^{-1}(G_i) = f^{-1}(G_i)$ for each $i \in I$.
\end{theorem}
\begin{proof}
  The proof is an easy modification of the proof of
  theorem~\ref{thm:approximation within a cover}. The necessary
  changes are to change condition (5) to
  \begin{enumerate}
  \item[(5')] $f^{-1}(G_{J_l}) = g_k^{-1}(G_{J_l})$ for each $l \leq
    k$ and $g_k$ is $\st \U{k}$-close to $f$ on $F_{J_k}$,
  \end{enumerate}
  and to change the sentence marked by (*) on the margin to
  \begin{quote}
    By (2), $\widetilde g_{k-1}$, defined on $B \cup F_{J_k}$, is
    approximable $\rel B$ by a closed embedding into $G_{J_k}$ that is
    $\st \U{k}$-close to $f$, has image disjoint from $D$ and maps
    $F_{J_k} \setminus B$ into $G_{J_k} \setminus \bigcup_{l < k}
    G_{J_l}$.
  \end{quote}

  By the assumption $\G{}$ is a $Z$-collection, hence the intersection
  of $\bigcup_{l < k} G_{J_l}$ with $G_{J_k}$ is a $Z$-set in
  $G_{J_k}$. Hence, by corollary~\ref{cor:z-approximation}, such
  approximation exists. The assumption that the approximation maps
  $F_{J_k} \setminus B$ into $G_{J_k} \setminus \bigcup_{l < k}
  G_{J_l}$ guarantees that $f^{-1}(G_{J_l}) = g_k^{-1}(G_{J_l})$ for
  each $l < k$.
\end{proof}

We will show how to construct $Z$-collections in $\N{n}$-spaces in
section~\ref{sec:adjustment to a Z-collection} in
chapter~\ref{ch:basic constructions}.


\chapter{Constructing closed $\N{n}$-covers}
\label{ch:swelling}

\section{Adjustment of a collection}

By theorem~\ref{thm:open subsets of n-spaces}, every open subset of
an $\N{n}$-space is an $\N{n}$-space. Obviously the property of being
an $\N{n}$-space is not inherited by closed subsets of $\N{n}$-spaces.
However by the strong universality of $\N{n}$-spaces, the inclusion of
an arbitrary open subset of an $\N{n}$-space is approximable by closed
embeddings and images thereof are closed $\N{n}$-subsets of the space.

Let $U$ be an open and $A$ a closed neighborhood of a point $x$ in an
$\N{n}$-space $X$ such that $A \subset U$. By
proposition~\ref{pro:approximation rel inverse image}, the inclusion
of $U$ into $X$ can be approximated $\rel A$ by a closed embedding.
The image of such an approximation (an \emph{adjustment} of~$U$) is a
closed $\N{n}$-neighborhood of $x$. In the sequel, we shall apply this
technique to collections of sets, in a way described by the following
definition.

\begin{definition*}
  Let $\F{} = \{ F_i \}_{i \in I}$ and $\G{} = \{ G_i \}_{i \in I}$ be
  arbitrary collections of subsets of a space~$X$ and let $\U{}$ be an open
  cover of $X$.
  \begin{enumerate}
  \item[(a)] We say that $\G{}$ is a
    \df{U-@$\U{}$-!adjustment}{$\U{}$-adjustment}\index{adjustment} of
    $\F{}$ if for each $J \subset I$ there is a homeomorphism of $F_J
    = \bigcap_{j \in J} F_j$ onto~$G_J = \bigcap_{j \in J} G_j$ that
    is $\U{}$-close to the inclusion $F_J \subset X$.
  \item[(b)] If each of these homeomorphisms can be taken to be the
    identity on a closed set $A$, then we say that $\G{}$ is an
    \df{U-@$\U{}$-!adjustment $\rel A$}{adjustment of~$\F{}$ $\rel
      A$}.\index{adjustment $\rel A$}
  \item[(c)] If additionally $G_i = F_i$ for every $i$ in a given set $J
    \subset I$, then we say that $\G{}$ is an \df{adjustment with a fixed
      $J$}{adjustment of $\F{}$ with fixed $J$}.
  \end{enumerate}
\end{definition*}

A convenient property of the notion of an adjustment is that for an arbitrary
class~$\mathcal{C}$ of topological spaces every adjustment of a
$\mathcal{C}$-collection is a $\mathcal{C}$-collection.

\begin{lemma}
  \label{lem:local adjustment}
  For each open subset $U$ of a metric space $X$ there exists an open
  cover~$\U{}$ of~$U$ such that if $A$ is a (closed) subset of $X$,
  then every $\U{}$-approximation of the inclusion $A \cap U \subset
  U$ by a (closed) embedding of $A \cap U$ into $U$ extends over~$A$,
  by the identity on $A \setminus U$, to a (closed) embedding of $A$
  into $X$.  Hence every (closed in $U$) $\U{}$-adjustment of a
  restriction $\F{} / U$ of a (closed) collection $\F{}$ of subsets
  of~$X$ extends in a natural way to a (closed) adjustment $\rel X
  \setminus U$ of the entire collection.
\end{lemma}
\begin{proof}
  The first statement follows from lemma~\ref{lem:local approximation}
  applied to the inclusion $A \subset X$. Let $\F{} = \{ F_i \}_{i \in
    I}$ be a (closed) collection of subsets of~$X$ and let $\G{} = \{
  G_i \}_{i \in I}$ be a (closed in $U$) $\U{}$-approximation of $\F{}
  / U = \{ F_i \cap U \}_{i \in I}$.  We have to show that the
  collection $\{ G_i \cup (F_i \setminus U) \}_{i \in I}$ is a
  (closed) $\U{}$-adjustment $\rel X \setminus U$ of $\F{}$. Since
  $\G{}$ is a $\U{}$-adjustment of $\F{} / U$, for each $J \subset I$
  there exists a homeomorphism $h_J$ from $F_J \cap U$ onto $G_J$ that
  is $\U{}$-close to the inclusion of $F_J \cap U$ into $U$. By the
  first statement, the map $h_J \cup (F_J \setminus U)$ is a (closed)
  embedding. This embedding is $\U{}$-close to the inclusion and its
  image is equal to $G_J \cup (F_J \setminus U)$. We are done.
\end{proof}

\section{Limits of sequences of adjustments}

\begin{definition}\label{def:definition of limit}
  For every sequence $\F{k} = \{ F^k_i \}_{i \in I}$ we let
  \[
  \lim_{k \to \infty} \F{k} = \{ \bigcap_{j \geq 1} \bigcup_{k \geq j} F^k_i
  \}_{i \in I}.
  \]
\end{definition}

\begin{proposition}
  \label{pro:sequence of adjustments}
  Assume that $\V{} = \{ V_k \}_{k \geq 1}$ is a locally finite
  collection of open subsets of a Polish space $X$. For every open
  cover~$\U{}$ of~$X$ there exists a sequence $\U{1}, \U{2}, \ldots$
  of open covers of $X$ such that if $\F{k}$ is a $\U{k}$-adjustment
  $\rel X \setminus V_k$ of $\F{k-1}$, then $\lim_{k \to \infty}
  \F{k}$ is a $\U{}$-adjustment of~$\F{0}$.
\end{proposition}
\begin{proof}
  Without a loss of generality we may assume that a second star of
  each element of~$\U{}$ intersects only finitely many elements of
  $\V{}$. Let $\U{1}, \U{2}, \ldots$ be a sequence of covers of $X$
  obtained via lemma~\ref{lem:limit close to the inclusion} applied
  to~$\U{}$. We will show that it satisfies our claim.

  Let $\F{k} = \{ F^k_i \}_{i \in I}$ and $\F{\infty} = \lim_{k \to
    \infty} \F{k}$. Let $F^k_J = \bigcap_{j \in J} F^k_j$. By the
  assumptions, for each $J \subset I$ there exists a homeomorphism
  $h^k_J \colon F^{k-1}_J \to F^{k}_J$ that is $\U{k}$-close to the
  inclusion of $F^{k-1}_J$ into $X$. By lemma~\ref{lem:limit close to
    the inclusion} and by the choice of $\U{k}$'s, the limit
  $h^\infty_J = \lim_{k \to \infty} h^k_J \circ \cdots \circ h^2_J
  \circ h^1_J$ exists, is continuous and $\U{}$-close to the inclusion
  of $F^0_J$ into $X$.

  Observe that if a map is $\U{}$-close to the identity, then to check
  that it is a homeomorphism (onto its image) it suffices to check it
  is a homeomorphism on $\st_\U{} U$ for every $U$ in $\U{}$. By the
  choice of $\U{}$, for each $U$ in $\U{}$ the second star $\st^2_\U{}
  U$ intersects only finitely many elements of $\V{}$. Since each
  $h^k_J \circ \cdots \circ h^2_J \circ h^1_J$ is $\U{}$-close to the
  identity, the sequence $h^J_k \circ \cdots \circ h^J_1 \circ h^J_0$
  restricted to $F^0_J \cap \st_\U{} U$ stabilizes after finitely many
  steps. By the assumption every $h^J_k$ is a homeomorphism, so
  $h^\infty_J$ restricted to $\st_\U{} U$ is a homeomorphism. We are
  done, because the image of $h^\infty_J$ is equal to $F^\infty_J$, by
  local finiteness of~$\V{}$.
\end{proof}

\section{Construction of a closed $\N{n}$-swelling}

Our first application of the technique of adjustments is a
construction of a closed interior $\N{n}$-cover given some constraints
on its size and structure.

\begin{lemma}\label{lem:closed adjustment}
  If $\G{}$ is a closed star-finite cover of an $\N{n}$-space~$Y$
  and~$\V{}$ is an open swelling of $\G{}$, then for each open
  cover~$\U{}$ of~$Y$ and each $Z$-set $Z$ in $Y$ there exists a
  closed $\U{}$-adjustment~$\F{}$ of~$\V{}$ that is a swelling
  of~$\G{}$ and such that $Z$ is a $Z(\F{})$-set.
\end{lemma}

Observe that the assumption that $\G{}$ is star-finite and has an open
swelling implies that it is locally finite.

\begin{proof} 
  Let $\G{} = \{ G_i \}_{i \in I}$, $\V{} = \{ V_i \}_{i \in I}$ and
  $V_J = \bigcap_{j \in J} V_j$ for each non-empty subset $J$ of $I$.
  Let $\mathcal{I} = \{ J \subset I \colon J \neq \emptyset \text{ and
  } V_J \neq \emptyset \}$.  Order $\mathcal{I}$ into a sequence $\{
  J_k \}_{k \geq 1}$ non increasing in the order by inclusion.  Such
  ordering exists because $\V{}$ is star-finite. Let $W_k = \bigcup_{i
    \in I \setminus J_k } V_i$. Let~$B_k$ be a boundary of $V_{J_k}
  \setminus W_k$ taken in $V_{J_k}$ (it is not necessarily closed in
  $Y$). Let
  \[
  U_k = (V_{J_k} \setminus \bigcup_{i \in I \setminus J_k} G_i) \cap
  \{ y \in Y \colon d(y, B_k) < d(y, Y \setminus V_{J_k}) \},
  \]
  where $d(y, A)$ is the distance of a point $y$ from a set $A$ (we
  let $d(y, \emptyset) = \infty$). By the assumptions, $\V{}$ is a
  swelling of $\G{}$, so $\bigcup_{i \in I \setminus J_k} G_i \subset
  W_k$.  Hence $U_k$ is an open neighborhood of $B_k$. We shall prove
  that if $U_k$ intersects $U_l$, then either $J_k \subset J_l$ or
  $J_l \subset J_k$. Let $y$ be an element of $U_k \cap U_l$. If $B_k
  \cap V_{J_l} = \emptyset$ and $B_l \cap V_{J_k} = \emptyset$, then
  $d(y, B_k) < d(y, Y \setminus V_{J_k}) \leq d(y, B_l) < d(y, Y
  \setminus V_{J_l}) \leq d(y, B_k)$ - a contradiction. Hence either
  $B_k$ intersects $V_{J_l}$ or $B_l$ intersects $V_{J_k}$. Assume the
  former case, i.e. $B_k \cap V_{J_l} \neq \emptyset$. By the
  definition, $B_k$ is disjoint from $W_k$.  Since $V_{J_l} =
  \bigcap_{i \in J_l} V_i$ intersects $B_k$, it is not a subset of
  $W_k = \bigcup_{i \in I \setminus J_k} V_i$, hence $J_l$ must be
  disjoint from $I \setminus J_k$, so it must be a subset of $J_k$.
  Same argument applies in the latter case, i.e. when $B_l \cap
  V_{J_k} \neq \emptyset$.
  
  Let $\U{1}, \U{2}, \ldots$ be a sequence of covers of $Y$ obtained
  via proposition~\ref{pro:sequence of adjustments} applied to a
  locally finite collection $\V{}$ and an open cover $\U{}$ of $Y$. If
  $P$ and $Q$ are open subsets of $Y$, then the boundary of $P
  \setminus Q$ taken in $P$ is equal to the boundary of $Q$ taken in
  $P \cup Q$. We let $P = V_{J_k}$ and $Q = W_k$, and obtain the
  equality $B_k = \Bd_{V_{J_k} \cup W_k} W_k$.  Since $U_k$ is a
  neighborhood of $B_k$, $W_k \setminus U_k$ is closed in $W_k \cup
  V_{J_k}$.  By theorem~\ref{thm:open subsets of n-spaces}, $W_k \cup
  V_{J_k}$ is an $\N{n}$-space. Therefore, by
  proposition~\ref{pro:approximation rel inverse image} applied to $A
  = W_k \setminus U_k$, $V = W_k$ and $Z_Y = Z \cap U_k$, the
  inclusion $W_k \subset W_k \cup V_{J_k}$ has a $\U{k}$-approximation
  $h_k \colon W_k \to W_k \cup V_{J_k}$ such that the following
  condition is satisfied.
  \begin{enumerate}
  \item[(h$_k$)] $h_k$ is a closed embedding into $W_k \cup V_{J_k}$,
  is~$\U{k}$-close to the inclusion, is equal to the inclusion on $W_k
  \setminus U_k$ and has image disjoint from $Z \cap U_k$.
  \end{enumerate}

  We shall recursively construct a sequence $\F{k} = \{ F^k_i \}_{i
    \in I}$ of covers of $Y$, starting with $\F{0} = \V{}$, such that
  for $k > 0$ the following conditions are satisfied.
  \begin{enumerate}
  \item[($1_k$)] $\F{k}$ is a $\U{k}$-adjustment of $\F{k-1}$ $\rel Y
    \setminus U_k$ with fixed $J_k$.
  \item[($2_k$)] $\F{k}$ restricted to $\bigcup_{0 < l \leq k}
    V_{J_l}$ is a closed cover, $Z$ is a $Z(\F{k})$-set and $\G{}$ is
    a shrinking of~$\F{k}$.
  \end{enumerate}

  Observe that $Z$ is a $Z(\F{0})$-set by remark~\ref{rem:ZU-sets} and
  $\G{}$ is a shrinking of $\F{0}$. Before we construct the sequence,
  we prove that for each $k \geq 0$ if ($1_l$) and ($2_l$) are
  satisfied for each $l$ greater than $0$ and less than or equal to
  $k$, then the following conditions hold:
  \begin{enumerate}\addtocounter{enumi}{2}
  \item[($3_k$)] $\bigcup_{i \in I \setminus J_{k+1}} F^{k}_i$ is a
    subset of $W_{k+1}$,
  \item[($4_k$)] $U_{k+1} \cup (V_{J_{k+1}} \setminus W_{k+1})$ is a
    subset of $F^k_{J_{k+1}}$.
  \end{enumerate}
  
  To prove ($3_k$) fix $j \in I \setminus J_{k+1}$. We will show by
  induction on $l$ that $F^l_j \subset W_{k+1}$ for each $0 \leq l
  \leq k$.  It is true for $l = 0$, because $F^0_j = V_j$. Let $l >
  0$. By ($1_l$), $\F{l}$ is equal to $\F{l-1}$ on the complement of
  $U_l$, hence $F^l_j \subset F^{l-1}_j \cup U_l$. By the inductive
  assumption $F^{l-1}_j \subset W_{k+1}$.  If $l \leq k$, then $J_l
  \setminus J_{k+1}$ is non-empty, by the order of $J_k$'s.  Hence
  $U_l \subset V_{J_l} \subset V_{{J_l} \setminus J_{k+1}} \subset
  W_{k+1}$. Hence $F^l_j \subset F^{l-1}_j \cup U_l \subset W_{k+1}$.

  We will show by induction on $l$ that $U_{k+1} \cup (V_{J_{k+1}}
  \setminus W_{k+1})$ is a subset of $F^l_{J_{k+1}}$ for each $0 \leq
  l \leq k$, thus proving ($4_k$). From the definition $U_{k+1}
  \subset V_{J_{k+1}} = F^0_{J_{k+1}}$, so the assertion is true for
  $l = 0$. Let $0 < l \leq k$. We either have $J_{k+1} \not\subset
  J_l$ or or $J_{k+1} \subset J_l$. By ($1_l$), $\F{l}$ is an
  adjustment of $\F{l-1}$ with fixed $J_l$ so in the latter case,
  $F^l_{J_{k+1}} = F^{l-1}_{J_{k+1}}$ and the assertion is true. In
  the former case, $U_l$ is disjoint from $U_{k+1}$ and $U_l \subset
  V_{J_l} \subset W_{k+1}$. Hence $U_{k+1} \cup (V_{J_{k+1}} \setminus
  W_{k+1})$ is disjoint from $U_l$. By ($1_l$), $\F{l}$ is equal to
  $\F{l-1}$ on the complement of $U_l$, so the inductive step is done.

  The construction. Fix $k > 0$ and assume that we already constructed
  $\F{k-1}$. Let
  \[
  F^k_i = \left\{ \begin{array}{ll}
      h_k(F^{k-1}_i) & i \in I \setminus J_k \\
      F^{k-1}_i & i \in J_k. \\
    \end{array} \right.
  \]

  Consider $J \subset I$. If $J \subset J_k$, then $F^k_J =
  F^{k-1}_J$. If $J \not\subset J_k$, then a direct computation using
  ($4_{k-1}$) and (h$_k$) shows that $F^k_J = h_k(F^{k-1}_J)$. Hence
  by ($3_{k-1}$), the restriction of $h_k$ to $F^{k-1}_J$ is a
  homeomorphism onto $F^k_J$.  By the definition of a
  $\U{k}$-adjustment, ($1_k$) is satisfied.

  By (h$_k$), $h_k$ is a closed embedding into $W_k \cup V_{J_k}$.
  Hence its restriction to $h^{-1}_k(\bigcup_{0 < l \leq k} V_{J_l}) =
  \bigcup_{0 < l < k} V_{J_l}$ is a closed embedding into $\bigcup_{0
    < l \leq k} V_{J_l}$. By ($2_{k-1}$), $\F{k-1}$ restricted to
  $\bigcup_{0 < l < k} V_{J_l}$ is closed. Hence for each $i \in I
  \setminus J_k$, $F^k_i \cap \bigcup_{0 < l \leq k} V_{J_l} =
  h_k(F^{k-1}_i \cap \bigcup_{0 < l < k} V_{J_l})$ is closed in
  $\bigcup_{0 < l \leq k} V_{J_l}$. If $i \in J_k$, then $F^k_i =
  F^{k-1}_i$. By ($4_{k-1}$), $V_{J_k} \setminus \bigcup_{0 < l < k}
  V_{J_l}$ is a subset of $F^{k-1}_i$. The intersection of $F^{k-1}_i$
  with $\bigcup_{0 < l < k} V_{J_l}$ is closed in $\bigcup_{0 < l < k}
  V_{J_l}$, so $F^{k-1}_i \cap \bigcup_{0 < l \leq k} V_{J_l} =
  (V_{J_k} \setminus \bigcup_{0 < l < k} V_{J_l}) \cup (F^{k-1}_i \cap
  \bigcup_{0 < l < k} V_{J_l})$ is closed in $\bigcup_{0 < l \leq k}
  V_{J_l}$.  Hence $\F{k}$ restricted to $\bigcup_{0 < l \leq k}
  V_{J_l}$ is a closed cover.

  By ($2_{k-1}$), $G_i \subset F^{k-1}_i$ for each $i \in I$. If $i
  \in J_k$, then $F^k_i = F^{k-1}_i$, so $G_i \subset F^k_i$. If $i
  \in I \setminus J_k$, then $G_i$ is disjoint from $U_k$ by the
  definition. By ($1_k$), $F^k_i$ is equal to $F^{k-1}_i$ on the
  complement of $U_k$, hence $G_i \subset F^k_i$. Therefore $\G{}$
  refines $\F{k}$.

  By ($2_{k-1}$), $Z$ is a $Z(\F{k-1})$-set. Hence for each $J \subset
  I$, $Z \cap F^{k-1}_J$ is a $Z$-set in $F^{k-1}_J$. If $J \subset
  J_k$, then $F^k_J = F^{k-1}_J$, so $Z \cap F^k_J$ is a $Z$-set in
  $F^k_J$ as well. If $J \not\subset J_k$, then $F^k_J =
  h_k(F^{k-1}_J)$. By (h$_k$), the image of $h_k$ is disjoint from $Z
  \cap U_k$ and $h_k$ is equal to the inclusion on the complement of
  $U_k$. Hence $Z \cap F^k_J = Z \cap (F^{k-1}_J \setminus U_k)$. By
  ($2_{k-1}$), the latter set is a $Z$-set in $F^{k-1}_J$, hence its
  image under $h_k$ is a $Z$-set in $h_k(F^{k-1}_J) = F^k_J$. But
  $h_k$ is the identity map on this set, hence $Z \cap F^k_J$ is a
  $Z$-set in $F^k_J$. Hence $Z$ is a $Z(\F{k})$-set.

  By ($1_k$), by proposition~\ref{pro:sequence of adjustments} and by
  the choice of $\U{k}$'s, the limit $\F{} = \lim_{k \to \infty}
  \F{k}$ exists and is a $\U{}$-adjustment of $\F{0} = \V{}$. By
  ($3_k$), $\F{}$ is closed. By ($2_k$), it is a swelling of $\G{}$.
  By ($2_k$) and by proposition~\ref{pro:ZF is a local property}, $Z$
  is a $Z(\F{})$-set. We are done.
\end{proof}

\begin{theorem}
  \label{thm:n-swelling} Every closed star-finite locally finite cover
  that refines an open cover~$\U{}$ of an $\N{n}$-space~$X$ has a swelling
  $\F{}$ that refines~$\U{}$ and is a closed locally finite interior
  $\N{n}$-cover. Moreover, we may require that a specified $Z$-set is a
  $Z(\F{})$-set.
\end{theorem}
\begin{proof}
  Let $\U{}$ be an open cover of an $\N{n}$-space $X$ and let $\G{} =
  \{ G_i \}_{i \in I}$ be a closed star-finite locally finite
  refinement of $\U{}$. We may assume that $\G{}$ is an interior
  cover~\cite[Exercise~7.1.G]{engelking1989}. By the same argument,
  there exists an open locally finite swelling $\V{} = \{ V_i \}_{i
    \in I}$ of~$\G{}$ such that the closure of each element of $\V{}$
  lies in an element of $\U{}$. Let~$Z$ be an arbitrary $Z$-set
  in~$X$. By lemma~\ref{lem:closed adjustment}, $\V{}$ can be adjusted
  to a closed cover $\F{}$ of $X$ such that $Z$ is a $Z(\F{})$-set and
  $\G{}$ is a refinement of $\F{}$. Obviously, if the adjustment is
  small enough, then $\F{}$ refines $\U{}$ and we are done.
\end{proof}


\chapter{Carrier and nerve theorems}
\label{ch:carrier and nerve theorems}

Almost entire chapter is cited word-for-word after~\cite{nagorko2006}
and it is included here for completeness.

\section{Regular covers}

An efficient way to investigate properties of a topological space
is to divide it into pieces and examine how they are glued
together. We show how to divide a general topological space and
endow it with a structure that resembles a triangulation.

In order to divide a space that belongs to a class~$\mathcal{C}$
of topological spaces we must decide what a piece is. We want it
to resemble a simplex as much as possible. As a simplex is an
archetype of an absolute extensor, the choice of absolute
extensors for~$\mathcal{C}$ as pieces is quite natural.

\begin{definition*}
  A space~$Y$ is an \df{absolute extensor for a space~$X$}{absolute extensor}
  for a space~$X$ if each map from a closed subset of~$X$ into~$Y$ extends
  over the entire space~$X$. The class of absolute extensors for all spaces
  from a class~$\mathcal{C}$ is denoted
  by~\df{AE(C)@$AE(\protect\mathcal{C})$}{$AE(\mathcal{C})$}. We write
  \df{AE(X)@$AE(X)$}{$AE(X)$} for $AE(\{ X \})$.
\end{definition*}

The following defines regular covers that endow a space with
structures similar to triangulations.

\begin{definition*}
  Let~$\mathcal{C}$ be a class of topological spaces. A locally finite locally
  finite-dimensional closed $AE(\mathcal{C})$-cover is said to be
  \df{cover!regular for a class~$\protect\mathcal{C}$}{regular for the
    class~$\protect\mathcal{C}$}. Recall that a cover is \df{cover!locally
    finite dimensional}{locally finite dimensional} if its nerve is such.
\end{definition*}

Examples of regular covers include a locally finite and locally finite
dimensional cover of an Euclidean space by its closed balls and a
cover of a finite simplicial complex by its simplices.

\section{Carrier theorem}
\label{sec:carrier theorem}

How to extend a partial map from a subcomplex over the entire CW complex? If
each map from the boundary of an Euclidean ball into the codomain extends over
the ball, then the answer is easy: order cells by inclusion and construct an
extension inductively. But the asphericity of the codomain (vanishing of all
its homotopy groups) is a rare luxury. The same technique would work though if
we were able to restrict ranges of the map on individual cells of the CW
complex to aspherical subspaces of the codomain. This idea leads to the notion
of a carrier and to the aspherical carrier theorem~\cite[II \S
9]{lundell1969}. We generalize this notion to arbitrary spaces using a cover
to replace the cell structure in the domain.

\begin{definition*}
  A \df{carrier}{carrier} is a function $C \colon \F{} \to \G{}$ from
  a cover~$\F{}$ of a space~$X$ into a collection~$\G{}$ of subsets of
  a topological space such that for each $\mathcal{A} \subset \F{}$ if
  $\bigcap \mathcal{A} \neq \emptyset$, then $\bigcap_{A \in
  \mathcal{A}} C(A) \neq \emptyset$. We say that a map~$f$ is
  \df{map!carried by a carrier}{carried by~$C$} if it is defined on a
  closed subset of~$X$ and $f(F) \subset C(F)$ for each $F \in \F{}$.
\end{definition*}

\begin{carrier theorem}
  Assume that $C \colon \F{} \to \G{}$ is a carrier such that~$\F{}$ is a
  closed cover of a space~$X$ and~$\G{}$ is an $AE(X)$-cover of another space.
  If~$\F{}$ is locally finite and locally finite dimensional, then each map
  carried by~$C$ extends to a map of the entire space~$X$, also carried
  by~$C$.
\end{carrier theorem}

Special cases of the carrier theorem follow from Michael's
selection theorem, as the multivalued map given by the formula
$F(x) = \bigcap_{F \owns x} C(F)$ is lower semi-continuous.

\begin{proof}
  Let $f_0$ be a map carried by $C$ and let $A = \dom f_0$. Let
  $\{ \E{\gamma} \}_{0 < \gamma < \Gamma}$ be a transfinite
  sequence of all subcollections of~$\F{}$ with non-empty
  intersections, such that the sequence $\{ \delta_\gamma =
  \bigcap \E{\gamma}\}_{0 < \gamma < \Gamma}$ is non decreasing in
  the order by inclusion. Its existence is guaranteed by the
  assumption of local finite dimensionality of~$\F{}$. Let
  $C_\gamma = \bigcap_{E \in \E{\gamma}} C(E)$ and observe that an
  arbitrary map~$f$ is carried by~$C$ if and only if
  $f(\delta_\gamma) \subset C_\gamma$ for each $0 < \gamma <
  \Gamma$. Let~$\delta_\Gamma = \emptyset$. We shall construct a
  transfinite sequence of maps $\{ f_\gamma \colon A \cup
  \bigcup_{0 < \alpha \leq \gamma} \delta_\alpha \to Y \}_{\gamma
  \leq \Gamma}$ such that~$f_\alpha$ extends~$f_\beta$ for all $0
  \leq \beta \leq \alpha \leq \Gamma$ and $f_\gamma(\delta_\gamma)
  \subset C_\gamma$ for each $0 < \gamma < \Gamma$. The
  map~$f_\Gamma$ will be an extension that we are looking for,
  since~$\bigcup_{\gamma < \Gamma} \delta_\gamma = X$.

  We proceed by transfinite induction. Fix $\gamma \leq \Gamma$
  and assume that for each~$\alpha < \gamma$ we already
  constructed~$f_\alpha$. The map $f'_\gamma = \bigcup_{\alpha <
  \gamma} f_\alpha$ is well defined, continuous and its domain is
  closed in~$X$ because maps~$f_\alpha$ agree on intersections of
  their domains and~$\F{}$ is closed and locally finite.
  If~$\gamma = \Gamma$, then~$\delta_\gamma = \emptyset$ and we may
  put $f_\Gamma = f'_\Gamma$. If $\gamma < \Gamma$, then by the
  order of $\delta_\gamma$ and by inductive assumptions~$f'_\gamma$
  maps $\delta_\gamma$ into~$C_\gamma$. The set~$C_\gamma$ is
  non-empty because~$C$ is a carrier and is an absolute extensor
  for~$\delta_\gamma$ because~$\G{}$ is an $AE(X)$-cover. So
  $f'_\gamma$ extends onto~$\delta_\gamma$ to a map~$f_\gamma$
  such that $f_\gamma(\delta_\gamma) \subset C_\gamma$ and our
  construction is finished.
\end{proof}

\begin{definition*}
  A \df{cover!regular for a space~$X$}{cover is regular for a space~$X$} if it
  is regular for the class~$\{ X \}$.
\end{definition*}

\begin{corollary}
  \label{cor:carried homotopy} If a closed cover~$\G{}$ of a
  space~$Y$ is regular for $X \times [0, 1]$, then every two $\G{}$-close maps
  from $X$ into $Y$ are $\G{}$-homotopic.  Moreover such homotopy exists with
  the additional property that if endpoints of its path lie in an element
  of~$\G{}$, then the entire path lies in it.
\end{corollary}
\begin{proof}
  Let $f$ and $g$ be two $\G{}$-close maps from $X$ into $Y$.  Let
  $\F{} = \{ F_G \}_{G \in \G{}}$ be the collection of subsets of $X
  \times [0, 1]$ defined by the formula
  \begin{displaymath}
    F_G = (f^{-1}(G) \cap g^{-1}(G)) \times [0, 1].
  \end{displaymath}
  It is a cover of $X \times [0, 1]$ because $f$ and $g$ are
  $\G{}$-close. Define a carrier~$C \colon \F{} \to \G{}$ by the
  formula $C(F_G) = G$ and a map $F \colon X \times \{ 0, 1 \} \to Y$
  by $F(x, 0) = f(x)$ and $F(x, 1) = g(x)$. By the definition, $F$ is
  carried by $C$ and by the carrier theorem, it admits an extension
  over the entire space $X \times [0, 1]$, also carried by $C$. This
  extension is a $\G{}$-homotopy that satisfies our claim.
\end{proof}

\section{Nerve theorem}
\label{sec:nerve theorems}

Nerve theorems give conditions under which the nerve of a cover is
equivalent to the underlying space. First examples of such
theorems are attributed to K.~Borsuk~\cite[p. 234]{borsuk1948}
(for closed covers) and A.~Weil~\cite[p. 141]{weil1952} (for open
covers), both for homotopy equivalences. Since then several
generalizations were made. First generalizations by
W.~Holszty\'nski~\cite{holsztynski1964} and
J.~N.~Haimov~\cite{haimov1979} relaxed conditions on the cover.
Next weak homotopy equivalences were studied in this context by
M.~McCord~\cite{mccord1967} and weak $n$-homotopy equivalences by
A.~Bj\"orner~\cite{bjorner2003}.

The nerve theorem that we shall use in the present paper is proved in
the next chapter. Here we develop tools used in its proof. As an easy
application, we prove a nerve theorem for closed covers
(theorem~\ref{thm:general nerve theorem}), which is more general than
results previously published in the literature.

\begin{definition*}
  A carrier is \df{carrier!invertible}{invertible} if it is bijective
  and its inverse is a carrier.
\end{definition*}

\begin{lemma}
  \label{lem:invertible carriers} If a collection $\F{} = \{ F_i \}_{i
    \in I}$ is point-finite, then its nerve is well defined and a
  function $F_i \mapsto \bst v(F_i)$ is an invertible carrier, where
  $\bst v(F_i)$ is a barycentric star of a vertex of the nerve of
  $\F{}$ that corresponds to $F_i$.
\end{lemma}
\begin{proof}
  This is a reformulation of lemma~\ref{lem:bst isomorphism}.
\end{proof}

\begin{lemma}
  \label{lem:regular stars} If a simplicial complex is locally finite
  dimensional and is endowed with either weak or metric topology, then
  its cover by barycentric stars of its vertices is regular for the
  class of metric spaces. If a simplicial complex is both locally
  countable and locally finite dimensional and is endowed with the
  metric topology, then its cover by barycentric stars of its vertices
  is regular for the class of normal spaces.
\end{lemma}
\begin{proof}
  Let $K$ be a locally finite dimensional simplicial complex and let
  $\mathcal{B}_K$ be a cover of $K$ by barycentric stars of its
  vertices. By lemma~\ref{lem:bst isomorphism}, the nerve of
  $\mathcal{B}_K$ is isomorphic to $K$. Hence $\mathcal{B}_K$ is
  locally finite dimensional. By lemma~\ref{lem:bst isomorphism},
  $\mathcal{B}_K$ is locally finite. What is left to prove is that
  every intersection of elements of $\mathcal{B}_K$ is either empty or
  is an absolute extensor for the class of metric spaces (or for the
  class of normal spaces if $K$ is countable and endowed with the
  metric topology). By theorem~11.7 on page~109 of~\cite{hu1965} and
  by theorem~10.4 on page~105 of~\cite{hu1965}, every simplicial
  complex is an absolute neighborhood extensor for the class of metric
  spaces.  By theorem~11.7 on page~109 of~\cite{hu1965}, a countable
  simplicial complex with the metric topology is an absolute
  neighborhood extensor for the class of normal spaces. By theorem~7.1
  on page~43 of~\cite{hu1965}, every contractible absolute
  neighborhood extensor is an absolute extensor.  Hence it suffices to
  prove that an intersection of a collection of barycentric stars of
  vertices of a simplicial complex is either empty, or contractible.
  Lemma~\ref{lem:contractible stars} finishes the proof.
\end{proof}

\begin{remark}
  Conditions listed in theorem~11.7 on page~109 of~\cite{hu1965} allow
  to state other variants of lemma~\ref{lem:regular stars}, i.e. if
  the complex is endowed with the metric topology but is not
  necessarily countable, then the cover is regular for the class of
  fully normal spaces, etc.
\end{remark}

Putting everything together we obtain a nerve theorem for homotopy
equivalences. We state the theorem for closed covers, which generalizes a
nerve theorem by J.~N.~Haimov~\cite{haimov1979}. An analogous theorem for open
covers may also be proved. We do not state it here as it turns out to be
equivalent to the nerve theorem by A.~Weil~\cite{weil1952}.

\begin{theorem}
  \label{thm:general nerve theorem}
  Assume that a closed cover~$\F{}$ of a normal space~$X$ is regular
  for the class of metric spaces.  If~$\F{}$ is star-countable,
  then~$X$ and the nerve of~$\F{}$ are homotopy equivalent.
\end{theorem}

The main theorem of~\cite{haimov1979} states the same conclusion under
the additional assumption that~$\F{}$ is star-finite and~$X$ is
paracompact. The topology on the nerve of~$\F{}$ is either weak or
metric; by theorem~$6.4$ on page~$81$ of~\cite{bessagapelczynski1975},
the identity map from a simplicial complex endowed with a weak
topology into itself, endowed with a metric topology, is a homotopy
equivalence.

\begin{proof}
  Let $B \colon \F{} \to \BF$ be an invertible carrier as defined in
  lemma~\ref{lem:invertible carriers}. By the carrier theorem and
  lemma~\ref{lem:regular stars} there exist $\kappa \colon X \to N(\F{})$ and
  $\lambda \colon N(\F{}) \to X$ carried by $B$ and $B^{-1}$ respectively.
  Then $\lambda \circ \kappa$ is carried by $B^{-1} \circ B$ so it is
  $\F{}$-close to $id_X$ and by corollary~\ref{cor:carried homotopy}~$\lambda$
  is a homotopy inverse of~$\kappa$.  Analogously~$\kappa$ is a homotopy
  inverse of~$\lambda$ so~$X$ and~$N(\F{})$ are homotopy equivalent.
\end{proof}


\chapter{Anticanonical maps and semiregularity}

\section{A nerve theorem and the notion of semiregularity}

The crux of the proof of theorem~\ref{thm:general nerve theorem} is
the proof of the existence of a map $\lambda$ carried by the carrier
$B^{-1}$ defined in lemma~\ref{lem:invertible carriers}. These maps
are very important in the sequel.

\begin{definition*}
  A map from the nerve of a closed cover~$\F{}$ into the underlying
  space is an \df{anticanonical map}{anticanonical map of $\F{}$} if
  it is carried by the carrier $B^{-1}$ defined in
  lemma~\ref{lem:invertible carriers}.
\end{definition*}

\begin{proposition}
  \label{pro:anticanonical map}
  If an $AE(n)$-cover is point-finite and has at most $n$-dimensional
  nerve, then it admits an anticanonical map.
\end{proposition}
\begin{proof}
  Let $\F{}$ be a point-finite $AE(n)$-cover with at most
  $n$-dimensional nerve. By lemma~\ref{lem:regular stars}, the cover
  $\BF$ of the nerve of $\F{}$ by barycentric stars of its vertices is
  locally finite (by the definition of a regular cover). By
  lemma~\ref{lem:invertible carriers}, $\BF$ is isomorphic to $\F{}$.
  In particular, it is locally finite dimensional. By the carrier
  theorem, the empty map extends over the nerve of $\F{}$ to a map
  carried by $B^{-1}$. This map is an anticanonical map of $\F{}$.
\end{proof}

The notion of an anticanonical map is dual to the notion of a
canonical map. For closed covers, we define canonical maps as
follows.

\begin{definition*}
  A map from a space into the nerve of its closed cover~$\F{}$ is a
  \df{canonical map}{canonical map of~$\F{}$} if it is carried by the
  carrier~$B$ defined in lemma~\ref{lem:invertible carriers}.
\end{definition*}

\begin{proposition}
  \label{pro:canonical map}
  Every closed locally finite and locally finite dimensional cover of
  a metric space admits a canonical map.
\end{proposition}
\begin{proof}
  Apply lemma~\ref{lem:regular stars} and the carrier theorem.
\end{proof}

If a cover satisfies assumptions of
propositions~\ref{pro:anticanonical map} and~\ref{pro:canonical map},
then it is regular for the class of at most $n$-dimensional metric
spaces (the assumption that its nerve is at most $n$-dimensional is
not needed). For brevity, we employ the following definition.

\begin{definition*}
  We say that a closed cover is \df{cover!$n$-regular}{$n$-regular} if
  it is a locally finite locally finite dimensional $AE(n)$-cover,
  i.e. if it is regular for the class of at most $n$-dimensional
  metric spaces.
\end{definition*}

\begin{nerve theorem}
  \label{thm:nerve theorem}
  The $n$-dimensional skeleton of the nerve of a closed $n$-regular
  cover of an at most $n$-dimensional metric space~$X$ is $n$-homotopy
  equivalent to~$X$.  In particular, if two at most $n$-dimensional
  metric spaces admit isomorphic closed $n$-regular covers, then they
  are $n$-homotopy equivalent.
\end{nerve theorem}
\begin{proof}
  We prove the second assertion first. Let $X$ and $Y$ be at most
  $n$-dimensional metric spaces and let $\F{} = \{ F_i \}_{i \in I}$
  and $\G{} = \{ G_i \}_{i \in I}$ be isomorphic closed $n$-regular
  covers of $X$ and $Y$ respectively. The map $C \colon \F{} \to \G{}$
  that assigns $G_i$ to $F_i$ for each $i \in I$ is an invertible
  carrier.  By the assumption, $\F{}$ is a closed locally finite
  locally finite dimensional cover of $X$ and $\G{}$ is an
  $AE(X)$-cover, because $X$ is at most $n$-dimensional. Hence, by the
  carrier theorem, there exists a map $\kappa \colon X \to Y$, carried
  by~$C$. By the symmetry of assumptions, there exists a map $\lambda
  \colon Y \to X$, carried by $C^{-1}$. Let $\varPhi$ be a map from an
  at most $(n-1)$-dimensional simplicial complex into $X$. The map
  $\lambda \circ \kappa$ is carried by $C^{-1} \circ C$, hence it is
  $\F{}$-close to the identity on $X$. Therefore, by
  corollary~\ref{cor:carried homotopy}, $\lambda \circ \kappa \circ
  \varPhi$ is homotopic to $\varPhi$, hence $\lambda$ is an
  $n$-homotopy inverse of~$\kappa$.  Analogously, $\kappa$ is an
  $n$-homotopy inverse of~$\lambda$, so~$\kappa$ and $\lambda$ are
  $n$-homotopy equivalences and $X$ is $n$-homotopy equivalent to $Y$.

  Let $\BF$ be the cover of the nerve $N(\F{})$ of $\F{}$ by
  barycentric stars of its vertices. Let~$\BFn$ denote the restriction
  of $\BF$ to the $n$-dimensional skeleton of the first barycentric
  subdivision of the nerve of $\F{}$. It is isomorphic to $\BF$, which
  is isomorphic to $\F{}$ by lemma~\ref{lem:invertible carriers}. By
  lemma~\ref{lem:regular stars}, $\BF$ is regular for the class of
  metric spaces.  Hence, by the cellular approximation theorem
  (see~\cite[Corollary 4.12]{hatcher2002}), $\BFn$ is $n$-regular.
  Therefore, by what we proved in the first paragraph, $X$ is
  $n$-homotopy equivalent to the $n$-dimensional skeleton of the first
  barycentric subdivision of the nerve of $\F{}$. By the cellular
  approximation theorem, inclusions of $n$-dimensional skeleton of
  $N(\F{})$ into $N(\F{})$ and of $n$-dimensional skeleton of the
  first barycentric subdivision of $N(\F{})$ into $N(\F{})$ are weak
  $n$-homotopy equivalences.  Hence, the inclusion of the former
  skeleton into the latter skeleton is a weak $n$-homotopy
  equivalence.  By Whitehead's characterization, both skeletons are
  $n$-homotopy equivalent. We are done.
\end{proof}

\begin{remark}
  \label{rem:nerve theorem}
  If a closed $n$-regular cover of an at most $n$-dimensional space has
  an at most $n$-dimensional nerve, then by
  proposition~\ref{pro:canonical map} and by
  proposition~\ref{pro:anticanonical map}, it admits both a canonical
  and an anticanonical map. By an argument used in the proof of the
  nerve theorem, these maps are $n$-homotopy equivalences and are
  $n$-homotopy inverses of each other.
\end{remark}

\begin{definition*}
  We say that a closed cover is
  \df{cover!$n$-semiregular}{$n$-semiregular} if it is a locally
  finite cover with at most $n$-dimensional nerve, which admits an
  anticanonical map that is an $n$-homotopy equivalence.
\end{definition*}

\begin{lemma}
  \label{lem:regular is semiregular}
  Every closed $n$-regular cover with at most $n$-dimensional nerve is
  $n$-semiregular.
\end{lemma}
\begin{proof}
  Apply proposition~\ref{pro:anticanonical map} and
  remark~\ref{rem:nerve theorem}.
\end{proof}

\begin{proposition}
  \label{pro:small cover with anticanonical map is semiregular}
  Every $ANE(n)$-space admits an open cover~$\U{}$ such that every
  anticanonical map of a closed locally finite locally finite
  dimensional cover that refines~$\U{}$ is an $n$-homotopy
  equivalence. In particular, if a closed locally finite cover
  refines~$\U{}$, has at most $n$-dimensional nerve and admits an
  anticanonical map, then it is $n$-semiregular.
\end{proposition}
\begin{proof}
  Let $X$ be an $ANE(n)$-space. By theorem~\ref{thm:ane
    characterization}, there exists an open cover $\U{}$ of~$X$ such
  that every two $\U{}$-close maps into $X$ are homotopic, provided
  that their domain is an at most $(n-1)$-dimensional metric space. We
  will show that $\U{}$ satisfies our claim.  Let $\F{}$ ba a closed
  locally finite locally finite dimensional cover of $X$ that
  refines~$\U{}$ and has an anticanonical map $\lambda \colon N(\F{})
  \to X$. By proposition~\ref{pro:canonical map}, there exisits a
  canonical map $\kappa \colon X \to N(\F{})$ of $\F{}$. By
  lemma~\ref{lem:regular stars} and corollary~\ref{cor:carried
    homotopy}, the composition $\kappa \circ \lambda$ is homotopic to
  the identity of $N(\F{})$. By the definitions, the composition
  $\lambda \circ \kappa$ is $\F{}$-close to the identity of $X$.
  Hence, for each map $f$ from an at most $(n-1)$-dimensional
  simplicial complex into $X$, the composition $\lambda \circ \kappa
  \circ f$ is $\U{}$-close to $f$, hence homotopic to $f$, by the
  choice of $\U{}$.  So, $\lambda \circ \kappa$ is $n$-homotopic to
  the identity on $X$.  Hence, $\kappa$ and $\lambda$ are $n$-homotopy
  inverses of each other.
\end{proof}

\begin{proposition}
  \label{pro:semiregular is bounded by X}
  If $\F{}$ is a closed $n$-semiregular cover of an at most
  $n$-dimensional $ANE(n)$-space $X$, then every map from an at most
  $(n-1)$-dimensional simplicial complex into an element of $\F{}$ is
  null-homotopic in~$X$.
\end{proposition}
\begin{proof}
  Let $f$ be a map from an at most $(n-1)$-dimensional simplicial
  complex into an element of $\F{}$. By proposition~\ref{pro:canonical
    map}, there exists a canonical map of~$\F{}$, which we denote
  by~$\kappa$. Let $\lambda$ be an anticanonical map of $\F{}$ that is
  an $n$-homotopy equivalence. By Whitehead's characterization, there
  exists a map $\mu \colon X \to N(\F{})$ such that $\lambda \circ
  \mu$ is $n$-homotopic to the identity on~$X$. The composition
  $\kappa \circ \lambda$ is carried by $B \circ B^{-1}$, so it is
  $\BF$-close to the identity on the nerve of $\F{}$ and by
  corollary~\ref{cor:carried homotopy}, it is homotopic to it. Hence,
  $\mu$ is homotopic to $\kappa \circ \lambda \circ \mu$, which is
  $n$-homotopic to~$\kappa$. The identity on $X$ is $n$-homotopic to
  $\lambda \circ \mu$, which is $n$-homotopic to $\lambda \circ
  \kappa$. Hence, by the definition of $n$-homotopy, $f$ is homotopic
  to $(\lambda \circ \kappa) \circ f$. But $\kappa \circ f$ is
  null-homotopic because its image lies entirely in a barycentric star
  of a vertex of the nerve of $\F{}$, by the definition of~$\kappa$
  and the assumption that the image of $f$ lies entirely in an element
  of~$\F{}$.  Thus~$f$, being homotopic to $\lambda \circ \kappa \circ
  f$, is null-homotopic.
\end{proof}

\section{A construction of regular covers}

We will now present a theorem (theorem~\ref{thm:pump up the regularity}) that
allows us to construct an $n$-regular cover from an $n$-semiregular cover. We
are still far away from its proof, which will be given in chapter~\ref{ch:the
  existence of regular covers}.

\begin{definition*}
  We say that a cover $\F{}$ is
  \df{cover!$n$-contractible}{$n$-contractible in a cover~$\G{}$} if
  it refines $\G{}$ and if every map of an at most $(n-1)$-dimensional
  simplicial complex into an element of~$\F{}$ is null-homotopic in an
  element of $\G{}$.
\end{definition*}

By a way of example, note that proposition~\ref{pro:semiregular is
  bounded by X} states that every closed $n$-semiregular cover of an
at most $n$-dimensional $ANE(n)$-space $X$ is $n$-contractible in~$\{
X \}$.

\begin{theorem}
  \label{thm:pump up the regularity} There exists a constant $N$ such
  that if a closed star-finite $k$-regular $n$-semiregular ($k <
  n$) interior $\N{n}$-cover $\F{}$ is $n$-contractible in a
  cover~$\E{}$, then there exists a closed $(k+1)$-regular
  $n$-semiregular interior $\N{n}$-cover isomorphic to $\F{}$ and
  refining $\st^N \E{}$. Moreover, we may require that the constructed
  cover is equal to $\F{}$ on a neighborhood of a given $Z(\F{})$-set.
\end{theorem}

The rest of this section is devoted to proofs of three corollaries of
theorem~\ref{thm:pump up the regularity} that will be used in the proof of the
main theorem. 

\begin{corollary}
  \label{cor:top level simplification} If $\F{}$ is a closed
  star-finite $n$-semiregular interior $\N{n}$-cover of a space~$X$
  and~$Z$ is a $Z(\F{})$-set, then there exists a closed $n$-regular
  interior $\N{n}$-cover of $X$ that is isomorphic to~$\F{}$ and is
  equal to~$\F{}$ on a neighborhood of $Z$.
\end{corollary}
\begin{proof}
  By proposition~\ref{pro:semiregular is bounded by X}, every closed
  $n$-semiregular cover of $X$ is $n$-contractible in~$\{ X \}$. Let
  $\F{0} = \F{}$. By theorem~\ref{thm:pump up the regularity} applied
  recursively with $\E{} = \{ X \}$, for each $0 < k < n$ there exists
  a closed star-finite $k$-regular $n$-semiregular interior
  $\N{n}$-cover $\F{k}$ of $X$ that is isomorphic to $\F{k-1}$ and is
  equal to $\F{k-1}$ on a neighborhood of $Z$. The cover $\F{n}$ is a
  cover that satisfies our claim.
\end{proof}

\begin{corollary}
  \label{cor:small cover simplification} For every cover $\U{}$ of an
  $\N{n}$-space $X$ there exists a cover $\V{}$ of $X$ such that if
  $\F{}$ is a closed star-finite $n$-semiregular interior
  $\N{n}$-cover that refines $\V{}$ and $Z$ is a $Z(\F{})$-set, then
  there exists a closed $n$-regular interior $\N{n}$-cover of $X$ that
  is isomorphic to $\F{}$, is equal to~$\F{}$ on an open neighborhood
  of $Z$ and refines $\U{}$.
\end{corollary}
\begin{proof}
  Let $N$ be the constant obtained via theorem~\ref{thm:pump up the
    regularity}. By theorem~\ref{thm:ane characterization}, for each
  open cover $\U{k}$ of $X$ there exists an open cover $\V{k}$ of $X$
  that is $n$-contractible in $\U{k}$. Define recursively a sequence
  $\U{n-1}, \V{n-1}, \ldots, \U{0}, \V{0}$ of open covers of $X$ such
  that $\st^N \U{n-1}$ refines $\U{}$, $\V{k}$ is $n$-contractible in
  $\U{k}$ and $\st^N \U{k-1}$ refines $\V{k}$. We will show that $\V{}
  = \V{0}$ satisfies our claim. Let $\F{0} = \F{}$. By the assumption,
  $\F{0}$ refines $\V{0}$, hence it is $n$-contractible in $\U{0}$.
  For $k = 1, 2, \ldots, n$, we apply theorem~\ref{thm:pump up the
    regularity} to a cover $\F{k-1}$ with $\E{} = \U{k-1}$ and obtain
  a closed $k$-regular $n$-semiregular interior $\N{n}$-cover $\F{k}$
  of $X$ that is isomorphic to $\F{k-1}$, refines $\st^N \U{k-1}$ (and
  hence $\V{k}$) and is equal to $\F{k-1}$ on a neighborhood of $Z$.
  The cover $\F{n}$ is a cover that satisfies our claim.
\end{proof}

\begin{definition*}
  We say that a cover is \df{cover!$m$-barycentric}{$m$-barycentric}
  if its nerve has a structure of an $m$th barycentric subdivision of
  a locally finite simplicial complex~$K$. We say that $K$ is a
  \df{primary complex}{primary complex} of the cover.
\end{definition*}

\begin{corollary}
  \label{cor:isomorphic cover simplification} For every $n$ and $k$
  there is a constant $m_{n,k}$ such that if $\E{}$ is an $n$-regular
  $m_{n,k}$-barycentric cover of an $\N{n}$-space $X$, $\F{}$ is a
  closed star-finite $n$-semiregular interior $\N{n}$-cover that
  refines $\st^k \E{}$ and $Z$ is a $Z(\F{})$-set, then there exists a
  closed $n$-regular interior $\N{n}$-cover of $X$ that is isomorphic
  to~$\F{}$, refines $\st^{m_{n,k}} \E{}$ and is equal to $\F{}$ on an
  open neighborhood of $Z$.
\end{corollary}

To prove corollary~\ref{cor:isomorphic cover simplification} we need
four lemmas stated below. The notion of a join is auxiliary and won't
be used later.

\begin{definition*}
  A \df{cover!$k$th join of a}{$k${\em th} join $\jn^k \F{}$} of an
  $m$-barycentric cover $\F{}$ with primary complex $K$ is a cover by
  those $(2^k-1)$th stars of elements of $\F{}$ that correspond to
  vertices of the $(m-k)$th barycentric subdivision of $K$.
\end{definition*}

\begin{remark}
  If a cover is $m$-barycentric, then its $k$th join is
  $(m-k)$-barycentric, when considered with the same primary complex.
  Also, every $m$-barycentric cover is $k$-barycentric for each $k <
  m$.  For each $(k+l)$-barycentric cover $\F{}$, we have $\jn^k \jn^l
  \F{} = \jn^{k+l} \F{}$.
\end{remark}

\begin{lemma}
  \label{lem:star prec join prec star}
  If $\F{}$ is $m$-barycentric, then $\st^{2^{m-1}-1} \F{} \prec \jn^{m} \F{}
  \subset \st^{2^m-1} \F{}$.
\end{lemma}
\begin{proof}
  The second relation follows directly from the definition. To prove
  the first relation, consider a discrete metric on the set $\F{}$,
  with the distance between two elements of~$\F{}$ defined to be the
  length of a shortest chain containing them, decreased by~$1$ (a
  chain is a sequence of sets in which every two consecutive elements
  intersect).  Then elements of $\st^{2^{m-1}-1} \F{}$ correspond to
  balls of radius $2^{m-1}-1$ in $\F{}$ and elements of $\jn^m \F{}$
  correspond to balls of radius $2^{m}-1$ in $\F{}$, but with centers
  in elements of~$\F{}$ corresponding to vertices of the primary
  complex $K$ of $\F{}$.  The nerve of $\F{}$ is $m$th barycentric
  subdivision of $K$, hence a distance of an element of $\F{}$ from an
  element that corresponds to a vertex of $K$ is not greater than
  $2^{m-1}$. We are done.
\end{proof}

\begin{definition*}
  We say that a subcomplex $L$ of a complex $K$ is \df{full subcomplex}{full}
  if it contains every simplex spanned in $K$ by its vertices.
\end{definition*}

\begin{lemma}
  \label{lem:intersections of bst-m}
  If $L_1, L_2, \ldots, L_k$ are subcomplexes of some other complex,
  then
   \[
   \bst^m L_1 \cap \ldots \cap \bst^m L_k = \bst^{m-1}(\bst
   L_1 \cap \ldots \cap \bst L_k),
   \]
   for each $m > 1$.
\end{lemma}
\begin{proof}
  Assume $L$ and $M$ are subcomplexes of some other complex. We shall
  prove that
  \em\begin{enumerate}
  \item If $L \cup M$ is full, then $\bst L \cap \bst M = \bst (L \cap
    M)$.
  \end{enumerate}\em
  We have to show that $\bst L \cap \bst M \subset \bst (L \cap M)$.
  Assume that $L \cup M$ is a full subcomplex of a complex~$K$.
  Observe that for each pair of vertices $v_1, v_2 \in K$ the
  interesction $\bst v_1 \cap \bst v_2$ is non-empty if and only if
  $v_1, v_2$ are connected by an edge. Let $v$ be a vertex in $\bst L
  \cap \bst M$. Then $v \in \bst v_1$ and $v \in \bst v_2$ for some
  $v_1 \in L$ and $v_2 \in M$. Then $v_1$ and $v_2$ are connected by
  an edge, which must be contained either in $L$ or in $M$, because $L
  \cup M$ is full. Hence either $v_1$ or $v_2$ is contained in $L \cap
  M$. Therefore $v$ is contained in $\bst (L \cap M)$ and we are done.

  A barycentric star of a subcomplex is always full, hence~(1) yields
  the equality
  \[
  \bst (\bst^{m-1} L_1) \cap \bst (\bst^{m-1} L_2 \cap \ldots \cap \bst^{m-1} L_k) =
  \bst (\bst^{m-1} L_1 \cap \ldots \bst^{m-1} L_k).
  \]
  An elementary induction on $m$ and $k$ finishes the proof (observe
  that $\bst(\bst^{m-1} L) = \bst^{m-1}(\bst L)$).
\end{proof}

\begin{lemma}
  \label{lem:bst homotopic}
  If $L$ is a full subcomplex of a locally finite dimensional simplicial
  complex~$K$, then $\bst L$ and $L$ are homotopy equivalent.
\end{lemma}


\begin{proof}
  Let $\mathcal{L}$ denote the collection of barycentric stars of
  vertices of $L$, taken in $K$. Since $L$ is a full subcomplex of
  $K$, $\mathcal{L}$ is isomorphic to the cover of $L$ by barycentric
  stars of its vertices, taken in $L$. Hence, by
  lemma~\ref{lem:invertible carriers}, the nerve of $\mathcal{L}$ is
  equal to $L$. By lemma~\ref{lem:regular stars} and by the definition
  of $\bst L$, $\mathcal{L}$ is a regular cover of $\bst L$. By
  theorem~\ref{thm:general nerve theorem}, $\bst L$ is homotopy
  equivalent to the nerve of $\mathcal{L}$. We are done.
\end{proof}

\begin{lemma}
  \label{lem:when join is regular} An $m$th join a closed $n$-regular
  $m$-barycentric cover of an at most $n$-dimensional metric space is
  $n$-regular.
\end{lemma}
\begin{proof}
  Let $\F{}$ be a closed $n$-regular $m$-barycentric cover of an at
  most $n$-dimensional metric space $X$. Let $\BF$ be a cover of the
  nerve of $\F{}$ by barycentric stars of its vertices. Let $\BFn$ be
  the restriction of $\BF$ to the $n$-dimensional skeleton of the
  first barycentric subdivision of the nerve of $\F{}$. We have to
  show that non-empty intersections of elements of $\jn^m \F{}$ are
  absolute extensors in dimension $n$. Fix any such intersection~$A$.
  Let $B$ be the intersection of the corresponding elements of $\jn^m
  \BFn$.  As we argued in the proof of the nerve theorem, $\BFn$ is
  $n$-regular and isomorphic to $\F{}$. Clearly, the restrictions of
  $\F{}$ to $A$ and of $\BFn$ to $B$ are $n$-regular and isomorphic.
  Hence, by the nerve theorem, $A$ is $n$-homotopy equivalent to $B$.
  By theorem~\ref{thm:dugundji characterization}, being an absolute
  extensor in dimension $n$ is $n$-homotopy invariant in the class of
  absolute neighborhood extensors in dimension $n$. Hence it suffices
  to show that $B$ is an absolute extensor in dimension $n$, which in
  turn would follow from the cellular approximation theorem
  (see~\cite[Corollary 4.12]{hatcher2002}) if the intersections of
  elements of $\jn^m \BF$ were contractible (cf.  proof of
  lemma~\ref{lem:regular stars}), which we are about to prove. Let $K$
  be a primary complex of $\BF{}$ such that the nerve of $\BF$ is the
  $m$th barycentric subdivision of $K$. Observe that $\jn^m \BF{}$ is
  a collection of $m$th barycentric stars of vertices of $K$. By
  lemma~\ref{lem:intersections of bst-m}, if $\{ v_i \}$ is a collection of vertices of
  $K$, then $\bigcap \bst^m v_i = \bst^{m-1} (\bigcap \bst v_i)$. By
  lemma~\ref{lem:bst homotopic}, $\bst^{m-1} (\bigcap \bst v_i)$ is
  homotopy equivalent to $\bigcap \bst v_i$.  The latter intersection
  is contractible by lemma~\ref{lem:invertible carriers}.  We are
  done.
\end{proof}

\begin{proof}[Proof of corollary~\ref{cor:isomorphic cover simplification}]
  By the definitions, if a cover refines an $n$-regular cover, then it
  is $n$-contractible in it. By lemma~\ref{lem:star prec join prec
    star} and lemma~\ref{lem:when join is regular}, the proof is
  reduced to an easy recursive application of theorem~\ref{thm:pump up
    the regularity} and we omit its details.
\end{proof}

The following technical lemma will be used in chapter~\ref{ch:proof of
  the main theorem}.

\begin{lemma}
  \label{lem:map carried into star}
  Assume that $\G{}$ is a $2$-barycentric and $n$-regular cover of a
  space~$Y$ and~$\F{}$ is a closed locally finite locally finite
  dimensional cover of an at most $n$-dimensional metric space $X$. If
  $\F{}$ and $\G{}$ are isomorphic, then every map $f$ defined on a
  closed subset of $X$ that maps elements of $\F{}$ into the
  corresponding elements of $\st \G{}$ extends over the entire domain
  to a map that maps elements of $\F{}$ into the corresponding
  elements of $\st^{7} \G{}$.
\end{lemma}
\begin{proof}
  Let $\F{} = \{ F_i \}_{i \in I}$ and $\G{} = \{ G_i \}_{i \in I}$.
  By lemma~\ref{lem:star prec join prec star} there exists a map $q$
  from the indexing set of $\st \G{}$ into the indexing set of $\jn^2
  \G{}$ such that $\st_\G{} G_i \subset \st^3 G_{q(i)}$. The function
  $F_i \mapsto \st^3 G_{q(i)}$ is a carrier and $f$ is carried by it.
  By lemma~\ref{lem:when join is regular}, $\jn^2 \G{}$ is
  $n$-regular, hence by the carrier theorem $f$ extends over the
  entire domain to a map $g$ such that $g(F_i) \subset \st^3
  G_{q(i)}$. Since $G_i \subset \st^3 G_{q(i)}$, we have $G_{q(i)}
  \subset \st^4 G_i$. Hence $\st^3 G_{q(i)} = \st^3_\G{} G_{q(i)}
  \subset \st^3_\G{} \st^4 G_i = \st^7 G_i$.  Hence $g$ maps elements
  of $\F{}$ into the corresponding elements of $\st^7 \G{}$.
\end{proof}

\section{A construction of semiregular covers}

First we briefly describe the technique of the construction (cf.
figure~\ref{fig:semiregular}). Consider an $\N{n}$-space $X$. By
lemma~\ref{lem:equivalence with a complex}, there exists a weak
$n$-homotopy equivalence of an at most $n$-dimensional countable
locally finite simplicial complex $L$ and $X$. By the strong
univerality of $X$, we may approximate it by a closed embedding. By
lemma~\ref{lem:set of weak equivalences is open}, we may require (and
we do) that the approximation is also a weak $n$-homotopy equivalence.
Let $K$ denote its image. It is a locally finite simplicial complex
embedded as a closed subset of $X$ and the inclusion $K \subset X$ is
a weak $n$-homotopy equivalence. By theorem~\ref{thm:subspace
  retraction}, there exists a retraction $r \colon X \to K$. Let
$\mathcal{B}_K$ be the cover of $K$ by barycentric stars of its
vertices. The nerve of $\mathcal{B}_K$ is isomorphic to $K$. Let $\F{}
= r^{-1}(\mathcal{B}_K)$.  It is isomorphic to $\mathcal{B}_K$, hence
its nerve is isomorphic to $K$.  The identity on $K$ is an
anticanonical map of $\F{}$. It is a weak $n$-homotopy equivalence by
the construction. By Whitehead's characterization, it is an
$n$-homotopy equivalence, hence $\F{}$ is $n$-semiregular.

{%
\begin{figure}[ht]
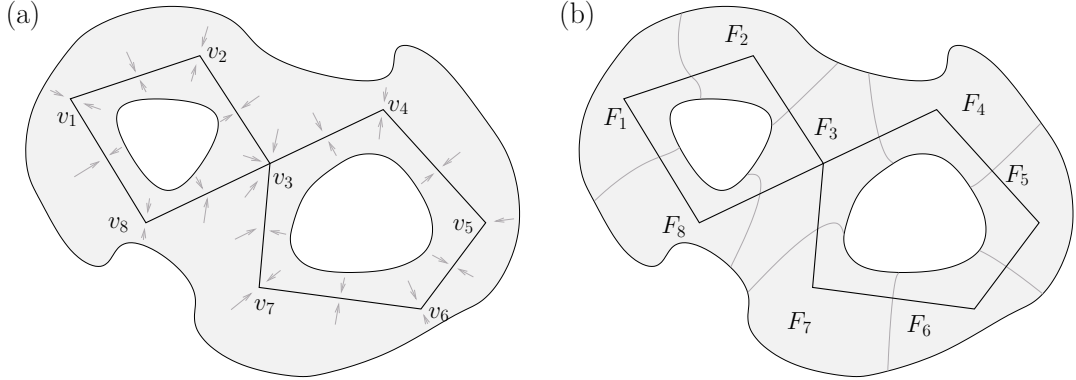

\begin{center}
$\begin{array}{c@{\hspace{5mm}}c}
\resizebox{0.45\textwidth}{!}{\input{psfigures/semiregular.pstex_t}} &
\resizebox{0.45\textwidth}{!}{\input{psfigures/semiregular2.pstex_t}}
\end{array}$
\end{center}
\caption{{Construction of a semiregular
cover.}}
\label{fig:semiregular}
\end{figure}
}

The construction described above allows us to construct a closed
$n$-semiregular cover of an $\N{n}$-space, but it doesn't give us a
control over the size of the constructed cover. It is because
theorem~\ref{thm:subspace retraction} does not give a way to limit the
size of fibers of $r$. Theorem~\ref{thm:retraction onto a complex}
allows us to do just that.

\begin{theorem}
  \label{thm:retraction onto a complex}
  For every open cover $\U{}$ of an $\N{n}$-space $X$ there exists a locally
  finite simplicial complex $K$ embedded as a closed subset of $X$ and a
  retraction of $X$ onto $K$ that is $\U{}$-close to the identity.
\end{theorem}

The proof of theorem~\ref{thm:retraction onto a complex} will be given in
section~\ref{sec:retraction onto a core}.

\begin{theorem}
  \label{thm:existence of semiregular covers} For each integer $m$,
  each $\N{n}$-space $X$ and each $Z$-set $Z$ in $X$ there exist
  arbitrarily small closed star-finite $m$-barycentric $n$-semiregular
  interior $\N{n}$-covers $\F{}$ of $X$ such that $Z$ is a
  $Z(\F{})$-set.
\end{theorem}
\begin{proof}
  Let $\U{}$ be an open cover of $X$. By theorem~\ref{thm:retraction
    onto a complex}, there exists a locally finite simplicial complex
  $K$ embedded as a closed subset of~$X$ and a retraction $r \colon X
  \to K$ that is $\U{}$-close to the identity. We may assume that the
  set of simplices of $K$ refines~$\U{}$ and that $K$ has a structure
  of an $m$th barycentric subdivision of another complex. Let $\F{} =
  r^{-1}(\mathcal{B}_K)$ and observe that $\F{}$ is a closed locally
  finite star-finite $m$-barycentric cover that refines $\st \U{}$. By
  theorem~\ref{thm:n-swelling}, there is a swelling~$\G{}$ of $\F{}$
  to a closed interior $\N{n}$-cover such that $Z$ is a $Z(\G{})$-set
  and $\G{} \prec \st \U{}$. As $K$ is the nerve of $\G{}$, the
  inclusion of $K$ into $X$ is an anticanonical map of $\G{}$. By
  proposition~\ref{pro:small cover with anticanonical map is
    semiregular}, if $\U{}$ is small enough, then $\G{}$ is
  $n$-semiregular. We are done.
\end{proof}


\chapter[The homeomorphism extension theorem]{Extending homeomorphisms by the use of a ``brick partitionings'' technique}
\label{ch:the homeomorphism extension theorem}

The classical approach to the problem of constructing a homeomorphism between
complete spaces is the brick partitionings technique.
Theorem~\ref{thm:homeomorphism extension} enhances this technique and sets the
strategy for the proof of the main theorem.

\begin{definition*}
  We write \df{F1F2@$\F{1} \prec_f \F{2}$}{$\F{1} \prec_f \F{2}$} if $\F{1}$
  and $\F{2}$ are covers of the same space and $f$ is a map from the indexing
  set of $\F{1}$ into the indexing set of $\F{2}$ such that for all $i$ the
  element of $\F{1}$ indexed by $i$ is a subset of the element of $\F{2}$
  indexed by $f(i)$.
\end{definition*}

\begin{theorem}[The homeomorphism extension theorem]
  \label{thm:homeomorphism extension}

  Assume that $\F{k}$ and $\G{k}$ are sequences of locally finite
  closed covers of complete spaces~$X$ and~$Y$. We do not assume that
  $X$ or $Y$ is separable or that $\F{k}$ or $\G{k}$ is countable.
  Assume that~$\F{k}$ and~$\G{k}$ are isomorphic for all~$k$ and $f_k$
  is a sequence of maps such that for some~$m$ and each $k$
  \begin{center}
    \begin{tabular}{r @{$\ \prec_{f_{2k+1}}\ $} c @{$\ \prec_{f_{2k}}\ $ } l }
      $\F{2k+1}$       & $\st^m \F{2k}$ & $\F{2k-1}$,  \\
      $\st^m \G{2k+1}$ & $\G{2k}$       & $\st^m \G{2k-1}$. \\
    \end{tabular}
  \end{center}
  If meshes of $\F{k}$ and of $\G{k}$ converge to zero, then there
  exists a unique homeomorphism from~$X$ onto~$Y$ that maps each
  element of $\F{k}$ into the corresponding element of $\st^{m+1}
  \G{k}$.
\end{theorem}

The name of \emph{the homeomorphism extension theorem} deserves an explanation
as the theorem does not mention explicitly extensions of homeomorphisms. It is
justified by the fact that the uniqueness of the homeomorphism implies that it
extends every partial homeomorphism that maps elements of~$\F{k}$'s into the
corresponding elements of~$\st^{m+1} \G{k}$'s.

\begin{proof}
  The uniqueness of a homeomorphism from $X$ onto $Y$ that maps elements
  of~$\F{k}$ into the corresponding elements of $\st^{m+1} \G{k}$ is a direct
  consequence of the assumption that meshes of $\G{k}$ converge to zero. Every
  two maps that satisfy this condition must be $\st^{m+1}\G{k}$-close for all
  $k$.

  The existence is proved in two steps. Let $\F{k} = \{ F^k_i \}_{i
    \in I_k}$ and $\G{k} = \{ G^k_i \}_{i \in I_k}$. In the first step
  we prove that there exists a function $f \colon X \to Y$ such that
  $f(F^k_i) \subset \st^{m+1} G^k_i$ for all $k \in \mathbb{N}$ and $i
  \in I_k$. In the second step we prove that it is a homeomorphism. 

  By the definition, an $m$th star of an element of a collection of
  sets is a union of some of the elements of the collection. We shall
  call summands of this union \emph{components} of the $m$th star.
  Note that the set of components of an $m$th star may differ from the
  set of elements of the collection that are subsets of the $m$th
  star. If $\mathcal{C}$ is a collection of sets and if $C_1, C_2$ are
  elements of $\mathcal{C}$, then $C_1$ is a component of
  $\st^m_\mathcal{C} C_2$ if and only if $C_2$ is a component of
  $\st^m_\mathcal{C} C_1$.

  For each $x \in X$ let $A_k(x) = \bigcup \{ \st^m G^k_i \colon i \in I_k, x
  \in F^k_i \}$. We first wish to apply the Cantor theorem to prove that
  $\bigcap_k A_k(x) \neq \emptyset$ for each $x \in X$.
  By local finiteness of $\G{k}$'s, each $A_k(x)$ is closed. The
  sequence of diameters of $A_k(x)$ converges to zero as $k$ goes to
  infinity since $A_k(x)$ is an union of $m$th stars of intersecting
  elements of $\G{k}$. The proof of the inclusion $A_{k+1}(x) \subset
  A_k(x)$ depends on the parity of $k$. From the definition of
  $A_{k+1}(x)$ it suffices to check that for each $i \in I_{k+1}$ such
  that $x \in F^{k+1}_i$ the inclusion $\st^m G^{k+1}_i \subset
  A_k(x)$ holds. 

  Assume that $k$ is even. Then $F^{k+1}_i \subset \st^m
  F^k_{f_{k+1}(i)}$. Then for some $j \in I_k$ such that $F^k_j$ is a
  component of $\st^m F^k_{f_{k+1}(i)}$, we have $x \in F^k_j$. Then
  from the definition $\st^m G^k_j \subset A_k(x)$. But $\F{k}$ and
  $\G{k}$ are isomorphic, so $G^k_{f_{k+1}(i)}$ is a component of
  $\st^m G^k_j$. In particular, $G^k_{f_{k+1}(i)} \subset \st^m
  G^k_j$.  Hence $\st^m G^{k+1}_i \subset G^k_{f_{k+1}(i)} \subset
  \st^m G^k_j \subset A_k(x)$, the first inclusion from the assumption
  that $k$ is even.

  Assume that $k$ is odd. Let $j \in I_{k+1}$ such that $F^{k+1}_j$ is
  a component of $\st^m F^{k+1}_i$. Then $F^{k+1}_i$ is a component of
  $\st^m F^{k+1}_j$. In particular, $F^{k+1}_i \subset \st^m
  F^{k+1}_j$. Hence $x \in F^{k+1}_i \subset \st^m F^{k+1}_j \subset
  F^k_{f_{k+1}(j)}$, the last inclusion from the assumption that $k$
  is odd. By the definition, $\st^m G^k_{f_{k+1}(j)} \subset A_k(x)$
  and by the assumptions $G^{k+1}_j \subset \st^m G^k_{f_{k+1}(j)}$.
  Therefore $G_j^{k+1} \subset A_k(x)$ so $\st^m G^{k+1}_i \subset
  A_k(x)$.

  Thus Cantor's theorem applies to give a unique function $x \mapsto
  f(x) \in \bigcap_k A_k(x)$ such that if $x \in F^k_i$ then $f(x) \in
  A_k(x) \subset \st^{m+1} G^k_i$ for all $k \in \mathbb{N}$ and $i
  \in I_k$. Let $x \in X$ and an open neighborhood $U$ of $f(x)$ in
  $Y$ be given. Take $k$ so large that the star of $f(x)$ in
  $\st^{m+1} \G{k}$ lies in $U$.  By local finiteness of $\F{k}$,
  there is open neighborhood $V$ of $x$ in $X$ that lies in a star of
  $x$ in a cover $\F{k}$. Then $f(V) \subset U$, showing that $f$ is
  continuous. By the symmetry of assumptions (which we obtain after
  throwing away $\F{1}$ and $\G{1}$ from the sequences), there exists
  a continuous map $g \colon Y \to X$ such that $g(G^k_i) \subset
  \st^{m+1} F^k_i$. By these conditions, $g(f(F^k_i)) \subset
  \st^{2m+2} F^k_i$ so $g \circ f = id_X$ because it is $\st^{2m+2}
  \F{k}$-close to $id_X$ for all $k$. Similarly $f \circ g = id_Y$, so
  $f$ is a homeomorphism from $X$ onto $Y$.
\end{proof}


\renewcommand{\P}[1]{{\protect\mathcal{P}_{#1}}}
\newcommand{\Q}[1]{{\protect\mathcal{Q}_{#1}}}

\chapter{Proof of the main results}
\label{ch:proof of the main theorem}

In this chapter we prove the main theorem and derive all theorems
stated in the introduction. We use theorems~\ref{thm:pump up the
  regularity} and~\ref{thm:retraction onto a complex}, whose proofs
will be given in the third part of the paper.

\begin{definition*}
  We say that a cover is a \df{closed partition}{closed partition of
    an $\N{n}$-space $X$} if it is a closed star-finite $n$-regular
  interior $\N{n}$-cover of $X$ with at most $n$-dimensional nerve.
\end{definition*}

\begin{lemma}
  \label{lem:existence of small partitions}
  For each open cover $\U{}$ of an $\N{n}$-space $X$ there exists a
  closed partition $\F{}$ of $X$ that refines $\U{}$. Moreover, if $m$
  is an integer and $Z$ is a $Z$-set in $X$, then we may require that
  $\F{}$ is $m$-barycentric and $Z$ is a $Z(\F{})$-set.
\end{lemma}
\begin{proof}
  Let $\V{}$ be an open cover of $X$ obtained via
  corollary~\ref{cor:small cover simplification} applied to $\U{}$. By
  theorem~\ref{thm:existence of semiregular covers} there exists a
  closed star-finite $m$-barycentric $n$-semiregular interior
  $\N{n}$-cover $\E{}$ that refines $\V{}$ and such that $Z$ is a
  $Z(\E{})$-set. Let $\F{}$ be a closed $n$-regular interior
  $\N{n}$-cover of~$X$ that refines~$\U{}$ and that is isomorphic to
  $\E{}$ and equal to $\E{}$ on an open neighborhood of $Z$.  The
  existence of~$\F{}$ follows again from corollary~\ref{cor:small
    cover simplification}, this time applied to $\E{}$, $Z$, $\V{}$
  and $\U{}$. By the construction, $\F{}$ is $m$-barycentric and $Z$ is
  a $Z(\F{})$-set. We are done.
\end{proof}

\begin{definition*}
  If $\F{}$ and $\G{}$ are isomorphic covers of spaces $X$ and~$Y$
  respectively and $h$ is a homeomorphism from a subset of $X$ onto a
  subset of $Y$, then we say that \df{covers!compatible with a
    homeomorphism}{$\F{}$ and $\G{}$ are compatible with $h$} if $h$
  maps elements of $\F{}$ into the corresponding elements of $\st
  \G{}$ and $h^{-1}$ maps elements of $\G{}$ into the corresponding
  elements of $\st \F{}$.
\end{definition*}

\begin{lemma}
  \label{lem:existence of isomorphic partition}
  If $\F{}$ is a closed partition of an $\N{n}$-space $X$ and $Y$ is
  an $\N{n}$-space $n$-homotopic to $X$, then there exists a closed
  parition $\G{}$ of $Y$ that is isomorphic to $\F{}$. Moreover, if an
  $n$-homotopy equivalence $f \colon X \to Y$ is a $Z$-embedding on a
  closed subset $A$ of $X$, then we may require that $\F{}$ and $\G{}$
  are compatible with $f_{|A}$ and that $f(A)$ is a $Z(\G{})$-set.
\end{lemma}
\begin{proof}
  By theorem~\ref{thm:ane characterization}, if a map into an
  $ANE(n)$-space is close enough to an $n$-homotopy equivalence, then
  it is an $n$-homotopy equivalence. Hence, by
  theorem~\ref{thm:z-approximation}, we may approximate $f$ $\rel A$
  by a closed embedding $h$ that is an $n$-homotopy equivalence. By
  theorem~\ref{thm:subspace retraction} there exists a retraction $r
  \colon Y \to h(X)$. Let $\F{} = \{ F_i \}_{i \in I}$ and let
  $\mathcal{R} = \{ r^{-1}(h(F_i)) \}_{i \in I}$. The collection
  $\mathcal{R}$ is a closed locally finite cover of $Y$. It is
  isomorphic to $\F{}$, hence it is star-finite and has at most
  $n$-dimensional nerve.  Covers $\F{}$ and $\mathcal{R}$ are
  compatible with $f_{|A}$ by the definition of $r$ and $h$. Let
  $\lambda$ be an anticanonical map of $\F{}$ and an $n$-homotopy
  equivalence of $N(\F{})$ and $X$. The map $f \circ \lambda$ is an
  anticanonical map of $\mathcal{R}$ and an $n$-homotopy equivalence
  of $N(\mathcal{R})$ and $Y$. Hence $\mathcal{R}$ is $n$-semiregular.
  By theorem~\ref{thm:n-swelling} there is a closed interior
  $\N{n}$-cover $\mathcal{H}$ that is a swelling of $\mathcal{R}$ and
  such that $h(A)$ is a $Z(\mathcal{H})$-set. We may assume that the
  swelling is so small that $\F{}$ and $\mathcal{H}$ are compatible
  with $f_{|A}$. The cover $\mathcal{H}$ is isomorphic to
  $\mathcal{R}$, hence it is star-finite and has at most
  $n$-dimensional nerve.  As a swelling of an $n$-semiregular cover,
  $\mathcal{H}$ admits an anticanonical map that is an $n$-homotopy
  equivalence.  Hence $\mathcal{H}$ is $n$-semiregular. We apply
  corollary~\ref{cor:top level simplification} to $\mathcal{H}$ and to
  the $Z(\mathcal{H})$-set $h(A)$ to obtain a closed $n$-regular
  interior $\N{n}$-cover $\G{}$ isomorphic to $\mathcal{H}$ and equal
  to $\mathcal{H}$ on an neighborhood of $h(A)$. This cover satisfies
  all prescribed conditions.
\end{proof}

\begin{lemma}
  \label{lem:existence of isomorphic cover}
  For each $n$ there exists an integer $m$ such that if~$\F{}$
  and~$\G{}$ are isomorphic closed $m$-barycentric partitions of
  $\N{n}$-spaces~$X$ and~$Y$ respectively and~$\P{}$ is a closed
  partition of~$X$ such that $\P{} \prec_p \F{}$ (where $p$ is a map
  from the indexing set of $\P{}$ into the indexing set of $\F{}$),
  then there exists a closed partition $\Q{}$ of $Y$ that is
  isomorphic to $\F{}$ and such that~$\Q{} \prec_p \st^m \G{}$.
  Moreover, if $h$ is a $Z(\G{})$-embedding of a closed subset $A$ of
  $X$ into $Y$ such that $\F{}$ and $\G{}$ are compatible with $h$,
  then we may require that $\P{}$ and $\Q{}$ are compatible with $h$
  as well.
\end{lemma}
\begin{proof}
  We let $m$ to be equal to $m_{n,9} + 1$ or $2$, whichever is
  greater, with $m_{n,9}$ obtained via corollary~\ref{cor:isomorphic
    cover simplification}. Let $\F{} = \{ F_i \}_{i \in I}$ and $\G{}
  = \{ G_i \}_{i \in I}$.

  By the carrier theorem, there is a map $f \colon X \to Y$ that maps
  elements of $\F{}$ into the corresponding elements of $\G{}$.  The
  composition of $f$ with a canonical map of~$\G{}$ is a canonical map
  of~$\F{}$. By the nerve theorem, canonical maps of $\F{}$ and
  of~$\G{}$ are $n$-homotopy equivalences. Hence $f$ is an
  $n$-homotopy equivalence as well. By lemma~\ref{lem:set of weak
    equivalences is open} and by Whitehead's characterization, every
  approximation of $f$ that is sufficiently close is an $n$-homotopy
  equivalence. By the definition, every approximation of $f$ within
  $\G{}$ maps elements of $\F{}$ into the corresponding elements of
  $\G{}$.  Hence by corollary~\ref{cor:approximation within a cover}
  we may assume that $f$ is a closed embedding that has image disjoint
  from $h(A)$. The function $h^{-1} \cup f^{-1}$ is well defined and
  continuous and maps elements of $\G{}$ into the corresponding
  elements of $\st \F{}$. By lemma~\ref{lem:map carried into star},
  $h^{-1} \cup f^{-1}$ extends over $Y$ to a map $g$ that maps
  elements of $\G{}$ into the corresponding elements of $\st^7 \F{}$.
  Note that $g^{-1}(F_i) \subset \st^{8} G_i$.

  Let $\P{} = \{ P_j \}_{j \in J}$. Let $\mathcal{R} = \{ g^{-1}(P_j)
  \}_{j \in J}$. We have $\mathcal{R} \prec_p \st^{8} \G{}$. The
  collection $\mathcal{R}$ is a closed locally finite cover of $Y$. It
  is isomorphic to $\P{}$, hence it is star-finite and has at most
  $n$-dimensional nerve. Covers $\P{}$ and $\mathcal{R}$ are
  compatible with~$h$ by the definition of $g$. Let $\lambda$ be an
  anticanonical map of $\P{}$ and an $n$-homotopy equivalence of
  $N(\P{})$ and $X$. The map $f \circ \lambda$ is an anticanonical map
  of $\mathcal{R}$ and an $n$-homotopy equivalence of $N(\mathcal{R})$
  and $Y$. Hence $\mathcal{R}$ is $n$-semiregular.
  By~\ref{thm:n-swelling} there is a closed, locally finite interior
  $\N{n}$-cover $\mathcal{S}$ that is a swelling of the cover
  $\mathcal{R}$ such that $h(A)$ is a $Z(\mathcal{S})$-set.  If the
  swelling is small enough, then $\mathcal{S} \prec \st^{9} \G{}$ and
  covers $\P{}$ and $\mathcal{S}$ are compatible with~$h$. Let $Z$ be
  a $Z(\mathcal{S})$-set that contains $h(A)$ and meets every
  non-empty intersection of elements of $\mathcal{S}$ (we may let $Z$
  to be the union of $h(A)$ and $\im f \circ \lambda$, by
  remark~\ref{rem:closed locally compact is ZF} and by
  proposition~\ref{pro:closed ZFsigma is ZF} this set is a
  $Z(\mathcal{S})$-set). We apply corollary~\ref{cor:isomorphic cover
    simplification} to $\mathcal{S}$, $Z(\mathcal{S})$-set $Z$ and
  $\E{} = \G{}$ to obtain a closed $n$-regular $\N{n}$-cover $\Q{}$
  that refines $\st^{m_{n,9}}$, is isomorphic to $\mathcal{S}$ and is
  equal to $\mathcal{S}$ on an neighborhood of $Z$. We have $\Q{}
  \prec \st^{m_{n,9}} \G{}$. By the assumption that $Z$ meets every
  non-empty intersection of elements of $\mathcal{S}$, we have $\Q{}
  \prec_p \st^{m_{n,9} + 1} \G{}$. The cover $\Q{}$ is a closed
  partition of $Y$. Since $\Q{}$ is equal to $\mathcal{S}$ on $h(A)$,
  partitions $\P{}$ and $\Q{}$ are compatible with $h$.  We are done.
\end{proof}

\begin{lemma}
  \label{lem:existence of a homeomorphism}
  For each $n$ there exists a constant $m$ such that if $\F{}$ is an
  $m$-barycentric closed partition of an $n$-dimensional N\"obeling
  manifold $X$ and $\G{}$ is an isomorphic closed partition of an
  $n$-dimensional N\"obeling manifold $Y$, then there exists a
  homeomorphism from $X$ onto $Y$ that maps elements of $\F{}$ into
  the corresponding elements of $\st^m \G{}$. Moreover, if $\F{}$ and
  $\G{}$ are compatible with a homeomorphism $h$ of a $Z$-set in $X$
  onto a $Z(\G{})$-set in $Y$, then we may require that the ambient
  homeomorphism extends $h$.
\end{lemma}
\begin{proof}
  Let $m$ be a constant obtained via lemma~\ref{lem:existence of
    isomorphic cover} applied to~$n$.  We will construct a sequence of
  covers and functions $\G{1}, \F{1}, \F{2}, f_2, \G{2}, \G{3}, f_3,
  \F{3}, \ldots$ (in this particular order), where $f_k$ is a map from
  the indexing set of $\F{k}$ (and of $\G{k}$) into the indexing set
  of $\F{k-1}$ (and of $\G{k-1}$), such that the following conditions
  are satisfied for each $k$.

  \begin{enumerate}
  \item $\F{k}$ and $\G{k}$ are isomorphic closed $m$-barycentric
    partitions of $X$ and $Y$ respectively.
  \item $Z_1$ is a $Z(\F{2k})$-set and $Z_2$ is a $Z(\G{2k-1})$-set.
  \item $\mesh \F{2k} < \frac{1}{k}$ and $\mesh \G{2k+1} <
    \frac{1}{k+1}$.
  \item
    \begin{tabular}{r @{$\ \prec_{f_{2k}}\ $ } l r @{$\
          \prec_{f_{2k+1}}\ $} l }
      $\st^{m} \F{2k}$ & $\F{2k-1}$, & $\F{2k+1}$ & $\st^{m} \F{2k}$ \\
      $\G{2k}$ & $\st^{m} \G{2k-1}$, & $\st^{m} \G{2k+1}$ & $\G{2k}$. \\
    \end{tabular}
  \item $\F{k}$ and $\G{k}$ are compatible with $h$.
  \end{enumerate}

  By corollary~\ref{cor:inverse equivalence}, there exists a closed
  embedding $h_1 \colon X \to Y$ and a closed embedding $h_2 \colon Y
  \to X$ such that $h_{1|Z_1} = h$, $h_{2|Z_2} = h^{-1}$ and both
  $h_1$ and $h_2$ are $n$-homotopy equivalences.

  We let $\G{1} = \G{}$ and $\F{1} = \F{}$. Let $k \geq 1$ and assume
  that we already constructed $\G{2k-1}$ and $\F{2k-1}$.  Let $\U{}$
  be an open cover of $X$ with mesh less than $1/k$ and such that
  $\st^{m} \U{}$ refines $\F{2k-1}$. By lemma~\ref{lem:existence of
    small partitions}, there exists a closed $m$-barycentric partition
  $\F{2k}$ of $X$ that refines $\U{}$ and such that $Z_1$ is a
  $Z(\F{2k})$-set. By the choice of $\U{}$, the mesh of $\F{2k}$ is
  less than $1/k$ and there exists a function $f_{2k}$ such that
  $\F{2k} \prec_{f_{2k}} \F{2k-1}$. We let $f_{2k}$ to be any such
  function. By the choice of $m$ and by lemma~\ref{lem:existence of
    isomorphic cover} applied with $\F{} = \F{2k-1}$, $\G{} =
  \G{2k-1}$, $\P{} = \F{2k}$ and $p = f_{2k}$, there exists a closed
  partition $\G{2k}$ of $Y$ that is isomorphic to $\F{2k}$, satisfies
  the identity $\st^m \G{2k} \prec_{f_{2k}} \G{2k-1}$ and such that
  $\F{2k}$ and $\G{2k}$ are compatible with $h$.

  The construction of $\G{2k+1}$ and $f_{2k+1}$ is similar to the
  construction of $\F{2k}$ and $f_{2k}$ and the construction of
  $\F{2k+1}$ is similar to the construction of $\G{2k}$.

  We constructed a sequence $\G{1}, \F{1}, \F{2}, f_2, \ldots$ that
  satisfies conditions (1) - (5) stated above. Hence by
  theorem~\ref{thm:homeomorphism extension}, there exists a
  homeomorphism from $X$ onto $Y$ that maps elements of $\F{k}$ into
  the corresponding elements of $\st^{m+1} \G{k}$. By (5), the
  restriction of this homeomorphism to $Z_1$ is $\st^{m+1}
  \G{k}$-close to $h$ for each $k$. By (3), meshes of $\st^{m+1}
  \G{k}$ converge to $0$, so this restriction must be equal to $h$,
  hence the constructed homeomorphism extends $h$. We are done.
\end{proof}

Now we are set to prove the main theorem. Let us recall the statement.
\begin{main theorem}
  A homeomorphism $h \colon Z_1 \to Z_2$ of a $Z$-set $Z_1$ in an
  $\N{n}$-space $X$ onto a $Z$-set $Z_2$ in an $\N{n}$-space $Y$
  extends to an ambient homeomorphism if and only if it extends to
  an $n$-homotopy equivalence of $X$ and $Y$.
\end{main theorem}
\begin{proof}
  \index{proof!main theorem} Only the ``if'' part requires a proof.
  Let $m$ be a constant obtained via lemma~\ref{lem:existence of a
    homeomorphism} applied to $n$. Let $h_2 \colon X \to Y$ be an
  extension of $h^{-1}$ to an $n$-homotopy equivalence. It exists by
  corollary~\ref{cor:inverse equivalence}.  By
  lemma~\ref{lem:existence of small partitions}, there exists a closed
  $m$-barycentric partition $\G{}$ of $Y$ such that $Z_2$ is a
  $Z(\G{})$-set. By lemma~\ref{lem:existence of isomorphic partition}
  applied to $\G{}$ and an $n$-homotopy equivalence $h_2$, there
  exists a closed partition $\F{}$ of $X$ that is isomorphic to $\G{}$
  and such that $Z_1 = h_2(Z_2)$ is a $Z(\F{})$-set and $\F{}$ and
  $\G{}$ are compatible with $h$. An application of
  lemma~\ref{lem:existence of a homeomorphism} finishes the proof.
\end{proof}

The proof of the local $Z$-set unknotting is similar. Let us recall
the statement.

\begin{local z-set unknotting theorem}
  For every open cover $\U{}$ of a N\"obeling manifold $X$ there
  exists an open cover $\V{}$ such that every homeomorphism $h \colon
  Z_1 \to Z_2$ between $Z$-subsets of the manifold that is
  $\V{}$-close to the inclusion $Z_1 \subset X$ extends to a
  homeomorphism of the entire manifold that is $\U{}$-close to the
  identity.
\end{local z-set unknotting theorem}
\begin{proof}
  \index{proof!local $Z$-set unknotting theorem} Let $\U{}$ be an open
  cover of an $n$-dimensional N\"obeling manifold~$X$.  Let $m$ be a
  constant obtained via lemma~\ref{lem:existence of a homeomorphism}
  applied to $n$.  Let $\mathcal{W}$ be an open cover of $X$ whose
  $(2m+1)$-th star refines $\U{}$. By lemma~\ref{lem:existence of
    small partitions}, there exists a closed $m$-barycentric partition
  $\G{}$ of $X$ that refines $\W{}$. We let $\V{}$ be an open cover
  whose star refines $\G{}$. We will show that it satisfies the
  assertion of the theorem.

  Let $h \colon Z_1 \to Z_2$ be a homeomorphism between $Z$-subsets of
  $X$ that is $\V{}$-close to the inclusion of $Z_1$ into $X$. Observe
  that neither $Z_1$ nor $Z_2$ has to be a $Z(\G{})$-set.  By
  theorem~\ref{thm:approximation within a cover}, there exists a
  closed $Z(\G{})$-embedding $i \colon Z_1 \to X$ that is $\V{}$-close
  to the inclusion of $Z_1$ into $X$. Let $\F{} = \G{}$.  Since $\st
  \V{}$ refines $\G{}$, both $i$ and $i \circ h^{-1}$ map elements of
  $\G{}$ into their stars. Hence partitions $\F{}$ and $\G{}$ are
  compatible with $i$ and with $i \circ h^{-1}$. By
  lemma~\ref{lem:existence of a homeomorphism}, there exists a
  homeomorphism $h_1$ of $X$ onto $X$ that extends $i$ and maps
  elements of $\G{}$ into the corresponding elements of $\st^m \G{}$
  and a homeomorphism $h_2$ of $X$ onto $X$ that extends $i \circ
  h^{-1}$ and maps elements of $\G{}$ into the corresponding elements
  of $\st^m \G{}$. Observe that $h_2^{-1}$ maps elements of $\G{}$
  into the corresponding elements of $\st^{m+1} \G{}$, hence $h_2^{-1}
  \circ h_1$ maps elements of $\G{}$ into the corresponding elements
  of $\st^{2m+1} \G{}$. By the construction, $h_2^{-1} \circ h_1$
  extends $h$ and is $\U{}$-close to the identity. We are done.
\end{proof}

\begin{z-set unknotting theorem}
  A homeomorphism between $Z$-sets in $n$-dimensional N\"obeling
  manifolds extends to an ambient homeomorphism if and only if it
  extends to an $n$-homotopy equivalence.
\end{z-set unknotting theorem}
\begin{proof}
  \index{proof!$Z$-set unknotting theorem} By
  theorem~\ref{thm:nobeling space is strongly universal}, the
  $n$-dimensional N\"obeling space is strongly universal in
  dimension~$n$. By corollary~\ref{cor:nobeling space is ae(n)}, it is
  absolute extensor in dimension $n$. Therefore theorems~\ref{thm:ane
    characterization} and~\ref{thm:strong universality is a local
    property} imply that N\"obeling manifolds are $\N{n}$-spaces.
  Therefore the $Z$-set unknotting theorem is just a special case of
  the main theorem, which we just proved.
\end{proof}

\begin{open embedding theorem}
  Every $n$-dimensional N\"obeling manifold is homeomorphic to
  an open subset of the $n$-dimensional N\"obeling space.
\end{open embedding theorem}
\begin{proof}
  \index{proof!open embedding theorem} By
  proposition~\ref{pro:equivalence with an open subset}, every
  $\N{n}$-space is $n$-homotopy equivalent to an open subset of
  $\nu^n$. Hence, by the main theorem, it is homeomorphic to an open
  subset of $\nu^n$.
\end{proof}

\begin{characterization theorem}
  An $n$-dimensional Polish space is a N\"obeling manifold if and only
  if it is an $ANE(n)$-space and is strongly universal in dimension
  $n$.
\end{characterization theorem}
\begin{proof}
  \index{proof!characterization theorem} By the definition, every open
  subset of~$\nu^n$ is an $n$-dimensional N\"obeling manifold. See the
  proof of the open embedding theorem.
\end{proof}

\begin{sum theorem}
  If a space $X$ is an union of two closed $n$-dimensional N\"obeling manifolds
  whose intersection is also an $n$-dimensional N\"obeling manifold, then $X$
  is an $n$-dimensional N\"obeling manifold.
\end{sum theorem}
\begin{proof}
  \index{proof!sum theorem} By the characterization theorem, every
  $n$-dimensional N\"obeling manifold is an $\N{n}$-space. By
  corollary~\ref{cor:pasting}, an union of two closed $\N{n}$-spaces
  that meet in an $\N{n}$-space, is an $\N{n}$-space. By the
  characterization theorem, every $\N{n}$-space is an $n$-dimensional
  N\"obeling manifold. We are done.
\end{proof}


\newpage\thispagestyle{empty}
\part{Constructing $n$-semiregular and $n$-regular $\N{n}$-covers}

\chapter{Basic constructions in $\N{n}$-spaces}
\label{ch:basic constructions}

\section{Adjustment to a $Z$-collection}
\label{sec:adjustment to a Z-collection}

See page~\pageref{def:z-collection} for the definition of a $Z$-collection.

\begin{lemma}
  \label{lem:adjustment to a Z-collection} Every star-finite locally
  finite closed $\N{n}$-cover~$\F{}$ admits arbitrarily small
  adjustment to a $Z$-collection.  Moreover, given arbitrary
  $Z(\F{})$-set and a closed shrinking of of~$\F{}$, we may assume
  that the $Z$-collection is equal to the shrinking on the
  $Z(\F{})$-set.
\end{lemma}
\begin{proof}
  Let $\F{} = \{ F_i \}_{i \in I}$ denote a star-finite locally finite
  closed $\N{n}$-cover of a space~$X$. Let $\mathcal{I} = \{ J \subset
  I \colon J \neq \emptyset, F_J \neq \emptyset \}$. Order
  $\mathcal{I}$ into a sequence $J_1, J_2, \ldots$ non increasing in
  the order by inclusion. Such ordering exists because $\F{}$ is
  star-finite.  Let $Z$ be a $Z(\F{})$-set and let $\G{} = \{ G_i
  \}_{i \in I}$ be a closed shrinking of $\F{}$.  Let $\U{}$ be an
  open cover of~$X$. We construct a sequence of maps $\{ \varphi_k
  \colon F_{J_k} \to X \}_{k \in \mathbb{N}}$ such that for each~$k$
  the following conditions are satisfied.
  \begin{enumerate}
  \item[($1_k$)] $\varphi_k$ is a $Z(\F{})$-embedding that is
    $\U{}$-close to the inclusion of $F_{J_k}$ into $X$,
  \item[($2_k$)] $\im \varphi_k \subset F_{J_k}$ and $\im \varphi_k
    \cap Z = G_{J_k} \cap Z$,
  \item[($3_k$)] if $J_{k} \cup J_{l} = J_{m}$ and $l < k$, then $\im
    \varphi_{k} \cap \im \varphi_{l} = \im \varphi_{m}$,
  \item[($4_k$)] if $J_k \subsetneq J_{l}$, then $\varphi_k \circ
    \varphi_{l} = \varphi_{l}$.
  \end{enumerate}
  
  Fix $k \in \mathbb{N}$ and assume that for each $l < k$ we
  constructed map $\varphi_l$ that satisfies conditions ($1_l$),
  ($2_l$) and ($3_l$). Let $Z_1 = F_{J_k} \cap (Z \cup \bigcup_{l < k}
  \im \varphi_l)$. By (1), $\varphi_l$ is a $Z(\F{})$-set for each $l
  < k$, so $Z_1$ is a $Z(\F{})$-set. Let $Z_2 = (G_{J_k} \cap Z) \cup
  \bigcup_{l < k \colon J_l \supset J_k} \im \varphi_l$ and observe
  that $Z_2 \subset Z_1$ and $Z_2 \cap Z = G_{J_k} \cap Z$. Let $A =
  F_{J_k} \sqcup_{Z_2} Z_1$ denote the disjoint sum of $F_{J_k}$ and
  $Z_1$ glued along $Z_2$. We treat $F_{J_k}$ and $Z_2$ as subspaces
  of~$A$.  Let $f \colon A \to X$ denote the map induced on $A$ by
  identities on $F_{J_k}$ and $Z_1$.  By the definition, $f$
  restricted to $Z_1$ is a $Z(\F{})$-embedding.  By
  theorem~\ref{thm:approximation within a cover}, there is a
  $\U{}$-approximation $g$ of $f$, $\rel Z_1$ and within $\F{}$, by a
  $Z(\F{})$-embedding. We let $\varphi_k$ to be the restriction of $g$
  to $F_{J_k}$. By the definition of $g$, $\varphi_k$ is a
  $Z(\F{})$-embedding that is $\U{}$-close to the inclusion of
  $F_{J_k}$ into $X$ and its image lies in $F_{J_k}$. By the
  definition of $g$, $\im \varphi_k \cap Z_1 = Z_2$, hence $\im
  \varphi_k \cap Z = \im \varphi_k \cap Z_1 \cap Z = Z_2 \cap Z =
  G_{J_k} \cap Z$. Hence ($1_k$) and ($2_k$) are satisfied. To verify
  $(3_k)$ observe that if $J_k \cup J_l = J_m$ and $l < k$, then $m <
  l$ and by $(3_l)$, $\im \varphi_m \subset \im \varphi_l$. By the
  definition, $\im \varphi_m \subset Z_2$, so $\im \varphi_m \subset
  \im \varphi_k$. Hence $\im \varphi_l \cap \im \varphi_k \supset \im
  \varphi_m$. We have $\im \varphi_l \cap \im \varphi_k = \im
  \varphi_l \cap \im \varphi_k \cap Z_1$, because $\im \varphi_l
  \subset Z_1$ as $l < k$. Hence $\im \varphi_l \cap \im \varphi_k =
  \im \varphi_l \cap Z_2$, because as we noted above $\im \varphi_k
  \cap Z_1 = Z_2$. Recall that $Z_2 = (G_{J_k} \cap Z) \cup
  \bigcup_{l' < k \colon J_{l'} \supset J_k} \im \varphi_{l'}$. By
  ($2_l$) and ($2_m$), $\im \varphi_l \cap G_{J_k} \cap Z = G_{J_l}
  \cap G_{J_k} \cap Z = G_{J_m} \cap Z \subset \im \varphi_m$. Let $l'
  < k$ such that $J_{l'} \supset J_k$ and observe that by ($3_l$) and
  ($3_{l'}$), $\im \varphi_{l} \cap \im \varphi_{l'} \subset \im
  \varphi_m$ as $J_m \subset J_l \cup J_{l'}$. Hence $\im \varphi_l
  \cap \im \varphi_k = \im \varphi_l \cap Z_2 \subset \im \varphi_m$,
  by the definition of $Z_2$. Hence ($3_k$) is satisfied. To verify
  ($4_k$) it is sufficient to see that if $J_k \subsetneq J_l$, then
  the restriction of $\varphi_k$ to $\im \varphi_l$ is equal to the
  inclusion of $\im \varphi_l$ into $X$.  But by the definition the
  restriction of $\varphi_k$ to $Z_2$ is equal to the inclusion and
  $\im \varphi_l \subset Z_2$ by the definition of $Z_2$.
  
  For every $i \in I$ let $k(i)$ denote an integer such that $J_{k(i)}
  = \{ i \}$.  Claim: the collection $\F{}' = \{ F'_i = \im
  \varphi_{k(i)} \}_{i \in I}$ satisfies the assertion. Observe that
  $F'_{J_k} = \im \varphi_{k}$, by (3), therefore $\F{}'$ is a
  $\U{}$-adjustment of $\F{}$, as by (1) every $\varphi_{k}$ is
  $\U{}$-close to the identity. Let $J_k \subsetneq J_{k_1}$.  Since
  $\varphi_{k_1}$ is a $Z(\F{})$-embedding, $F'_{J_{k_1}} = \im
  \varphi_{k_1}$ is a $Z$-set in $F_{J_k}$. Therefore
  $\varphi_{k}(F'_{J_{k_1}})$ is a $Z$-set in $F'_{J_k} = \im
  \varphi_{k}$, but by (4) $\varphi_{k}(F'_{J_{k_1}}) = F'_{J_{k_1}}$,
  therefore $\F{}'$ is a $Z$-collection. By (2) $\F{}'$ is equal to
  $\G{}$ on $Z$, so $\F{}'$ satisfies the assertion.
\end{proof}

\section{Fitting closed $\N{n}$-neighborhoods}

\begin{definition*}
  Let $\F{} = \{ F_i \}_{i \in I}$ denote a collection of subsets of a
  space~$X$. We say that a set $A \subset X$ \df{fits}{fits} $\F{}$ if it is
  an $AE(n)$-space and if for each $J \subset I$ the intersection $A \cap F_J$
  is an $AE(n)$-space.
\end{definition*}

\begin{lemma}\label{lem:fitting approximation}
  Assume that a subset $B$ of a space $X$ fits a collection $\F{} = \{
  F_i \}_{i \in I}$ and $\F{}$ restricted to $B$ is both a
  $Z$-collection and an $\N{n}$-collection. Assume that $f \colon B
  \to X$ is an embedding such that $f^{-1}(F_i) = F_i \cap B$.  Then
  $f(B)$ fits $\F{}$ and $\F{}$ restricted to $f(B)$ is both a
  $Z$-collection and an $\N{n}$-collection.
\end{lemma}
\begin{proof}
  Let $\F{} = \{ F_i \}_{i \in I}$ and let $F_J = \bigcap_{i \in J}
  F_i$ for each non-empty subset $J$ of $I$. We have $F_J \cap f(B) =
  F_J \cap \im f = f(f^{-1}(F_J)) = f(F_J \cap B)$, the last equality
  by the assumption that $f^{-1}(F_i) = F_i \cap B$ for each $i \in
  I$. By the assumption $B$ fits $\F{}$, so $F_J \cap B$ is $AE(n)$
  for each $J \subset I$ and so is $F_J \cap f(B)$. Therefore $\F{}$
  fits $f(B)$. The homeomorphism $f \colon B \to f(B)$ maps collection
  $\{ F_i \cap B \}_{i \in I}$ onto a collection $\{ F_i \cap f(B)
  \}_{i \in I}$. Since the first collection is a $Z$-collection and an
  $\N{n}$-collection in $B$, the second collection is a $Z$-collection
  and an $\N{n}$-collection in $f(B)$.
\end{proof}

\begin{lemma}
  \label{lem:relative approximation within a Z-collection}
  Assume that $\G{} = \{ G_i \}_{i \in I}$ is an $\N{n}$-collection
  and a $Z$-collection in an $\N{n}$-space~$Y$, $V$ is an open subset
  of $Y$ and $f$ is a map from an at most $n$-dimensional Polish
  space~$X$ into~$Y$. If $f$ is a closed embedding into $V$ on
  $f^{-1}(V)$ and $A$ is a closed (in $Y$) subset of $V$, then $f$ is
  approximable $\rel f^{-1}(A)$ by closed embeddings $g$ such that
  $g^{-1}(G_i) = f^{-1}(G_i)$ for each $i \in I$.  Moreover, if $Z$ is
  a subset of $X$ and $f_{|Z}$ is a $Z(\G{})$-embedding, then such
  approximations exist~$\rel f^{-1}(A) \cup Z$.
\end{lemma}
\begin{proof}
  The proof follows the proof of proposition~\ref{pro:approximation
    rel inverse image}. Let $\V{}$ be an open cover of $Y$. Let $U = Y
  \setminus A$ and let $\U{}$ be an open cover of~$U$ obtained via
  lemma~\ref{lem:local approximation} applied to $U$ and $Y$. Without
  a loss of generality, we may assume that $\U{}$ refines $\V{}$. By
  proposition~\ref{pro:ZF is a local property}, the restriction $\G{}
  / U$ of $\G{}$ to $U$ is a $Z$-collection in $U$ and the restriction
  of $f$ to $Z \cap f^{-1}(U)$ is a $Z(\G{} / U)$-embedding into $U$.
  By theorem~\ref{thm:ane characterization} and by
  theorem~\ref{thm:strong universality is a local property}, $U$ is an
  $\N{n}$-space and $\G{} / U$ is an $\N{n}$-collection, hence by
  theorem~\ref{thm:approximation within a Z-collection}, there exists
  a $\U{}$-approximation $g$ of $f_{|f^{-1}(U)}$ $\rel Z \cap
  f^{-1}(U)$ by a closed embedding into $U$ such that $g^{-1}(G_i) =
  f^{-1}(G_i)$ for each $i \in I$. By lemma~\ref{lem:local
    approximation} and by the choice of $\U{}$, the map $g \cup
  f_{|f^{-1}(U)}$ is a closed embedding into $U$ and is $\V{}$-close
  to $f$.
\end{proof}

\begin{lemma}
  \label{lem:regular neighborhoods}
  Let a finite collection $\F{}$ of subsets of an $\N{n}$-space $X$ be both an
  $\N{n}$-collection and a $Z$-collection. Then every $Z(\F{})$-set $A \subset
  X$ that fits~$\F{}$ has arbitrarily small closed $\N{n}$-neighborhoods~$B$
  that fit $\F{}$.  Moreover, we may require that $\F{}$ restricted to $B$ is
  both a $Z$-collection and an $\N{n}$-collection of subsets of $B$.
\end{lemma}

Observe that we do not assume above that~$A$ is an $\N{n}$-space.

\begin{proof}
  The proof is by induction on the cardinality of $\F{}$. Assume that
  $\F{}$ is empty. We have to prove that if $A$ is an absolute
  extensor in dimension $n$ and a $Z$-set in $X$, then it has
  arbitrarily small closed $\N{n}$-neighborhoods that are absolute
  neighborhood extensors in dimension~$n$. Fix an open neighborhood
  $U$ of $A$ and let $V$ be another open neighborhood of $A$ whose
  closure lies in $U$. By lemma~\ref{lem:extension property} and by
  remark~\ref{rem:small n-neighborhoods}, there is a closed
  $\N{n}$-neighborhood $B$ of $A$ with a property that every partial
  map from an $n$-dimensional Polish space into $B$ extends over the
  entire domain to a map into $V$.  By theorem~\ref{thm:nobeling space
    is strongly universal}, there exists a closed embedding $f \colon
  B \to \nu^n$. By the choice of $B$, there is an extension $g \colon
  \nu^n \to V$ of $f^{-1}$.  By corollary~\ref{cor:approximation
    within a cover}, there is an approximation $h$ of $g$ within an
  $\N{n}$-cover $\{ X, B \}$ by a closed embedding with image disjoint
  from $A$. By Whitehead's characterization and by lemma~\ref{lem:set
    of weak equivalences is open}, if the approximation is close
  enough (which we assume), then the inclusion of $h(f(B))$ into $B$
  is an $n$-homotopy equivalence.  Consider a space $C = B
  \sqcup_{h(f(B))} h(\nu^n)$, a disjoint sum of $B$ and $h(\nu^n)$
  glued along $h(f(B))$. By corollary~\ref{cor:pasting}, it is an
  $\N{n}$-space. By lemma~\ref{lem:excision}, it is an $AE(n)$-space.
  Let $i \colon C \to B$ denote a map induced by inclusions $B \subset
  X$ and $h(\nu^n) \subset X$.  It is an embedding on some
  neighborhood of $A$. Hence, by proposition~\ref{pro:approximation
    rel inverse image}, it can be approximated by a closed embedding
  $\rel$ some neighborhood of~$A$.  The image of this embedding is a
  closed $\N{n}$-neighborhood of $A$ that is an $AE(n)$-space. We may
  assume that both approximations done above are close to the
  identity, so this image lies in~$U$.  The basis step of induction is
  done.
  
  Let $\F{} = \{ F_i\}_{i \in I}$, fix $j \in I$, let $J = I \setminus
  \{ j \}$ and assume that the theorem is proved for collections of
  cardinality $J$. Fix a pair of open neighborhoods $U, V$ of $A$ as
  above. Let $\F{}' = \{ F_i \}_{i \in J}$.  By the inductive
  assumption, there is a closed $\N{n}$-neighborhood $B'$ of $A$ in
  $X$ that fits $\F{}'$ and lies in $V$. Let $\F{}'' = \{ F_i \cap F_j
  \}_{i \in J}$.  Since $\F{}$ is both a $Z$-collection and an
  $\N{n}$-collection, $\F{}''$ is a $Z$-collection and an
  $\N{n}$-collection in $F_j$. Since $A$ is a $Z(\F{})$-set, $A \cap
  F_j$ is a $Z(\F{}'')$-set. By the inductive assumption, there is a
  closed $\N{n}$-neighborhood $B''$ of $A \cap F_j$ in $F_j$ that fits
  $\F{}''$ and lies in $F_j \cap \Int_X B'$. Let $B''' = B'
  \sqcup_{B''} F_j$ be a disjoint union of $B'$ and $F_j$ glued along
  $B''$. Let $f \colon B''' \to X$ be a map induced by inclusions on
  $B'$ and $F_j$. By the assumption that $\F{}$ is a $Z$-collection,
  $F_j$ is a $Z(\F{}')$-set. By the construction, $f$ is an embedding
  on an open neighborhood of $A$. By lemma~\ref{lem:relative
    approximation within a Z-collection} applied to a $Z$-collection
  $\F{}'$, map $f$ and a set $F_j \subset B'''$ ($f$ restricted to
  $F_j$ is a $Z(\F{}')$-embedding), there exists an approximation of
  $f$ by a closed embedding $g$ that is equal to $f$ on an union of
  $F_j$ with a neighborhood of $A$ and that satisfies equality
  $g^{-1}(F_i) = f^{-1}(F_i)$ for each $i \in J$. Let $B = g(B')$. It
  is a closed $\N{n}$-neighborhood of $A$. We will show that it fits
  $\F{}$ and that the restriction of $\F{}$ to $B$ is both a
  $Z$-collection and an $\N{n}$-collection in $B$, thus ending the
  inductive step.

  By the construction, $B'$ fits $\F{}'$ and the restriction of
  $\F{}'$ to $B'$ is both an $\N{n}$-collection and a $Z$-collection.
  Hence by lemma~\ref{lem:fitting approximation} applied to an
  embedding $g_{|B'}$, $B$ fits $\F{}'$ and the restriction of $\F{}'$
  to $B$ is both an $\N{n}$-collection and a $Z$-collection. Let $J$
  be a subset of $I$ such that $j \in J$. We have $F_J \cap B = (F_j
  \cap B) \cap F_J = B'' \cap F_J$. Since by the construction $B''$
  fits $\F{}''$, the intersection $B'' \cap F_J$ is an $AE(n)$-space
  and so is $B \cap F_J$. Hence $B$ fits the entire collection $\F{}$.
  Again, by the construction $\F{}''$ restricted to $B''$ is an
  $\N{n}$-collection, hence the entire collection $\F{}$ restricted to
  $B$ is an $\N{n}$-collection. To see that this restriction is a
  $Z$-collection in $B$ observe that $B''$ is a $Z$-set in $B$
  (because $F_j$ is a $Z$-set in $X$ and $B''$ is a closed subset of
  $F_j$ that lies in the interior of $B$), $\F{}''$ is a
  $Z$-collection in $B''$ (by the construction) and $\F{}'$ is a
  $Z$-collection in $B$ (by lemma~\ref{lem:fitting approximation}).
\end{proof}

\section{Patching of holes}

\begin{lemma}
  \label{lem:disc attaching} Assume that $B$ is a closed
  $\N{n}$-subspace and a $Z$-set in an $\N{n}$-space~$X$. If $\Delta$ is a
  $k$-dimensional simplex embedded in $X$ such that $\partial \Delta = \Delta
  \cap B$ and $N$ is an open neighborhood of $\Delta$, then there exists a
  closed (in $X$) $\N{n}$-subset $A$ of $N$ such that $A \cap B$ is an
  $\N{n}$-space, $\Delta \subset A$ and the inclusions $\partial \Delta
  \subset A \cap B$ and $\Delta \subset A$ are $n$-homotopy equivalences.
\end{lemma}
\begin{proof}
  Let $\delta \subset \nu^n$ be a standard $k$-dimensional simplex
  isometrically embedded into a $k$-dimensional hyperplane that lies in
  $\nu^n$. Fix an $\varepsilon > 0$ and let $U = \{ x \in \nu^n \colon d(x,
  \delta) < \varepsilon \}$ and $V = \{ x \in \nu^n \colon d(x, \partial
  \delta) < \varepsilon \}$. Assume that~$\varepsilon$ is so small that the
  inclusion $\partial \delta \subset V$ is an $n$-homotopy equivalence (cf.
  proposition~\ref{pro:intersection of nobeling space with an open set}). The
  inclusion $\delta \subset U$ is also an $n$-homotopy equivalence. By
  corollary~\ref{cor:z-approximation} applied with $D = \delta \setminus
  \partial \delta$, the inclusion $V \subset U$ can be approximated $\rel
  \partial \delta$ by a closed embedding $g \colon V \to U$ with image
  disjoint from $\delta \setminus \partial \delta$. Let $C = g(V)$. Since the
  inclusion $\partial \delta \subset V$ is an $n$-homotopy equivalence, the
  inclusion $g(\partial \delta) = \partial \delta \subset C = g(V)$ is an
  $n$-homotopy equivalence. Therefore, by theorem~\ref{thm:subspace
    retraction}, there is a retraction $C \to \partial \delta$. Since $\delta$
  is an absolute extensor and $\delta \cap C = \partial \delta$, this
  retraction can be extended to a retraction $r \colon U \to \delta$. Note
  that by the construction, $r(C) \subset \partial \delta$.  Let $H \colon
  \delta \to \Delta$ be any homeomorphism. By theorem~\ref{thm:approximation
    within a cover}, $(H \circ r) \colon U \to \Delta$ can be approximated by a
  closed embedding $h \colon U \to X$ $\rel \delta$ such that $h(C) \subset
  B$. By the assumption that $B$ is a $Z$-set, and by
  corollary~\ref{cor:z-approximation} applied with $D = B \setminus h(C)$ and
  $\rel F \cup \delta$, we may assume that $h(C) = h(U) \cap B$. Observe that
  $\Delta \subset h(U)$, $\partial \Delta \subset h(C) = h(U) \cap B$ and
  these inclusions are $n$-homotopy equivalences.  Hence letting $A = h(C)$
  finishes the proof.
\end{proof}

\begin{remark}
  \label{rem:disc attaching}
  Assume that $\Delta$, $A$ and $B$ are closed subsets of $A \cup B$ such that
  $\Delta$, $A$, $B$ and their intersections are at most $n$-dimensional
  $ANE(n)$-spaces. If inclusions $\partial \Delta \subset A \cap B$, $\Delta
  \subset A$ hold and are $n$-homotopy equivalences, then the inclusion $\Delta
  \cup B \subset A \cup B$ is an $n$-homotopy equivalence.
\end{remark}
\begin{proof}
  By corollary~\ref{cor:excision} applied with $X = A$, $A_1 = \Delta$ and
  $A_2 = A \cap B$, the inclusion $\Delta \cup (A \cap B) \subset A$ is an
  $n$-homotopy equivalence. Hence, by lemma~\ref{lem:excision} applied with
  $A_1 = \Delta \cup B$ and $A_2 = A$, the inclusion $\Delta \cup B \subset A
  \cup B$ is an $n$-homotopy equivalence.
\end{proof}



\chapter{Core of a cover}
\label{ch:core of a cover}

Our goal is a proof of theorem~\ref{thm:pump up the regularity} and 
we will attain it by ``patching of holes'' in a cover. In this chapter we
introduce the notion of an $n$-core of a cover that aids the ``patching''
process and allows us to control it.

\begin{definition*}
  Assume that $\F{} = \{ F_i \}_{i \in I}$ is a cover of a space $X$.
  A pair $\langle K, \mathcal{K} \rangle$ consisting of an at most
  $n$-dimensional locally finite simplicial complex $K$ embedded as a
  closed subset of $X$ and of an interior cover $\mathcal{K} = \{ K_i
  \}_{i \in I}$ of~$K$ by subcomplexes is an
  \df{n-core@$n$-core}{$n$-core of $\F{}$} if for every $J \subset I$
  the inclusion $K_J \subset F_J$ holds and induces an
  $n$-homotopy equivalence. If additionally $\mathcal{K}$ is equal
  to~$\F{}$ on $K$, then we say that the core is
  \df{n-core@$n$-core!exact}{exact}. We shall say that~$K$ is an exact
  $n$-core of $\F{}$ if $\langle K, \{ F_i \cap K \}_{i \in I}
  \rangle$ is an $n$-core of $\F{}$.
\end{definition*}


Observe that the conditions of the above definition imply that
collections $\F{}$ and $\mathcal{K}$ must be isomorphic, because the
empty space is not $n$-homotopy equivalent to an non-empty space.

\begin{definition*}
  Let $\F{}$ be a cover of a space $X$, let $\langle K, \mathcal{K} \rangle$
  be an $n$-core of $\F{}$ and let $\U{}$ be an open cover of $X$. We say that
  a cover~$\G{}$ of $X$ is a \df{U-@$\U{}$-!deformation}{$\U{}$-deformation
    of~$\F{}$ with fixed $\langle K, \mathcal{K} \rangle$} if $\langle K,
  \mathcal{K} \rangle$ is an $n$-core of $\G{}$ and if~$G_i$ is a subset
  of~$\st_\U{} F_i$ for each pair of corresponding elements $F_i \in \F{}$ and
  $G_i \in \G{}$ (note that $\F{}$ and $\G{}$ must be isomorphic because they
  have the same $n$-core).  Sometimes we shall write that $\G{}$ is a
  $\U{}$-deformation of $\F{}$ or that $\G{}$ is a deformation of $\F{}$, if
  $\U{}$ and $\langle K, \mathcal{K} \rangle$ are known from the context.
\end{definition*}

Theorem~\ref{thm:construction of an exact core} is the main theorem of this
chapter.

\begin{theorem}
  \label{thm:construction of an exact core} Every closed star-finite interior
  $\N{n}$-cover $\F{}$ with an $n$-core $\langle K, \mathcal{K} \rangle$
  admits arbitrarily small deformation to a closed interior $\N{n}$-cover of
  which $\langle K, \mathcal{K} \rangle$ is an exact $n$-core. Moreover, if
  $L$ is a subcomplex of $K$ and $\mathcal{F}$ is equal to $\mathcal{K}$ on $K
  \setminus L$, then we may additionally require that the deformed cover is
  equal to $\F{}$ on a complement of arbitrarily small open neighborhood of
  $L$.
\end{theorem}

{%
\begin{figure}[ht]
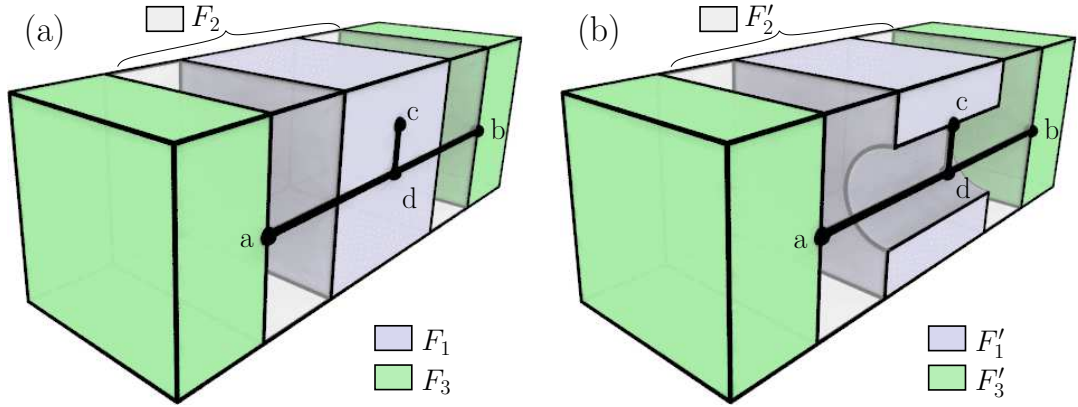

\begin{center}
$\begin{array}{c@{\hspace{5mm}}c}
\resizebox{0.45\textwidth}{!}{\input{psfigures/core1.pstex_t}} &
\resizebox{0.45\textwidth}{!}{\input{psfigures/core2.pstex_t}}
\end{array}$
\end{center}
\caption{{Sample application of
  theorem~\ref{thm:construction of an exact core}.}}
\label{fig:strict core}
\end{figure}
}

We prove it in section~\ref{sec:an exact core}. Here we illustrate its typical
usage. Consider a cover~$\F{}$ drawn on figure~\ref{fig:strict core}(a). It
contains three sets: $F_1$ (dark one), $F_2$ (the large transparent box - it
contains $F_1$) and $F_3$ (light one). Note that $\F{}$ is not $1$-regular:
$F_3$ is disconnected. We want to ``patch'' $F_3$ by attaching arc
$\overline{ab}$ to it. But $\overline{ab}$ crosses $F_1$ and the cover $\{
F_1, F_2, F_3 \cup \overline{ab} \}$, although it is $1$-regular, is not
isomorphic to $\F{}$.  Theorem~\ref{thm:construction of an exact core} gives a
solution to this problem, as follows.

Consider a graph $K = \overline{ab} \cup \overline{cd}$ and its cover
$\mathcal{K} = \{ \{ c \} , K, \{ a, b \} \}$. A pair $\langle K, \mathcal{K}
\rangle$ is a $1$-core of $\F{}$. Assume that $\F{}$ is a closed interior
$\N{1}$-cover and apply theorem~\ref{thm:construction of an exact core} to it.
We get a cover $\{ F'_1, F'_2, F'_3 \}$ drawn on picture~\ref{fig:strict
  core}(b) that has $\langle K, \mathcal{K} \rangle$ as an exact $1$-core. Now
the construction described above works: if we attach $\overline{ab}$ to $F'_3$
we get a $1$-regular cover that is isomorphic to $\F{}$.

\section{The existence of an $n$-core}

First, to show that theorem~\ref{thm:construction of an exact core} has a
point, we prove that every closed star-finite locally finite $\N{n}$-cover has
an $n$-core.  To this end we introduce an auxiliary notion of an outer
$n$-core.

\begin{definition}
  \label{def:n-core} Assume that $\F{} = \{ F_i \}_{i \in I}$ is a
  cover of a space~$X$. A triple $\langle K, \mathcal{K}, k \rangle$ is an
  \df{n-core@$n$-core!outer}{outer $n$-core of $\F{}$} if
  \begin{enumerate}
  \item $K$ is an at most $n$-dimensional locally finite simplicial complex,
  \item $\mathcal{K} = \{ K_i \}_{i \in I}$ is a cover of $K$ by subcomplexes,
  \item $k \colon K \to X$ is a map such that $k(K_i) \subset F_i$ and
    for every $J \subset I$ the restriction of $k$ to $K_J$ is an
    $n$-homotopy equivalence of $K_J$ and $F_J$.
  \end{enumerate}
\end{definition}

For example, if $\F{}$ is an $n$-regular cover and $\lambda$ is its
anticanonical map, then the triple $\langle N(\F{}), \BF, \lambda \rangle$ is
an outer $n$-core of~$\F{}$.

\begin{lemma}
  \label{lem:construction of a core} Assume that $\F{} = \{ F_i
  \}_{i \in I}$ is a star-finite $ANE(n)$-cover of a space $X$ and $\langle L,
  \mathcal{L}, l \rangle$ is a triple that satisfies all conditions of
  definition~\ref{def:n-core} with the exception that the restriction of $l$
  to $L_J$ need not to be an $n$-homotopy equivalence of $L_J$ and $F_J$. Then
  there is an outer $n$-core $\langle K, \mathcal{K}, k \rangle$ of $\F{}$
  such that $L$ is a subcomplex of $K$, $\mathcal{K}$ restricted to $L$ is
  equal to $\mathcal{L}$ and $k$ is an extension of $l$. Moreover,
  if~$\mathcal{L}$ is an interior cover, then we may additionally require that
  $\mathcal{K}$ is also an interior cover.
\end{lemma}
\begin{proof}
  Let $\mathcal{I} = \{ J \subset I \colon F_J \neq \emptyset \}$.
  Order $\mathcal{I}$ into a sequence $\{ J_k \}$ non increasing in
  the order by inclusion. Such ordering exists because $\F{}$ is
  star-finite. By a recursive application of
  lemma~\ref{lem:equivalence with a complex} we construct a sequence
  of triples $\langle W_k, \W{k},w_k \rangle$, starting with $\langle
  W_0, \W{0},w_0 \rangle = \langle L, \mathcal{L},l \rangle$, such
  that
  \begin{enumerate}
    \item $W_{k-1}$ is a subcomplex of an at most $n$-dimensional
    locally finite complex~$W_k$,
    \item $\W{k} = \{ W^k_i \}_{i \in I}$ is a cover of $W_k$ by
    subcomplexes and
    \[
        W^k_i = \left\{
        \begin{array}{ll}
            W^{k-1}_i & i \notin J_k \\
            W^{k-1}_i \cup (W_k \setminus W_{k-1}) &i \in J_k, \\
        \end{array}
        \right.
    \]
    \item $w_k$ is an extension of $w_{k-1}$ onto $W_k$ such that
    $w_k$ restricted to $W^k_{J_k}$ is an $n$-homotopy equivalence
    of $W^k_{J_k}$ and $F_{J_k}$.
  \end{enumerate}
  
  We take $K = \bigcup W_k$ (with the direct limit topology), $k = \bigcup
  w_k$ and $\mathcal{K} = \{ K_i = \bigcup W^k_i \}_{i \in I}$. Obviously the
  triple $\langle K, \mathcal{K}, k \rangle$ has all required properties; the
  complex $K$ is locally finite because of the assumption that $\F{}$ is
  star-finite.
  
  Now assume that $\mathcal{L}$ is an interior cover. A simple modification of
  the construction described above assures that the constructed cover
  $\mathcal{K}$ is also an interior cover. Let $W_{k - 1/2}$ denote a complex
  $W_{k - 1} \cup C_k$, where $C_k = W^{k-1}_{J_k} \times \{ 0, 1\} \cup
  {W^{k-1}_{J_k}}^{(n-1)} \times [0, 1]$ - a cylinder taken over
  $W^{k-1}_{J_k}$ with removed interiors of $n$-dimensional simplices. Let $p
  \colon W_{k - 1/2} \to W_{k - 1}$ be a natural projection. Let $\W{k -
    1/2} = p^{-1}(\W{k-1})$ and $w_{k-1/2} = w_{k-1} \circ p$. Construct
  $\langle W_k, \W{k}, w_k \rangle$ as above, but over $\langle W_{k - 1/2},
  \W{k - 1/2}, w_{k-1/2} \rangle$ instead of $\langle W_{k-1}, \W{k-1},
  w_{k-1} \rangle$ and in such a way, that the closure of $W_k \setminus
  W_{k-1}$ intersects only $W^{k-1}_{J_k} \times \{ 1 \}$. It is easy to
  verify that with such construction if $\mathcal{L}$ is an interior cover,
  then $\mathcal{K}$ is also an interior cover.
\end{proof}

By proving theorem~\ref{thm:embedded core} we will show that every
closed star-finite locally finite $\N{n}$-cover has an $n$-core.

\begin{lemma}
  \label{lem:close approximation}
  For each closed locally finite $ANE(n)$-cover $\F{} = \{ F_i \}_{i
    \in I}$ of a space~$X$ there exists an open cover~$\U{}$ of~$X$
  satisfying the following condition.

  \begin{quote}
    For each $J \subset I$ and each map $\varphi$ from at most
    $(n-1)$-dimensional space into $F_J = \bigcap{i \in J} F_i$, every
    map into $F_J$ that is $\U{}$-close to $\varphi$ is homotopic to
    it (in $F_J$).
  \end{quote}
\end{lemma}
\begin{proof}
  By theorem~\ref{thm:ane characterization}(3) for every $J \subset I$
  there is an open cover $\U{J}$ of $F_J$ such that every two
  $\U{J}$-close maps into $F_J$ from at most $(n-1)$-dimensional
  Polish space are homotopic. Since~$\F{}$ is closed and locally
  finite there is an open cover $\U{}$ of $X$ such that for every $J
  \subset I$ the cover~$\U{}$ restricted to $F_J$ refines $\U{J}$.
  This cover satisfies our claim.
\end{proof}

\begin{theorem}
  \label{thm:embedded core} 
  Let $\F{} = \{ F_i \}_{i \in I}$ be a closed star-finite locally
  finite $\N{n}$-cover of a space $X$. Let $L$ be a locally finite
  simplicial complex embedded as a closed subset of $X$. Let
  $\mathcal{L} = \{ L_i \}_{i \in I}$ be an interior cover of $L$ by
  subcomplexes such that $L_i \subset F_i$ for each $i \in I$. Then
  there exists an $n$-core $\langle K, \mathcal{K} \rangle$ of $\F{}$
  such that $L$ is a subcomplex of $K$ and $\mathcal{K}$ is equal to
  $\mathcal{L}$ on $L$. Moreover, if $Z$ is a $Z(\F{})$-set, then we
  may additionaly require that $K \setminus L$ is disjoint from $Z$.
\end{theorem}
\begin{proof}
  Let $\F{} = \{ F_i \}_{i \in I}$ be a closed star-finite locally
  finite $\N{n}$-cover of a space~$X$. By lemma~\ref{lem:close
    approximation}, there exist an open cover~$\U{}$ of~$X$ that
  satisfies the following condition.

  \vspace{ 1mm}
  \begin{tabular}{rl}
    (*) &
    \begin{tabular}{p{\propwidth}}\noindent\em
    For each $J \subset I$ and each map $\varphi$ from at most
    $(n-1)$-dimensional space into $F_J$, every map into $F_J$ that is
    $\U{}$-close to $\varphi$ is homotopic to it (in $F_J$).
    \end{tabular}
  \end{tabular}
  \vspace{ 1mm}
  
  Let $l$ denote the inclusion of $L$ into $X$. By
  lemma~\ref{lem:construction of a core} the triple $\langle L,
  \mathcal{L}, l \rangle$ can be enlarged to an outer $n$-core
  $\langle W, \mathcal{W}, w \rangle$ of $\F{}$. By
  theorem~\ref{thm:approximation within a cover} and by
  remark~\ref{rem:closed locally compact is ZF} there exists a
  $\U{}$-approximation of $w$ within $\F{}$ $\rel L$ by a closed
  embedding $v$ with image disjoint from $Z \setminus L$. Then
  $v_{|W_J}$ is $\U{}$-close to $w_{|W_J}$ so by (*), the inclusion of
  $v(W_J)$ into $F_J$ is an $n$-homotopy equivalence (see the proof of
  Whitehead's characterization).  Therefore $\langle v(W),
  v(\mathcal{W}) \rangle$ is an $n$-core of $\F{}$.  Since $v_{|L}$ is
  equal to~$l$, $L$ is a subcomplex of $v(W)$ and $v(\mathcal{W})$ is
  equal to $\mathcal{L}$ on~$L$.
\end{proof}

\begin{remark}
  \label{rem:core of a small adjustment}
  Let $\F{}$ be a closed star-finite locally finite $\N{n}$-cover of a
  space $X$. Let $\U{}$ be a cover satisfying the condition stated in
  lemma~\ref{lem:close approximation}. Assume that $\langle K,
  \mathcal{K} \rangle$ is an $n$-core of $\F{}$, where $\mathcal{K} =
  \{ K_i \}_{i \in I}$. If $\G{} = \{ G_i \}_{i \in I}$ is a
  $\U{}$-adjustment of $\F{}$ such that $K_i \subset G_i$ for each $i
  \in I$, then $\langle K, \mathcal{K} \rangle$ is an $n$-core of
  $\G{}$. To prove it fix $J \subset I$ and let $h \colon G_J \to F_J$
  be a homeomorphism that is $\U{}$-close to the identity. The
  inclusion of $h(K_J)$ into $F_J$ is $\U{}$-close to the identity on
  $K_J$, hence it is $n$-homotopic to it.  Therefore the inclusion of
  $K_J = h^{-1}(h(K_J))$ into $G_J = h^{-1}(F_J)$ is a weak
  $n$-homotopy equivalence.
\end{remark}

\section{An $n$-core of a limit of a sequence of deformations}

The purpose of this section is to prove an analogue of
proposition~\ref{pro:sequence of adjustments} for sequences of
deformations.  Observe that, even if deformations are done on a
discrete family of open subsets, the limit of the sequence doesn't
have to be a deformation of the original cover. Figure~\ref{fig:limit}
gives an appropriate example: the singleton $K$ is weak $1$-homotopy
equivalent to the set $F_n$ obtained from $F_0$ by disconnecting $n$
``bridges''. But it is not weak $1$-homotopy equivalent to the set
that we get in the limit of this construction, because $F_\infty =
\bigcap F_n$ is disconnected.

{%
\begin{figure}[ht]
\begin{center}
\resizebox{0.60\textwidth}{!}{\input{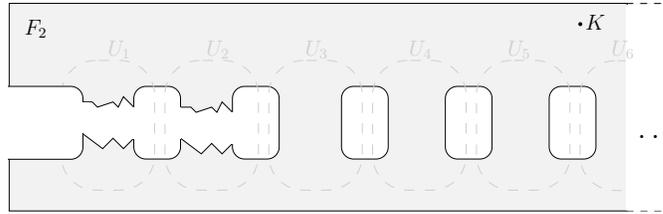}}
\end{center}
\caption{{$K$ is not a $1$-core of $F_\infty = \bigcap
  F_n$.}}
\label{fig:limit}
\end{figure}
}

\begin{lemma}
  \label{lem:limit of deformations}
  Let $\{ U_m \}_{m \geq 1}$ be a discrete family of open subsets of a space
  $X$ and for each $m$ let $V_m$ be an open subset of $U_m$. Let $\F{m}$ ($m
  \geq 0$) be a closed star-finite locally finite $\N{n}$-cover of $X
  \setminus \bigcup_{k > m} V_k$ and assume that $\F{m}$'s are pairwise
  isomorphic. Assume that $\F{m-1}$ is equal to $\F{m}$ on $X \setminus U_m$
  and that it is a shrinking of $\F{m}$ on $U_m \setminus V_m$. If $\langle K,
  \mathcal{K} \rangle$ is an $n$-core of $\F{m}$ for each $m$, then it is an
  $n$-core of $\lim_{m \to \infty} \F{m}$.
\end{lemma}

Lemma~\ref{lem:limit of deformations} is a direct consequence of the
following lemma.

\begin{lemma}
  \label{lem:passing to the limit} Let $\{ U_m \}_{m \geq 1}$ be a
  discrete family of open subsets of a space $X$ and for each $m$ let $V_m$
  be an open subset of $U_m$. Let $\F{m}$ ($m \geq 0$) be a closed
  star-finite locally finite $\N{n}$-cover of $X \setminus \bigcup_{k > m}
  V_k$ and assume that $\F{m}$'s are pairwise isomorphic. If $\F{m}$ is equal
  to $\F{m-1}$ on the complement of $U_m$, then $\F{} = \lim_{m \to \infty}
  \F{m}$ is a closed star-finite locally finite $\N{n}$-cover of $X$.  Assume
  that $\F{m}$ is a swelling of $\F{m-1}$ on $U_m \setminus V_m$ and let
  $\langle K_m, \mathcal{K}_m \rangle$ be an $n$-core of $\F{m}$, where
  $\mathcal{K}_m = \{ K^m_i \}_{i \in I}$. If $K_{m-1}$ is a subcomplex of
  $K_m$, $\mathcal{K}_{m-1}$ is a shrinking of $\mathcal{K}_m$ and $K =
  \bigcup K_m$ is a locally finite complex embedded as a closed subset of $X$,
  then $\langle K, \mathcal{K} \rangle$ is an $n$-core of $\F{}$, where
  $\mathcal{K} = \{ K_i \}_{i \in I}$ and $K_i = \bigcup_{m \geq 1} K^m_i$.
\end{lemma}
\begin{proof}
  To prove that $\F{}$ is a closed locally finite $\N{n}$-cover of $X$ it is
  sufficient to observe that the sequence $\F{m}$ stabilizes locally and
  consists of closed locally finite $\N{n}$-covers. Also $\F{}$ must be
  star-finite because $\F{m}$'s are pairwise isomorphic. What is left is to
  check that $\langle K, \mathcal{K} \rangle$ is an $n$-core of $\F{}$. Let
  $\F{m} = \{ F^m_i \}_{i \in I}$ and $\F{} = \{ F_i \}_{i \in I}$ where $F_i
  = \bigcap_{j} \bigcup_{k \geq j} F^k_i$. Fix $J \subset I$. By the
  definition it is sufficient to prove that the inclusion of $K_J$ into $F_J$
  is a weak $n$-homotopy equivalence.

  Fix $k < n$ and a map $\varphi \colon S^k \to F_J$. Since $S^k$
  is compact $\im \varphi$ intersects only finitely many $U_m$'s
  so there is $M$ such that $\im \varphi \cap U_m = \emptyset$ for
  every $m \geq M$. Therefore $\varphi$ is homotopic with a map
  into $K_M \cap F^M_J$ in $F^M_J \setminus \bigcup_{m > M} V_m$.
  But $F^M_J \setminus \bigcup_{m > M} V_m \subset F_J$ and $K_M
  \cap F^M_J \subset K \cap F_J$ so the inclusion induces
  epimorphisms on homotopy groups of dimensions less than $n$.

  Fix a map $\varphi \colon S^k \to K \cap F_J$ and assume that it
  has an extension $\varPhi \colon D^{k+1} \to F_J$. By
  compactness and local finiteness of $K$ there is $M$ such that
  $\im \varphi \subset K_M$ and $\im \varPhi \subset F_J \setminus
  \bigcup_{m > M} U_m \subset F^M_J$. But $K_M$ is a core of
  $\F{M}$ therefore $\varphi$ must be contractible in $K_M \cap
  F^M_J$. But $K_M \cap F^M_J \subset K \cap F_J$ so the inclusion
  induces monomorphisms on homotopy groups of dimensions less than $n$
  and we are done.
\end{proof}

\section{Proof of theorem~\lowercase{\ref{thm:construction of an exact core}}}
\label{sec:an exact core}

Let $K$ be a locally finite simplicial complex embedded as a closed
subset of an $\N{n}$-space~$X$ and let $L$ be a subcomplex of $K$. Let
$U_L$ be an open neighborhood of $L$. Let $\mathcal{K} = \{ K_i \}_{i
  \in I}$ be an interior cover of~$K$ by subcomplexes. Let $\F{} = \{
F_i \}_{i \in I}$ be a star-finite interior $\N{n}$-cover of $X$ such
that $\langle K, \mathcal{K} \rangle$ is an $n$-core of $\F{}$. We fix
a triangulation~$\tau$ of $K$ such that for each simplex $\delta$ in
$\tau$ there exist $i$ and $j$ in $I$ such that $\delta \subset \Int_X
F_i$ (this is possible because $\F{}$ is an interior cover) and
$\delta \subset \Int_K K_j$ (if a triangulation does not have this
property, it suffices to take its first barycentric subdivision,
because $\mathcal{K}$ is an interior cover). Let $\tau_L$ be a subset
of $\tau$ that contains simplices of~$L$.

For each $l \in \{ 0, 1, \ldots, n+1 \}$ and $m \in \{ 0, 1, \ldots, \infty
\}$ we let $\F{l,m} = \{ F^{l,m}_i \}_{i \in I}$ denote a collection (not yet
defined) of subsets of $X$. Let $(1_{l,m})$, $(2_{l,m})$, $(3_{l,m})$ and
$(4^{\delta, i}_{l,m})$ denote the following statements, with $0 \leq l \leq
n+1$, $0 \leq m \leq \infty$, $\delta \in \tau$ and $i \in I$.
\begin{equation}\tag{$1_{l,m}$}
  \F{l,m} \text{ is a closed star-finite interior } \N{n}\text{-cover of } X.
\end{equation}
\begin{equation}\tag{$2_{l,m}$}
  \langle K, \mathcal{K} \rangle \text{ is an } n\text{-core of } \F{l,m}.
\end{equation}
\begin{equation}\tag{$3_{l,m}$}
  \text{For each } \delta \in \tau \text{ there is } i \in I \text{
  such that } \delta \subset \Int_X F^{l,m}_i.
\end{equation}
\begin{equation}\tag{$4^{\delta,i}_{l,m}$}
  \delta \subset \Int_K K_i \Rightarrow \delta \subset \Int_X
  F^{l,m}_i.
\end{equation}

On a very general level, the proof is organized as follows. We let
$\F{0,0} = \F{}$ and construct $\F{l,m}$ recursively.  In $(l,m)$th
step of the construction, with $0 \leq l \leq n+1$ and $0 \leq m <
\infty$, we construct a collection $\F{l, m+1}$. If $l < n+1$, then
$\F{l, m+1}$ is a deformation of $\F{l,m}$ with fixed $n$-core
$\langle K, \mathcal{K}\rangle$. If $l = n+1$, then $\F{l, m+1}$ is an
adjustment of $\F{l,m}$. In $l$th limit step of the construction, with
$0 \leq l \leq n+1$, we let $\F{l, \infty} = \lim_{m \to \infty}
\F{l,m}$. We also let $\F{l+1, 0}$ denote $\F{l, \infty}$ for each $0
\leq l < n+1$. The collection $\F{n+1, \infty}$ shall satisfy the
assertion of the theorem.

A general constraint on the construction is that conditions
$(1_{l,m})$, $(2_{l,m})$ and $(3_{l,m})$ are satisfied for all $l$ and
$m$. Observe that $(1_{0,0})$, $(2_{0,0})$ and $(3_{0,0})$ are
satisfied by the assumptions of the theorem (we did let $\F{0,0} =
\F{}$).

The first part of the proof includes the first $n+1$ limit steps (i.e.
the construction of $\F{0, \infty}, \F{1,\infty}$ up to $\F{n,
  \infty}$). In $l$th limit step of the construction, with $0 \leq l
\leq n$, we construct a cover $\F{l, \infty}$ that satisfies
$(4^{\delta, i}_{l, \infty})$ for each at most $l$-dimensional simplex
$\delta \in \tau$ and each $i \in I$.  A general constraint in the
first part of the construction is that in its $(l,m)$th step, with $0
\leq l \leq n$ and $0 \leq m < \infty$, \emph{the fourth condition is
  preserved}, i.e.  the implication $(4^{\delta,i}_{l,m}) \Rightarrow
(4^{\delta,i}_{l,m+1})$ is true for each at most $l$-dimensional
simplex $\delta$ of triangulation $\tau$ and each $i \in I$. The first
part of the proof concludes with a construction of $\F{n, \infty}$
that satisfies $(1_{n,\infty})$, $(2_{n,\infty})$, $(3_{n,\infty})$
and, for each $\delta \in \tau$ and each $i \in I$,
$(4^{\delta,i}_{n,\infty})$ (no other property of $\F{n,\infty}$ is
important in the second part of the proof).

The second part of the proof includes the last limit step. It starts with a
cover $\F{n+1, 0} = \F{n, \infty}$ satisfying $(1_{n+1,0})$, $(2_{n+1,0})$,
$(3_{n+1,0})$ and $(4^{\delta,i}_{n+1,0})$ for each $\delta \in \tau$ and $i
\in I$. A cover $\{ V_m \}_{m \geq 0}$ of $K$, consisting of open subsets of
$X$, is picked and in $(n+1, m)$th step of the construction an adjustment
$\rel X \setminus V_m$ of $\F{n+1, m}$ is done in such a way that
$\F{n+1,m+1}$ is equal to $\mathcal{K}$ on $V_m$. Under some additional
constraints, the limit cover $\F{n+1,\infty} = \lim_{m \to \infty} \F{n+1,m}$
satisfies the assertion of the theorem.

\subsection*{First part.}

In $l$th limit step, with $0 \leq l \leq n$, we shall take $\F{l, \infty} =
\lim_{m \to \infty} \F{l, m}$ and argue that $\langle K, \mathcal{K} \rangle$
is an $n$-core of $\F{l,\infty}$. To this end we shall apply
lemma~\ref{lem:limit of deformations}. To be able to do that we must construct
an auxiliary $n$-core of $\F{0,0}$ that is disjoint from $K$. We start with
its construction.

Let $K' = K \times \{ 0, 1 \} \cup K^{(n-1)} \times [0,1]$ and let $p
\colon K' \to K$ be a projection along $[0,1]$. By
theorem~\ref{thm:approximation within a cover}, there exists an
approximation $\tilde p \colon K' \to X$ of $p$ within $\F{0,0}$ $\rel
K \times \{ 0 \}$ by a closed embedding. Let $L = \tilde p(K \times \{
1 \})$ and $M = \tilde p(K')$. Note that $K = \tilde p(K \times \{ 0
\}) \subset M$. Let $\mathcal{L} = \{ L_i \}_{i \in I}$ and
$\mathcal{M} = \{ M_i \}_{i \in I}$, where $M_i = \tilde
p(p^{-1}(K_i))$ and $L_i = M_i \cap L$. Note that $\langle L,
\mathcal{L} \rangle$ and $\langle M, \mathcal{M} \rangle$ are
$n$-cores of $\F{0,0}$ and $L$ is disjoint from some open neighborhood
of $K$. Also note that if $M_i \subset F^{l,k}_i$ for some $0 \leq l
\leq n$ and $0 \leq k \leq \infty$, and each $i \in I$, then $\langle
L, \mathcal{L} \rangle$ is an $n$-core of $\F{l,k}$ if and only if
$\langle K, \mathcal{K} \rangle$ is an $n$-core of $\F{l,k}$.

Let $\delta_0, \delta_1, \ldots$ be a sequence of all $0$-dimensional
simplices of $\tau_L$. Let $J_k = \{ i \in I \colon \delta_k \subset
\Int_X F^{0,0}_i \}$. Fix an open neighborhood $U_k$ of $\delta_k$ in
$X$ whose closure lies in the interior of $F^{0,0}_{J_k}$ and meets
only those simplices of~$K$ that contain~$\delta_k$. Without a loss of
generality we may assume that $\{ U_k \}$ is a discrete family and
each $U_k$ is disjoint from $L$. Let $h_k \colon U_k \to U_k$ be an
approximation of $id_{U_k}$ within $\F{0,0}$ by a closed embedding
with image disjoint from $\delta_k$. We assume that the approximation
is such that $h_k \cup id_{X \setminus U_k}$ is a closed embedding
(see lemma~\ref{lem:local adjustment}). We let $V_k = U_k \setminus
h_k(U_k)$. Note that the choice of sets $V_k$ will be important only
when passing to the limit and for now the reader may think that $V_k =
U_k$.

We are going to construct a sequence $\F{0, 1}, \F{0, 2}, \ldots$ in such a
way that the following conditions will be satisifed for each $k > 0$.
\begin{enumerate}\addtocounter{enumi}{4}
\item[($5_k$)] $\F{0, k}$ is a deformation of $\F{0,k-1}$ $\rel X \setminus
  V_m$, for $m$ such that $2m < k \leq 2m+2$,
\item[($6_k$)] $M_i$ is a subset of $F^{0,k}_i$ for each $i \in I$.
\item[($7_{2k}$)] the condition $(4^{\delta_{k-1}, i}_{0, 2k})$ is satisfied,
\end{enumerate}
Also, as we stated in the introductory part of the proof, for each $m$
conditions $(1_{0,m})$, $(2_{0,m})$ and $(3_{0,m})$ will be satisfied and the
fourth condition will be preserved in each step of the construction.

Assume that we already constructed $\F{0,1}, \F{0,2}, \ldots,
\F{0,2k}$.  Let $F_i = F^{0, 2k}_i \cap V_k$ and $J'_k = \{ i \in I
\setminus J_k \colon F^{0,2k}_i \cap V_k \neq \emptyset \}$. By
star-finiteness of $\F{0,2k}$, the set $J'_k$ is finite. By
lemma~\ref{lem:adjustment to a Z-collection}, there exists an
adjustement of a collection $\{ F_i \}_{i \in J'_k}$ to a
$Z$-collection $\{ F'_i \}_{i \in J'_k}$ that is equal to $\{ M_i
\}_{i \in J'_k}$ on $M \cap V_k$. Let
\begin{equation*}
  F^{0,2k+1}_i = \left\{
  \begin{array}{ll}
    F'_i \cup (F^{0,2k}_i \setminus V_k) & i \in J'_k \\
    F^{0,2k}_i & i \in I \setminus J'_k.
  \end{array} \right.
\end{equation*}
By lemma~\ref{lem:local adjustment}, if the adjustement was small
enough (which we assume), then the collection $\F{0,2k+1}$ is an
adjustment of $\F{0,2k}$ $\rel X \setminus V_k$. Therefore
$\F{0,2k+1}$ is a closed star-finite $\N{n}$-collection. It is an
interior cover because the adjustment was done with fixed $J_k$ $\rel
X \setminus V_k$ and the closure of~$V_k$ lies in the interior
of~$F^{0,2k}_{J_k}$ by the definition of~$U_k$. By
remark~\ref{rem:core of a small adjustment}, if the adjustement was
small enough (which we assume), then the condition $(2_{0,2k+1})$
holds.  Condition $(3_{0,2k+1})$ is satisfied and the fourth condition
is preserved by the construction, because the closure of $V_k$ is
disjoint from every simplex not containing~$\delta_k$ and if $\delta_l
\supset \delta_k$, then $J_l \subset J_k$ and the adjustement is done
with fixed~$J_k$. The condition $(5_{2k+1})$ is satisfied because
$\F{0,2k+1}$ is an adjustement $\rel X \setminus V_k$ of $\F{0,2k}$
and $\langle K, \mathcal{K} \rangle$ is its $n$-core by
$(2_{0,2k+1})$. The condition $6_{2k+1}$ follows from the construction
of $\{ F'_i \}_{i \in J'_k}$. Note that there is no condition
$(7_{2k+1})$ to check.

By lemma~\ref{lem:regular neighborhoods}, there is a closed (in $X$)
$\N{n}$-neighborhood $A_k$ of $\delta_k$ that lies in $V_k$ and such that for
each $J''_k \subset J'_k$ (hence for each $J''_k \subset I$) the intersection
$A_k \cap F^{0,2k+1}_{J''_k}$ is an $(AE(n) \cap \N{n})$-space. Let
\begin{equation*}
F^{0,2k+2}_i = \left\{
    \begin{array}{ll}
      F^{0,2k+1}_i \cup A_k & i \in J'_k \\
      F^{0,2k+1}_i & i \in I \setminus J'_k.
    \end{array}\right.
\end{equation*}
The collection $\F{0,2k+2}$ is $\N{n}$-collection by
corollary~\ref{cor:pasting}, hence $(1_{0,2k+2})$ is satisfied. By
lemma~\ref{lem:excision}, the pair $\langle K, \mathcal{K} \rangle$ is
its $n$-core (we assume that $A_k$ is so small that $\F{0,2k+2}$ is
isomorphic to $\F{0,2k+1}$). Hence $(2_{0, 2k+2})$ is satisfied.
Conditions $(3_{2k+2})$ and $(6_{2k+2})$ are satisfied and the fourth
condition is preserved, because $F^{2k+1}_i \subset F^{2k+2}_i$ for
each $i \in I$. The condition $(5_{2k+2})$ is satisfied by
$(2_{0,2k+2})$ and because $A_k$ is a subset of $V_k$.  By the choice
of $A_k$, $(4^{\delta_k,i}_{0,2k+2})$ is satisfied for each $i \in I$,
hence $(7_{2k+2})$ is satisfied.

Let $\F{0, \infty} = \lim_{k \to \infty} \F{0, k}$. We are going to
show that $(1_{0,\infty})$, $(2_{0,\infty})$, $(3_{0,\infty})$ and,
for each $0$-dimensional simplex $\delta \in \tau$ and each $i \in I$,
$(4^{\delta,i}_{0,\infty})$ are satisfied. Since the sequence
$\F{0,0}, \F{0,1}, \ldots$ stabilizes locally and $(1_{0,k})$,
$(3_{0,k})$ and $(7_{2k})$ hold for each $k$, only verification of
$(2_{0,\infty})$ needs an argument. What we have to show is that
$\langle K, \mathcal{K} \rangle$ is an $n$-core of $\F{0, \infty}$.
Let $H_k \colon X \to X$ be a closed embedding defined as the identity
on $X \setminus \bigcup_{l \geq k} U_l$ and by the formula $H_k(x) =
h_l(x)$ for each $x \in U_l$ with $l \geq k$. Let $\G{k} = \{ G^k_i
\}_{i \in I}$ and let $G^k_i = H_k(F^{0, 2k}_i)$. Observe that $\G{k}$
is a closed interior $\N{n}$-cover of $X \setminus \bigcup_{l \geq k}
V_l$ and $\langle L, \mathcal{L} \rangle$ is its $n$-core. Also,
$\G{k}$ is equal to $\G{k-1}$ on $X \setminus U_k$ and $\G{k-1}$ is a
shrinking of $\G{k}$ on $U_k \setminus V_k$. By lemma~\ref{lem:limit
  of deformations}, $\langle L, \mathcal{L} \rangle$ is an $n$-core of
$\lim_{k \to \infty} \G{k}$. But $\lim_{k \to \infty} \G{k} = \lim_{k
  \to \infty} \F{0, k}$, so $\langle L, \mathcal{L} \rangle$, and also
$\langle K, \mathcal{K} \rangle$, is an $n$-core of $\F{0, \infty}$.

Now we are going to describe a similar construction for simplices of dimension
$l > 0$. It presents additional difficulty because simplices of positive
dimensions have non-empty boundaries.

Let $0 < l \leq n$ and assume that we already constructed $\F{l, 0}$
and $(1_{l,0})$, $(2_{l,0})$, $(3_{l,0})$ and $(4^{\delta,i}_{l,0})$
are satisfied for each at most $(l-1)$-dimensional simplex $\delta \in
\tau$ and each $i \in I$. Let $\delta_0, \delta_1, \ldots$ be a
sequence of all $l$-dimensional simplices in $\tau_L$. For each $k$
let $\delta'_k$ be a shrinking of $\delta_k$, by a homothety with a
ratio close to $1$ and with a center equal to the barycenter of
$\delta_k$, such that the following conditions are satisfied with $m =
0$.
\begin{enumerate}\addtocounter{enumi}{7}
\item $\{ \delta'_k \}$ is a discrete collection in $X$,
\item[$(9_m)$] if $\delta_k \subset \Int_K K_i$, then $\Cl_X (\delta_k
  \setminus \delta'_k) \subset \Int_X F^{l,m}_i$.
\end{enumerate}
The last condition can be satisfied because by the assumption every proper
face $\delta$ of $\delta_k$ satisfies ($4^{\delta,i}_{l,0}$) for each $i \in
I$.

Let $J_k = \{ i \in I \colon \delta_k \subset \Int_X F^{l,0}_i \}$.
For each $k$ fix an open neighborhood $U_k$ of $\delta'_k$ whose
closure lies in the interior of $F^{l,0}_{J_k}$ and is disjoint from
every simplex of $\tau$ that is disjoint from $\delta'_k$. Without a
loss of generality we may assume that $\{ U_k \}$ is a discrete family
in $X$ and each $U_k$ is disjoint from $L$.  Define $V_k$ in the same
manner as in the $0$-dimensional case. Let $J'_k = \{ i \in I
\setminus J_k \colon F^{l,0}_i \cap V_k \neq \emptyset \}$. Let
$\widehat J_k = \{ i \in J'_k \colon \delta_k \subset \Int_K K_i \}$.
Let $m_0 = 0$ and $m_{k+1} = m_k + 4|\widehat J_k|$, where $|\widehat
J_k|$ denotes the number of elements of $\widehat J_k$ (the set $J_k$
is finite because $\F{l, 0}$ is locally finite).  For each~$k$ order
$\widehat J_k$ into a sequence $j^0_k, j^1_k, \ldots, j^{|\widehat
  J_k|-1}_k$.

We are going to construct a sequence $\F{l, 1}, \F{l, 2}, \ldots$ in such a
way that the following conditions will be satisfied for each $k > 0$.
\begin{enumerate}\addtocounter{enumi}{9}
\item[($10_k$)] $\F{l,k}$ is a deformation of $\F{l,k-1}$ $\rel X \setminus
  V_m$, for $m$ such that $m_m < k \leq m_{m+1}$,
\item[($11_k$)] $M_i$ is a subset of $F^{l,k}_i$ for each $i \in I$,
\item[($12_{k,m}$)] the condition $(4^{\delta_k, j^m_k}_{l, m_k + 4m+4})$ is
  satisfied, where $0 \leq m < |\widehat J_k|$.
\end{enumerate}
Also, for each $m$, conditions $(1_{l,m})$, $(2_{l,m})$, $(3_{l,m})$ and
$(9_m)$ will be satisfied and the fourth condition will be preserved in each
step of the construction.

Fix $0 \leq k < \infty$ and $0 \leq \hat m < |\widehat J_k|$ and let $m = m_k
+ 4\hat m$. Assume that we already constructed $\F{l, m}$. Let $j$ denote
$j^{\hat m}_k$. We shall construct covers $\F{0,m+1}, \ldots, \F{0,m+4}$. Our
goal is to meet condition $(12_{k,\hat m})$ and verify that other conditions
specified above are satisfied. Let $\hat\delta_k$ be a shrinking of
$\delta'_k$ by a homothety with a ratio less than $1$ and a center equal to
the barycenter of $\delta'_k$, such that $(9_m)$ is satisfied with
$\hat\delta_k$ substituted for $\delta'_k$. Let $\hat J = \{ i \in J'_k \colon
i \neq j^q_k, q \leq \hat m \}$ and let $\hat V_k$ be an open neighborhood of
$\hat\delta_k$, whose closure lies in the intersection of $(V_k \setminus
\delta_k) \cup (\delta'_k \setminus \partial\delta'_k)$ with $\Int_X
F^{l,m}_{J_k \cup \{ j^q_k \colon q < \hat m \} }$. Let $m'$ be an integer
greater than or equal to $m$ and less than $m+4$. We shall construct $\F{l,
  m'+1}$ as a collection equal to $\F{l, m'}$ on $\rel X \setminus \hat V_k$
and such that $F^{l,m'+1}_i = F^{l,m'}_i$ for each $i \in J_k \cup \{ j^q_k
\colon q < \hat m \}$. This way, the condition $(9_{m'+1})$ will be satisfied
and the condition $(10_{m'+1})$ will follow from $(2_{l,m'+1})$.

Let $P = \hat V_k \cap \Int_X F^{l,m}_{j}$. Consider a collection $\{ F_i \}_{i \in
  \hat J}$, where $F_i = F^{l, m}_i \cap P$. By lemma~\ref{lem:adjustment to a
  Z-collection}, it can be adjusted to a $Z$-collection $\{ F'_i \}_{i \in \hat
  J}$ that is equal to $\{ M_i \}_{i \in \hat J}$ on $M \cap P$.  Let
\begin{equation*}
  F^{l,m+1}_i = \left\{
  \begin{array}{ll}
    F'_i \cup (F^{l,m}_i \setminus P) & i \in \hat J \\
    F^{l,m}_i & i \in I \setminus \hat J.
  \end{array} \right.
\end{equation*}
As before, if we assume that the adjustement is small enough, then conditions
$(1_{l,m+1})$, $(2_{l,m+1})$, $(3_{l,m+1})$ and $(11_{m+1})$ are satisfied and
the fourth condition is preserved. By our earlier remark, $(9_{m+1})$ and
$(10_{m+1})$ are satisfied as well.

Let $V'_k$ be an open subset of $\hat V_k$ such that $\Cl_X V'_k \cap
\Cl_X (\delta_k \setminus \hat\delta_k) = \emptyset$ and $V'_k \cup P = \hat
V_k$.  Consider a collection $\{ F_i \}_{i \in \hat J \cup \{ j \}}$ such that
$F_i = F^{l, m+1} \cap V'_k$. By lemma~\ref{lem:adjustment to a Z-collection},
it can be adjusted to a $Z$-collection $\{ F'_i \}_{i \in \hat J \cup \{ j
  \}}$ that is equal to $\{ M_i \}_{i \in \hat J \cup \{j\}}$ on $M \cap
V'_k$. Let
\begin{equation*}
  F^{l,m+2}_i = \left\{
  \begin{array}{ll}
    F'_i \cup (F^{l,m+1}_i \setminus V'_k) & i \in \hat J \cup \{ j \}  \\
    F^{l,m+1}_i & i \in I \setminus \hat J, i \neq j.
  \end{array} \right.
\end{equation*}
Again, if the adjustement is small enough, then the collection $\F{l,m+2}$
satisfies all prescribed conditions. We claim that the following conditions
are satisfied as well.
\begin{enumerate}\addtocounter{enumi}{12}
\item $\{ F^{l, m+2}_i \cap \hat V_k \}_{i \in \hat J}$ is an
  $\N{n}$-collection and a $Z$-collection in $\hat V_k$,
\item $\{ F^{l, m+2}_i \cap \hat V_k \cap F^{l, m+2}_{j} \}_{i \in \hat J}$ is
  an $\N{n}$-collection and a $Z$-collection in $\hat V_k \cap F^{l,m+2}_{j}$,
\item $\Cl_X (\delta_k \setminus \hat\delta_k) \subset \Int_X F^{l, m +
    2}_{j}$.
\end{enumerate}

By the construction, $\{ F^{l,m+2}_i \}_{i \in \hat J}$ is equal to
$\mathcal{M}$ on $\hat V_k \cap M$, so $\hat\delta_k$ fits $\{
F^{l,m+2}_i \}_{i \in \hat J}$ and by $(13)$ $\{ F^{l,m+2}_i \cap \hat
V_k \}_{i \in \hat J}$ is a $Z$- and $\N{n}$-collection in $\hat V_k$.
Hence, by lemma~\ref{lem:regular neighborhoods}, there is a closed (in
$X$) $\N{n}$-subset $A_k$ of $\hat V_k$ that fits $\{ F^{l, n+2}_i
\cap \hat V_k \}_{i \in \hat J}$ and such that $\hat\delta_k \subset
\Int_X A_k$ and $\{ F^{l,m+2}_i \}_{i \in \hat J}$ restricted to $A_k$
is both a $Z$- and an $\N{n}$-collection.  Observe that by the choice
of $\hat V_k$ and $\hat J$, $A_k$ fits $\{ F^{l, m+2}_i \cap \hat V_k
\}_{i \in I, i \neq j}$. By lemma~\ref{lem:fitting approximation}, we
may adjust $A_k$ to a closed (in $X$) $\N{n}$-subset $A'_k$ of $\hat
V_k$ such that $A'_k \cap \delta_k = \hat\delta_k$, $\hat\delta_k
\setminus \partial \hat\delta_k \subset \Int_X A'_k$ and $\{
F^{l,m+2}_i \}_{i \in \hat J}$ restricted to $A'_k$ is both a $Z$- and
an $\N{n}$-collection. Let $W_k = (F^{l, m+2}_j \cap \Int_X A'_k) \cup
(A'_k \cap \Int_X F^{l, m+2}_j)$. It is an open neighborhood of
$\hat\delta_k$ in $A'_k \cap F^{l, m+2}_{j}$ and it is an
$\N{n}$-space. By (14) and by lemma~\ref{lem:regular neighborhoods},
there is a closed $\N{n}$-subset $B_k$ of $W_k$ that fits $\{ F^{l,
  m+2}_i \cap W_k\}_{i \in I}$ and such that $\hat \delta_k \subset
\Int_{W_k} B_k$. If it is small enough, then it is closed in $X$. Let
$W'_k$ be an open subset of $V_k$ such that $W'_k \cap (A'_k \cap
F^{l,m+2}_j) = (A'_k \cap F^{l,m+2}_j) \setminus B_k$ and $\Cl_X W'_k
\cap \delta_k = \emptyset$. By lemma~\ref{lem:fitting approximation},
we may adjust $A'_k$ and $F^{l,m+2}_{j}$ $\rel X \setminus W'_k$ to
obtain sets $A''_k$ and $\widetilde F^{l, m+2}_{j}$ such that $\F{l,
  m+3}$ defined in the following way
\begin{equation*}
F^{l,m+3}_i = \left\{
    \begin{array}{ll}
      \widetilde F^{l, m+2}_{j} & i = j \\
      F^{l,m+2}_i & i \in I \text{ and } i \neq j.
    \end{array}\right.
\end{equation*}
is an adjustment of $\F{l, m+2}$. Observe that $A''_k$ fits $\F{l, m+3}$,
because its intersection with $\widetilde F^{l,m+2}_j$ is equal to $B_k$.
Observe that $\hat\delta_k \setminus \partial\hat\delta_k \subset \Int_X
A''_k$ and the restriction of $\F{l,m+3}$ to $A''_k$ is an $\N{n}$-cover. By
the construction, $\Cl_X \delta_k \setminus \hat\delta_k \subset \Int A''_k$.

Let
\begin{equation*}
F^{l,m+4}_i = \left\{
    \begin{array}{ll}
      F^{l,m+3}_i \cup A''_k & i \in J'_k \\
      F^{l,m+3}_i & i \in I \setminus J'_k.
    \end{array}\right.
\end{equation*}

Since $\delta_k \subset \Int_X F^{l,m+4}_j \cup \Int_X A''_k$, ($12_{k,\hat
  m}$) is satisfied.  Verification of other conditions is similar to the
verification done in the previous steps.

The limit argument is similar to the $0$-dimensional case. Let $\F{l, \infty}
= \lim_{k \to \infty} \F{l, k}$. We are going to show that $(1_{l,\infty})$,
$(2_{l,\infty})$, $(3_{l,\infty})$ and, for each at most $l$-dimensional
simplex $\delta \in \tau$ and each $i \in I$, $(4^{\delta,i}_{l,\infty})$ are
satisfied. Since the sequence $\F{l,0}, \F{l,1}, \ldots$ stabilizes locally
and $(1_{l,k})$, $(2_{l,k})$ and $(12_{k,m})$ are satsified for all $k$ and
$m$, all we have to show is that $\langle K, \mathcal{K} \rangle$ is an
$n$-core of $\F{l, \infty}$. Let $H_k \colon X \to X$ be a closed embedding
defined as the identity on $X \setminus \bigcup_{l \geq k} U_l$ and by the
formula $H_k(x) = h_m(x)$ for each $x \in U_m$ with $m \geq k$. Let $\G{k} =
\{ G^k_i \}_{i \in I}$ and let $G^k_i = H_k(F^{l, m_{k+1}}_i)$. Observe that
$\G{k}$ is a closed interior $\N{n}$-cover of $X \setminus \bigcup_{m \geq k}
V_m$ and $\langle L, \mathcal{L} \rangle$ is its $n$-core. Also, $\G{k}$ is
equal to $\G{k-1}$ on $X \setminus U_k$ and $\G{k-1}$ is a shrinking of
$\G{k}$ on $U_k \setminus V_k$. By lemma~\ref{lem:limit of deformations},
$\langle L, \mathcal{L} \rangle$ is an $n$-core of $\lim_{k \to \infty}
\G{k}$. But $\lim_{k \to \infty} \G{k} = \lim_{k \to \infty} \F{l, k}$, so
$\langle L, \mathcal{L} \rangle$, and also $\langle K, \mathcal{K} \rangle$,
is an $n$-core of $\F{l, \infty}$.

\subsection*{Second part.}

In the first part we constructed a closed interior $\N{n}$-cover
$\F{n, \infty}$ of $X$ that satisfies conditions ($1_{n,\infty}$),
($2_{n,\infty}$), ($3_{n,\infty}$) and ($4^{\delta,i}_{n,\infty}$) for
each $\delta \in \tau$ and $i \in I$. By the construction,
$\F{n,\infty}$ is equal to $\F{0,0}$ on the complement of an open
neighborhood $U_L$ of $L$. Let $\U{}$ be an open cover of $X$
satisfying conditions of remark~\ref{rem:core of a small adjustment}.
If $\F{n+1,\infty}$ is a $\U{}$-adjustment $\rel X\setminus U_L$ of
$\F{n,\infty}$ to a closed interior $\N{n}$-cover such that $K_i = K
\cap F^{n+1,\infty}_i$ for each $i \in I$, then, by the choice of
$\U{}$, $\langle K, \mathcal{K} \rangle$ is an exact $n$-core of
$\F{n+1,\infty}$. We shall construct such $\F{n+1,\infty}$, finishing
the proof.

Order $\tau_L$ into a sequence $\{ \delta_k \}_{k \geq 1}$ non
decreasing in the order by inclusion. Let $J_k = \{ i \in I \colon
\delta_k \subset \Int_K K_i \}$. By the choice of $\tau$ made at the
beginning of the proof, $J_k$ is non-empty for each $k \geq 1$. Since
$\F{n, \infty}$ satisfies ($4^{\delta_k, i}_{n, \infty}$) for each $k
\geq 1$ and each $i \in I$, the inclusion $\delta_k \subset \Int_X
F^{n, \infty}_{J_k}$ holds for each $k$ (this is the crucial property
obtained in the first part of the proof).  Therefore there exists a
locally finite collection $\{ V_k \}_{k \geq 1}$ of open subsets of
$U_L$ such that $L \subset \bigcup_{k \geq 1} V_k$ and the following
conditions are satisfied for each $k$.

\begin{enumerate}\addtocounter{enumi}{15}
\item $K \cap V_k \subset K_{J_k}$,
\item if $\Cl_X V_l$ intersects $\Cl_X V_k$, then either $\delta_l$ is a face
  of $\delta_k$ or $\delta_k$ is a face of $\delta_l$,
\item $\Cl_X V_k \subset \Int_X F^{n,\infty}_{J_k}$.
\end{enumerate}

Let $\F{n+1, 0} = \F{n, \infty}$. Let $\U{k}$ be a sequence of open
covers of $X$ obtained via lemma~\ref{pro:sequence of adjustments}
applied to a locally finite collection $\{ V_k \}_{k \geq 1}$ and an
open cover $\U{}$. We recursively define a sequence $\F{n+1, k}$ of
closed interior $\N{n}$-covers of~$X$ such that the following
conditions are satisfied for each $k \geq 1$.

\begin{enumerate}\addtocounter{enumi}{18}
\item $\F{n+1,k}$ is a $\U{k}$-adjustment of $\F{n+1,k-1}$ $\rel X
  \setminus V_k$ with fixed $J_k$,
\item $\F{n+1,k}$ is equal to $\mathcal{K}$ on $K \cap \bigcup_{l \leq k}
  V_l$.
\end{enumerate}

Fix $k \geq 1$ and assume that we already constructed $\F{n+1, k-1}$.
First, we claim that the inclusion $\Cl_X V_k \subset \Int_X
F^{n+1,l}_{J_k}$ holds for each $0 \leq l \leq k-1$. The proof is by
induction on $l$. The basis step of induction ($l = 0$) follows
directly from (18) and the identity $\F{n+1,0} = \F{n, \infty}$. Fix
$l > 0$ and assume that we already proved that $\Cl_X V_k \subset
\Int_X F^{n+1,l-1}_{J_k}$.  By~(19), $\F{n+1,l}$ is an adjustment
of~$\F{n,l-1}$ $\rel X \setminus V_l$, hence if $V_l$ is disjoint from
$Cl_X V_k$, then $\Cl_X V_k \cap \Int_X F^{n+1,l-1}_{J_k} = \Cl_X V_k
\cap \Int_X F^{n+1,l}_{J_k}$, hence $\Cl_X V_k \subset \Int_X
F^{n+1,l}_{J_k}$. If $\Cl_X V_l$ intersects $V_k$, then by (17) and by
the assumption that the sequence $\{ \delta_k \}$ is non-decreasing in
the order by inclusion, $\delta_l$ is a face of $\delta_k$. In this
case, $J_k$ is a subset of $J_l$ by the definition. By (19),
$\F{n+1,l}$ is an adjustment of $\F{n,l}$ with fixed $J_l$, hence
$\Int_X F^{n+1,l}_{J_k} = \Int_X F^{n+1, l-1}_{J_k}$ and the inclusion
$\Cl_X V_k \subset \Int_X F^{n+1,l}_{J_k}$ follows from the inductive
assumption. The inductive step is done. In particular, we have proven
that $\Cl_X V_k \subset \Int_X F^{n+1,k-1}_{J_k}$.

By lemma~\ref{lem:adjustment to a Z-collection}, there exists an
adjustment of a collection $\{ F^{n+1,k-1}_i \cap V_k \}_{i \in I
  \setminus J_k}$ of subsets of $V_k$ to a $Z$-collection that is
equal to $\{ K_i \cap V_k \}_{i \in I \setminus J_k}$ on $K \cap V_k$.
By lemma~\ref{lem:local adjustment}, we may require (and we do) that
the adjustment is so small, that it extends to an adjustment $\rel X
\setminus V_k$ of $\{ F^{n+1,k-1}_i \}_{i \in I \setminus J_k}$ to a
closed $\N{n}$-collection of subsets of $X$. We let $\F{n+1,k}$ be a
union of this collection with the collection $\{ F^{n+1,k-1} \}_{i \in
  J_k}$. Observe that $\F{n+1,k}$ is a $\U{k}$-adjustment of
$\F{n+1,k-1}$ $\rel X \setminus V_k$ with fixed $J_k$ directly from
the construction. It is a closed interior $\N{n}$-cover of $X$ because
$\Cl_X V_k \subset \Int_X F^{n+1,k}_{J_k}$. Hence (19) is satisfied.
To see that $\F{n+1,k}$ satisfies condition (20) recall that by (16),
$K \cap V_k \subset K_{J_k}$, hence $F^{n+1,k}_i \cap K \cap V_k = K_i
\cap V_k$ for each $i \in J_k$. If $i \in I \setminus J_k$, then
$F^{n+1,k}_i \cap K \cap V_k = K_i \cap V_k$ directly from the
construction.

Let $\F{n+1, \infty} = \lim_{k \to \infty} \F{n+1, k}$. By the choice
of $\U{k}$'s, it is a $\U{}$-adjustment of $\F{n+1, 0}$. Since the
collection $\{ V_k \}$ is locally finite, $\F{n+1,\infty}$ is a closed
interior $\N{n}$-cover. By the choice of $\U{}$, $\langle K,
\mathcal{K} \rangle$ is an $n$-core of $\F{n+1,\infty}$. By (20), by
the assumption that $L \subset \bigcup_{k \geq 1} V_k$ and by the
assumption that $\F{n+1, 0}$ is equal to $\mathcal{K}$ on $K \setminus
L$, $\langle K, \mathcal{K} \rangle$ is an exact $n$-core of $\F{n+1,
  \infty}$. The cover $\F{n+1,\infty}$ is equal to $\F{0,0}$ on the
complement of $U_L$ directly from the construction. We are done.

\section{Retraction onto a core and a proof of theorem~\lowercase{\ref{thm:retraction onto
    a complex}}}
\label{sec:retraction onto a core}

We shall now obtain theorem~\ref{thm:retraction onto a complex} as an
easy corollary of the results obtained earlier in this chapter.

\begin{lemma}
  \label{lem:retraction onto a core} If $K$ is an exact $n$-core of a
  closed star-finite locally finite $ANE(n)$-cover~$\F{}$ of an at most
  $n$-dimensional space $X$, then there exists a retraction of $X$ onto~$K$
  that is $\F{}$-close to the identity.
\end{lemma}
\begin{proof}
  Let $\F{} = \{ F_i \}_{i \in I}$ and for each non-empty $J \subset
  I$ let $F_J$ denote, as usual, the intersection $\bigcap_{i \in J}
  F_i$. Let $J_1, J_2$ be subsets of~$I$. Assume that $F_{J_2}$ is
  non-empty. We will show that the inclusion of $F_{J_2} \cap (K \cup
  \bigcup_{i \in J_1} F_i)$ into $F_{J_2}$ is a weak $n$-homotopy
  equivalence. Without a loss of generality, we may assume that $J_1$
  is disjoint from $J_2$ and that all elements indexed by $J_1$
  intersect $F_{J_2}$. By the definition of an exact $n$-core, for
  each $i \in I$ the intersection of $F_i$ with $K$ is a subcomplex of
  $K$. By the assumption that $\F{}$ is an $ANE(n)$-cover, for each
  non-empty $J \subset I$ the intersection $F_J$ is an $ANE(n)$-space.
  By the assumption, $\F{}$ is closed and locally finite, hence
  theorem~\ref{thm:sum theorem for ane-spaces} implies that $F_{J_2}
  \cap (K \cup \bigcup_{i \in J_1} F_i)$ is an $ANE(n)$-space for
  each~$J_1$ and~$J_2$.  Fix $J_2$ such that $F_{J_2}$ is non-empty
  and assume that the assertion is already proven for all $J \subset
  I$ such that $J_2 \subsetneq J$ and $F_J \neq \emptyset$ (there are
  only finitely many of such $J$'s as $\F{}$ is star-finite). Let
  $J_1$ be a subset of $I$ that is disjoint from $J_2$ and is such
  that elements indexed by $J_1$ intersect $F_{J_2}$.  Again, there
  are only finitely many of those. We shall prove the assertion by an
  induction over the cardinality of $J_1$.  By the assumption that $K$
  is an exact $n$-core of $\F{}$, the inclusion of $K \cap F_{J_2}$
  into $F_{J_2}$ is an $n$-homotopy equivalence. By Whitehead's
  characterization, this inclusion is a weak $n$-homotopy equivalence.
  Hence the assertion is verified for $J_1 = \emptyset$.  Let $J_1
  \neq \emptyset$. Fix $j \in J_1$. We assumed that $J_1$ is disjoint
  from $J_2$, so $j \not\in J_2$. Let $A_1 = F_{J_2} \cap (K \cup
  \bigcup_{i \in J_1, i \neq j} F_i)$ and $A_2 = F_{J_2} \cap F_j$. As
  we argued earlier, both $A_1$ and $A_2$ are $ANE(n)$-spaces. By the
  assumption, $X$ is at most $n$-dimensional, and so are $A_1$ and
  $A_2$. By the inductive assumptions, the inclusions $A_1 \cap A_2
  \subset A_2$ and $A_1 \subset F_{J_2}$ are weak $n$-homotopy
  equivalences. Hence by corollary~\ref{cor:excision}, the inclusion
  $A_1 \cup A_2 \subset F_{J_2}$ is a weak $n$-homotopy equivalence.
  This is the assertion that we wanted to prove.

  By theorem~\ref{thm:subspace retraction}, for each $J \subset I$
  such that $F_J$ is non-empty there exists a retraction $r_J \colon
  F_J \to F_J \cap (K \cup \bigcup_{i \in I \setminus J} F_i)$. Let
  $\mathcal{I} = \{ J \subset I \colon F_J \neq \emptyset \}$. Order
  $\mathcal{I}$ into a sequence $\{ J_k \}_{k \geq 1}$ non decreasing
  in the order by inclusion. Such ordering exists because $\F{}$ is
  point finite. Let $r_0$ be the identity on $X$. Let $r_k \colon
  \bigcup_{i \in I \setminus J_k} F_i \cup K \cup F_{J_k} \to
  \bigcup_{i \in I \setminus J_k} F_i \cup K$ be a map defined by the
  following equation.
  \[
    r_k = \left\{ 
    \begin{array}{ll}
      r_{J_k}(x) & x \in F_J \\
      x & x \in K \cup \bigcup_{i \in I \setminus J_k} F_i.
    \end{array}
    \right.
  \]
  Obviously each $r_k$ is a contiuous retraction and $r_k(F_i) \subset
  F_i$ for each $i \in I$. By the order of $J_k$'s for each $x$ in $X$
  the sequence of compositions $r_1, r_2 \circ r_1, \ldots$ is well
  defined. By star-finiteness of $\F{}$ it stabilizes for each $x$
  after finitely many steps. Hence the limit
  \[
    r(x) = \lim_{k \to \infty} r_k \circ \cdots \circ r_2 \circ r_1(x).
  \]
  is well defined for each $x$ in $X$. It is a retraction that we were
  looking for.
\end{proof}

\begin{proof}[Proof of theorem~\ref{thm:retraction onto a complex}]
  Let $\U{}$ be an open cover of an $\N{n}$-space~$X$.
  By~\cite[theorem 5.3.10]{engelking1989} and the fact that every
  separable metric space is strongly paracompact, there exists a
  closed star-finite locally finite cover of $X$ that refines $\U{}$.
  By theorem~\ref{thm:n-swelling}, it has a closed star-finite
  swelling to an interior $\N{n}$-cover $\F{}$ refining $\U{}$. By
  theorem~\ref{thm:embedded core}, $\F{}$ has an $n$-core $\langle K,
  \mathcal{K} \rangle$. By theorem~\ref{thm:construction of an exact
    core}, $\F{}$ can be deformed to a closed star-finite interior
  $\N{n}$-cover $\G{}$, which refines $\U{}$ and has an exact
  $n$-core~$K$. By lemma~\ref{lem:retraction onto a core}, there is a
  retraction $r \colon X \to K$ that is $\G{}$-close to the identity.
  It satisfies assertion of the theorem because $\G{}$ refines~$\U{}$.
\end{proof}


\chapter{Proof of theorem~\lowercase{\ref{thm:pump up the regularity}}}
\label{ch:the existence of regular covers}

The entire chapter is devoted to the proof of theorem~\ref{thm:pump up the
  regularity}. Let us recall the statement.

\begin{pump up theorem}
  There exists a constant $N$ such that if a $k$-regular
  $n$-semiregular ($k < n$) closed countable star-finite interior
  $\N{n}$-cover $\F{}$ is $n$-contractible in a closed cover $\E{}$,
  then there exists a $(k+1)$-regular $n$-semiregular interior
  $\N{n}$-cover isomorphic to $\F{}$ and refining $\st^N \E{}$.
  Moreover, we may additionally require that the constructed cover is
  equal to $\F{}$ on some open neighborhood of a given $Z(\F{})$-set
  $Z$.
\end{pump up theorem}

We shall prove that the assertion is true with $N = 472$, but the
actual value is of no importance and no effort was made to minimize
it.

\section{Patching of small holes}
\label{ssec:patching of small holes}

We fix $n$ and $k < n$. Let $\F{} = \{ F_i \}_{i \in I}$ be a cover of
a space $X$. We say that a map~$\varphi$ from a $k$-dimensional sphere
$S^k$ into $X$ is \dfi{a $k$-hole in $\F{}$} if the image of $\varphi$
lies in $F_J = \bigcap_{i \in J} F_i$ for some non-empty subset $J$
of~$I$. When we will want to specify $J$, we will say that $\varphi$
is \df{hole!in a cover}{a $k$-hole in~$F_J$}. If~$\varphi$ is a
$k$-hole in $F_J$ and is null-homotopic in $F_J$, then we say that it
is a \df{hole!trivial}{trivial $k$-hole}.  Observe that if $\F{}$ is
$k$-regular and all $k$-holes in $\F{}$ are trivial, then it is
$(k+1)$-regular.  We shall describe a construction that enables us to
\emph{patch} holes.  That is, to make them trivial in a newly
constructed cover. We begin with a special case. If $\varphi$ is a
$k$-hole in $F_J$ and $\varphi$ is contractible in $F_j$ for some $j
\in J$, then we say that $\varphi$ is a \df{hole!small}{small
  $k$-hole}. We describe a construction of patching of a small
$k$-hole.

\begin{definition} \label{def:patching of holes} Let $\langle K,
  \mathcal{K} \rangle$ be a pair consisting of a locally finite
  countable simplicial complex~$K$ and a cover $\mathcal{K} = \{ K_i
  \}_{i \in I}$ of~$K$ by subcomplexes.  Let $k < n$ and let $B$ be a
  $(k+1)$-dimensional ball that is embedded as a subcomplex of~$K$. We
  will say that a pair $\langle L, \mathcal{L} \rangle$ is
  \df{obtained by patching}{obtained from $\langle K, \mathcal{K}
    \rangle$ by patching of $B$}, if the following conditions are
  satisfied.

  \begin{enumerate}
  \item 
    \begin{enumerate}
    \item $L$ is a locally finite countable simplicial complex and $K$
      is a subcomplex of $L$,
    \item $\mathcal{L} = \{ L_i \}_{i \in I}$ is a cover of $L$ by
      subcomplexes and $\mathcal{L}$ is equal to $\mathcal{K}$ on $K$.
    \end{enumerate}
  \item There exists a complex 
    \[
    C = \left\{ \begin{array}{ll}
        B \times [0, 1], & \text{if } k < n - 1 \\
        B \times \{ 0, 1\} \cup \partial B \times [0, 1], & \text{if }
        k = n - 1
      \end{array} \right.
    \]
    such that $L = K \cup C$ and the intersection of $K$ with $C$ is
    equal to $B$ (in $K$) and to $B \times \{ 0 \}$ (in $C$). We
    identify $B$ with $B \times \{ 0 \}$.
  \item Recall that $\mathcal{K} = \{ K_i \}_{i \in I}$ and
    $\mathcal{L} = \{ L_i \}_{i \in I}$. For each $i \in I$, $L_i$ is
    the smallest subset of $L$ for which the following conditions are
    satisfied.
    \begin{enumerate}
    \item if $B \subset K_i$, then $C \subset L_i$.
    \item if $\partial B \subset K_i$, then $\partial B \times [0, 1]
      \cup B \times \{ 1 \} \subset L_i$.
    \item $\begin{array}{l}
        \text{if } k < n - 1, \text{then } (B \cap K_i) \times [0, 1/2] \subset L_i\\
        \text{if } k = n - 1, \text{then } (\partial B \cap K_i) \times [0, 1/2] \subset L_i
      \end{array}
      $.
    \end{enumerate}
    Note that these conditions uniquely determine $L_i$'s.
  \end{enumerate}
\end{definition}

{%
\begin{figure}[ht]
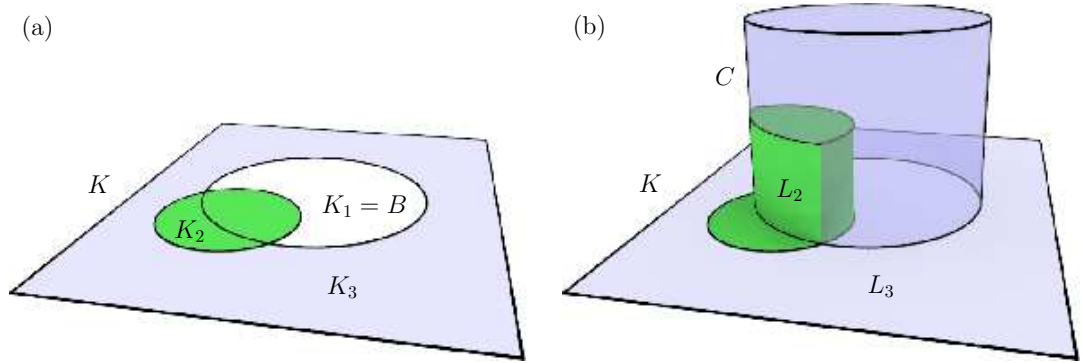

\begin{center}
$\begin{array}{c@{\hspace{5mm}}c}
\resizebox{0.45\textwidth}{!}{\input{psfigures/small_hole1.pstex_t}} &
\resizebox{0.45\textwidth}{!}{\input{psfigures/small_hole2.pstex_t}}
\end{array}$
\end{center}
\caption{{A picture discussed in
  example~\ref{ex:patching of holes}}}
\label{fig:patching in cores}
\end{figure}
}

\begin{example}\label{ex:patching of holes}
  Let $K$ be a square shown on figure~\ref{fig:patching in cores}(a)
  and let $B$ be the white disc embedded as a subcomplex of $K$. We
  let $K_1 = B$, $K_2$ to be the dark disc and $K_3$ to be the
  complement of the interior of $K_1$. Let $\mathcal{K} = \{ K_1, K_2,
  K_3 \}$.  Figure~\ref{fig:patching in cores}(b) shows a pair
  $\langle L, \mathcal{L} \rangle$ obtained from $\langle K,
  \mathcal{K} \rangle$ by patching of $B$. The cylinder $C = B \times
  [0, 1]$ is attached to $K$ along its lower base (in the example we
  assume that $k < n - 1$) to form $L = K \cup C$. The set $L_1$ is
  equal to $C$.  The set $L_2$ is an union of $K_2$ with a cylinder of
  height $1/2$ and base $K_1 \cap K_2$. The set $L_3$ is an union of
  $K_3$ with $\partial B \times [ 0, 1 ] \cup B \times \{ 1 \}$.

  The cover $\mathcal{K}$ of $K$ is $1$-regular, but not $2$-regular
  as both $K_3$ and $K_1 \cap K_3$ are not simply connected. Observe
  that $\partial B \subset K_1 \cap K_3$ is a small hole in $K_1 \cap
  K_3$, as it is null-homotopic in $K_1$. Observe that $\partial B
  \subset L_1 \cap L_3$ is null-homotopic in $L_1 \cap L_3$. The small
  hole $\partial B \subset K_1 \cap K_3$ in $\mathcal{K}$ is
  ``patched'' in $\mathcal{L}$ and the cover $\mathcal{L}$ is
  $2$-regular.
\end{example}

\begin{remark}\label{rem:L is an interior cover}
  Assume the notation of definition~\ref{def:patching of holes}. By
  (1b), $\mathcal{L}$ is a cover of $L$ and by (3), $\mathcal{L}$
  covers $L$ if and only if $B \subset K_i$ for some $i \in I$. It
  follows from (3c) that if $\mathcal{K}$ is an interior cover, then
  $\mathcal{L}$ is an interior cover as well. This is the only reason
  for which (3c) is required. By the conditions enforced on
  $\mathcal{L}$, $\mathcal{L}$ is always isomorphic to $\mathcal{K}$.
\end{remark}

\begin{definition}\label{def:patching of holes in covers}
  Let $\F{} = \{ F_i \}_{i \in I}$ be a closed star-finite interior
  $\N{n}$-cover of a space $X$. Let $K$ be an exact $n$-core of
  $\F{}$. Let~$B$ be a subcomplex of~$K$ that is homeomorphic to a
  $(k+1)$-dimensional ball and such that $B \subset F_i$ for some $i
  \in I$. Let $U$ be an open neighborhood of $B$ in $X$. We say that a
  closed star-finite interior $\N{n}$-cover $\G{} = \{ G_i \}_{i \in
    I}$ of~$X$ is \df{cover!obtained by patching of $B$}{obtained from
    $\F{}$ by patching of $B$ in $U$}, if the conditions (0)-(3)
  stated below are satisfied.

  \begin{enumerate}\addtocounter{enumi}{-1}
    \item
    \begin{enumerate}
    \item $\G{}$ is isomorphic to $\F{}$ and is equal to $\F{}$ on $K
      \cup (X \setminus U)$,
    \item there exists an exact $n$-core $L$ of $\G{}$ such that $L
      \setminus K \subset U$.
    \end{enumerate}
  \end{enumerate}

  Recall that $\F{} / K = \{ F_i \cap K \}_{i \in I}$ and $\G{} / L =
  \{ G_i \cap L \}_{i \in I}$ (definition~\ref{def:restricted cover}).
  We require that the pair $\langle L, \G{} / L \rangle$ is obtained
  from $\langle K, \F{} / K \rangle$ by patching of $B$, i.e. that

  \begin{enumerate}
  \item 
    \begin{enumerate}
    \item $L$ is a locally finite countable simplicial complex and $K$
      is a subcomplex of $L$,
    \item $\mathcal{G} / L$ is a cover of $L$ by subcomplexes and
      $\mathcal{F}$ is equal to $\mathcal{G}$ on~$K$.
    \end{enumerate}
  \item There exists a complex 
    \[
    C = \left\{ \begin{array}{ll}
        B \times [0, 1], & \text{if } k < n - 1 \\
        B \times \{ 0, 1\} \cup \partial B \times [0, 1], & \text{if }
        k = n - 1
      \end{array} \right.
    \]
    such that $L = K \cup C$ and the intersection of $K$ with $C$ is
    equal to $B$ (in $K$) and to $B \times \{ 0 \}$ (in $C$). We
    identify $B$ and $B \times \{ 0 \}$.
  \item For each $i \in I$, $G_i \cap L$ is the smallest subset of $L$
    for which the following conditions are satisfied.
    \begin{enumerate}
    \item if $B \subset F_i$, then $C \subset G_i$.
    \item if $\partial B \subset F_i$, then $\partial B \times [0, 1]
      \cup B \times \{ 1 \} \subset G_i$.
    \item $\begin{array}{l}
        \text{if } k < n - 1, \text{then } (B \cap F_i) \times [0, 1/2] \subset G_i\\
        \text{if } k = n - 1, \text{then } (\partial B \cap F_i) \times [0, 1/2] \subset G_i
      \end{array}
      $.
    \end{enumerate}
    Note that these conditions uniquely determine $L\cap G_i$'s.
  \end{enumerate}
\end{definition}

\begin{remark}
  Assume the notation of definition~\ref{def:patching of holes in
    covers}.  Let $J \subset I$ and assume that $B \subset F_i$ for
  some $i \in J$. If $\partial B \subset F_J$ is a small $k$-hole in
  $F_J$, then $\partial B \subset G_J$ is a trivial $k$-hole in $G_J$.
  There are no "new" non-trivial $k$-holes in $\G{}$ (cf.
  example~\ref{ex:patching of holes}; $\langle L, \G{} /
  L \rangle$ is an exact $n$-core of $\G{}$).
\end{remark}

\begin{lemma}\label{lem:patching of a small hole}
  If $\F{}$ is a closed star-finite interior $\N{n}$-cover of a
  space~$X$, $K$ is an exact $n$-core of $\F{}$, $B$ is a
  $(k+1)$-dimensional ball embedded as a subcomplex of $K$, $B$ is
  contained in some element of $\F{}$ and $U$ is a neighborhood of
  $B$, then there exists a closed star-finite interior $\N{n}$-cover
  $\G{}$ of $X$ that is obtained from $\F{}$ by patching of~$B$
  in~$U$.
\end{lemma}
\begin{proof}
  Let $\F{} = \{ F_i \}_{i \in I}$. The construction of $\G{}$
  proceeds in two stages.

  \vspace{ 1mm}
  \noindent \emph{Cleaning up.} In the first stage we construct a
  closed star-finite interior $\N{n}$-cover $\mathcal{H} = \{ H_i
  \}_{i \in I}$ of $X$ and its exact $n$-core $M$ that satisfy
  conditions of definition~\ref{def:patching of holes in covers} (with
  $H_i$ substituted for $G_i$ and $M$ substituted for $L$), with the
  exception that condition (3b) is replaced by condition

  \begin{enumerate}
  \item[(3b')] if $\partial B \subset F_i$, then $\partial B \times
    [0, 1] \subset H_i$.
  \end{enumerate}

  Let
  \[
    C_0 = \left\{ \begin{array}{ll}
      B \times [0, 1], & \text{if } k < n - 1, \\
      B \times \{ 0, 1 \} \cup \partial B \times [0, 1], & \text{if } k = n - 1
    \end{array}\right..
  \]

  Identify $B \times \{ 0 \}$ with $B$ and let $p$ be the projection
  of $C_0$ onto $B$ along $[0, 1]$. Let $\imath$ be the inclusion of
  $K$ into $X$. The map $p \cup \imath$ is a well-defined map into
  $X$, and its domain is an most $n$-dimensional countable simplicial
  complex. By theorem~\ref{thm:approximation within a cover} and by
  remark~\ref{rem:closed locally compact is ZF}, we can approximate $p
  \cup \imath$ $\rel K$ within $\F{}$ by a closed embedding $f$. Take
  an approximation close enough, such that the image $f(C_0)$ lies in
  $U$ and let~$C = f(C_0)$. Observe that $C \cap K = B$. Let
  $\mathcal{K} = \{ K_i \}_{i \in I}$ be a cover of $K \cup C$ by
  subcomplexes, such that $\mathcal{K}$ is equal to $\F{}$ on $K$ and
  such that for each $i \in I$, $K_i$ is a smallest subset of $K$ for
  which the following conditions are satisfied.
  \begin{enumerate}
  \item[(a)] if $B \subset F_i$, then $C \subset K_i$.
  \item[(b')] if $\partial B \subset F_i$, then $\partial B \times [0, 1]
    \subset K_i$.
  \item[(c)] $\begin{array}{l}
      \text{if } k < n - 1, \text{then } (B \cap F_i) \times [0, 1/2] \subset K_i\\
      \text{if } k = n - 1, \text{then } (\partial B \cap F_i) \times
      [0, 1/2] \subset K_i
    \end{array}
    $.
  \end{enumerate}
  Observe that $\mathcal{K}$ refines $\F{}$ because $f$ was an
  approximation within $\F{}$. By the assumption $K$ is an exact
  $n$-core of $\F{}$, so $\langle K \cup C, \mathcal{K} \rangle$ is an
  $n$-core of $\F{}$.  By theorem~\ref{thm:construction of an exact
    core} applied to $\langle K \cup C, \mathcal{K} \rangle$ and
  $\F{}$ there is a deformation of $\F{}$ $\rel X \setminus U$ to an
  interior $\N{n}$-cover $\mathcal{H}$ with an exact $n$-core $\langle
  K \cup C, \mathcal{K} \rangle$. By theorem~\ref{thm:construction of
    an exact core}, we may require (and we do), that if $B \times \{ 1
  \}$ is a subset of $K_i$, then it is a subset of $\Int_X H_i$. The
  cover $\mathcal{H}$ satisfies conditions of
  definition~\ref{def:patching of holes in covers} (with $H_i$
  substituted for $G_i$ and $M$ substituted for $L$), with the
  exception that condition (3b) is replaced by condition (3b').

  \vspace{ 1mm}
  \noindent \emph{Patching.}  Let $\Delta = B \times \{ 1 \}$ be the
  upper base of $C$. Let $J = \{ i \in I \colon B \subset F_i \}$ and
  $J_\partial = \{ i \in I \colon \partial B \subset F_i \}$. By the
  construction, $\Delta$ is a subset of the interior of $H_J$ and is
  disjoint from $H_i$ for each $i \in I \setminus J_\partial$. Hence
  there exists an open neighborhood $V$ of $\Delta$ that is disjoint
  from $\bigcup_{i \in I \setminus J_\partial} H_i$ and whose closure
  lies in the intersection of $U$ with the interior of $H_J$. By
  theorem~\ref{thm:z-approximation} and by lemma~\ref{lem:local
    adjustment} there exists an approximation $\rel C$ of the
  inclusion $V \cap \bigcup_{i \in J_\partial \setminus J } H_i
  \subset V$ by a $Z$-embedding $h$ into $V$ such that $id_{X
    \setminus V} \cup h$ is continuous. Let $\mathcal{H}' = \{ H'_i
  \}_{i \in I}$, where
  \[
  H'_i = \left\{ \begin{array}{ll}
      H_i & \text{if } i \in J \\
      (H_i \setminus V) \cup h(V \cap H_i) & \text{if } i \in I \setminus J.
    \end{array}\right.
  \]
  Observe that $\mathcal{H}'$ satisfies conditions (0) - (3ab'c) with
  $H'_i$ substituted for $G_i$ and $V \cap \bigcup_{i \in I \setminus
    J} H'_i$ is a $Z$-set in $V$.

  By lemma~\ref{lem:disc attaching} applied to $\Delta$, a $Z$-set
  $\tilde B = H'_{J_\partial} \cap V$ and a space~$V$, there exists a
  closed $\N{n}$-subset $A$ of $V$ such that $A \cap H'_{J_\partial}$
  is an $\N{n}$-space, $\Delta \subset A$ and the inclusions $\partial
  \Delta \subset A \cap \tilde B$, $\Delta \subset A$ are $n$-homotopy
  equivalences. We may assume that $A$ is so small that it is closed
  in $X$. By corollary~\ref{cor:z-approximation}, the inclusion $A
  \cup H'_{J_\partial} \subset V$ can be approximated $\rel \Delta
  \cup H'_{J_\partial}$ by a closed embedding with image disjoint from
  $D = (\bigcup_{i \in I \setminus J} H'_i \cup C) \setminus (\Delta
  \cup H'_{J_\partial})$. Let $A'$ denote the image of $A$ under this
  embedding.

  Let $\G{} = \{ G_i \}_{i \in I}$, with $G_i = H'_i$ for $i \in I
  \setminus J_\partial$ and $G_i = H'_i \cup A'$ for $i \in
  J_\partial$. By the choice of $V$, $A' \subset H'_i$ for each $i \in
  J$ and by the construction $H'_i \cap A' = H'_J \cap A$ for each $i
  \in J_\partial \setminus J$ and the latter intersection is an
  $\N{n}$-space, so $\G{}$ is a closed interior $\N{n}$-cover by
  corollary~\ref{cor:pasting}.  By the choice of $V$ condition (0a) is
  satisfied.  Let $L = K \cup C$. Since $A' \cap L = \Delta$,
  conditions (1) - (3abc) are satisfied.  What is left is to show that
  $L$ is an exact $n$-core of $\G{}$, i.e. that (0b) holds.

  Consider a set $J' \subset I$. We have to show that $L \cap G_{J'}
  \subset G_{J'}$ is an $n$-homotopy equivalence. If $J' \not\subset
  J_\partial$, then $G_{J'} = H'_{J'}$ and $L \cap G_{J'} = L \cap
  H'_{J'}$.  By the construction $L$ is an exact $n$-core of
  $\mathcal{H}'$, so we are done. If $J' \subset J_\partial$ then $L
  \cap G_{J'} = (L \cap H'_{J'}) \cup \Delta$ and $G_{J'} = H'_{J'}
  \cup A'$.  We shall show that inclusions $(L \cap H'_{J'}) \cup
  \Delta \subset H'_{J'} \cup \Delta$ and $H'_{J'} \cup \Delta \subset
  H'_{J'} \cup A'$ are $n$-homotopy equivalences. The latter inclusion
  is an $n$-homotopy equivalence by remark~\ref{rem:disc attaching}.
  The former inclusion is an $n$-homotopy equivalence by
  lemma~\ref{lem:excision} applied with $A_1 = (L \cap H'_{J'}) \cup
  \Delta$ and $A_2 = H'_{J'}$. We are done.
\end{proof}

\section{$\circledcirc$-Contractibility}

Throughout the section we fix a positive integer $n$ and a
nonnegative integer $k < n$. We let $B^{k+1} \subset
\mathbb{R}^{k+1}$ be the $(k+1)$-dimensional unit ball.

\begin{definition*}
  For each positive integer~$l$ and each positive $m \leq l$ we
  let~\df{Aml@$A^m_l$}{$A^m_l$} denote the annulus $\{ x \in B^{k+1}
  \colon \frac{m-1}{l} \leq |x| \leq \frac{m}{l} \}$.
\end{definition*}

\begin{definition}\label{def:o-contractible}
  Let $\F{}$ and $\G{}$ be covers of a space $X$. Let $\varphi \colon
  S^k \to X$ be a $k$-hole in~$\F{}$. We say that $\varphi$ is
  \df{map!o-contractible@$\circledcirc_\F{}$-contractible}{$\circledcirc_\F{}$-contractible
    in~$\G{}$} (pronounced ``o''-$\F{}$-contractible), if it extends to
  a map~$\varPhi$ from $B^{k+1}$ into an element of $\G{}$ such that
  for some positive integer~$l$ the collection of annuli $\{ A^m_l
  \}_{1 \leq m \leq l}$ refines $\{ \varPhi^{-1}(F) \}_{F \in \F{}}$.
  We shall call $\varPhi$ an
  \df{o-contraction@$\circledcirc_\F{}$-contraction}{$\circledcirc_\F{}$-contraction
    of $\varphi$}.
\end{definition}

\begin{remark}\label{rem:interior cover by annuli}
  If $\varPhi$ is an $\circledcirc_\F{}$-contraction of a map $\varphi
  \colon S^k \to X$, then for some positive integer~$l$ the collection
  $\{ \bigcup_{j \leq l, \varPhi(A^j_l) \subset F_i} A^j_l \}_{i \in
    I}$ is a cover of $B^{k+1}$ by subcomplexes.  We can always
  replace $\varPhi$ by another $\circledcirc_\F{}$-contraction such
  that this collection (possibly with different $l$) is an interior
  cover of $B^{k+1}$. The new $\circledcirc_\F{}$-contraction can be
  obtained by a reparametrization of $\varPhi$. In particular, we may
  assume that its image is equal to the image of $\varPhi$.
\end{remark}

\begin{remark}
  In definition~\ref{def:o-contractible}, the sole role of a cover
  $\G{}$ is to limit the size of the image of $\varPhi$. The role of
  $\F{}$ is more interesting and may be interpreted in terms of small
  holes as follows. Let $\F{}$ be a cover of a space $X$ and let
  $\varphi$ be a $k$-hole in $\F{}$. Assume that $\varphi$ is
  $\circledcirc_\F{}$-contractible and let $\varPhi$ be a map
  satisfying conditions of definition~\ref{def:o-contractible} for a
  positive integer $l$. Let $\varphi_m$ be a map from $S^k$ into $X$
  defined by the formula $\varphi_m(x) = \varPhi(mx/l)$. The hole
  $\varphi_1$ is small. The hole $\varphi_{2}$ might not be small, but
  if we patch $\varphi_1$ (a precise meaning of it will be given in
  the proof of lemma~\ref{lem:ok contractible}), then it would become
  small. Hence we may think that $\varphi_{2}$ is a ``small hole of
  second order''. Pushing this argument further, we see that
  $\varphi_l$, which is equal to $\varphi$, is a ``small hole of $l$th
  order'' and that it can be patched after recursively patching of
  $\varphi_1, \varphi_2, \ldots$ up to $\varphi_{l-1}$.
\end{remark}

\begin{definition*}
  We say that a cover $\F{} = \{ F_i \}_{i \in I}$ is
  \df{cover!$ok-contractible@$\circledcirc_k$-contractible}{$\circledcirc_k$-contractible
    in a cover $\G{}$} if it refines $\G{}$ and if for every non-empty
  subset $J$ of $I$ such that $F_J = \bigcap_{i \in J} F_i$ is
  non-empty, the $k$th homotopy group $\pi_k(F_J)$ has a set of
  generators that are $\circledcirc_\F{}$-contractible in~$\G{}$, for
  some choice of base-points of $F_J$'s.
\end{definition*}

\begin{lemma}
  \label{lem:ok contractible} If $\F{}$ is a $k$-regular ($k < n$)
  closed star-finite interior $\N{n}$-cover that is
  $\circledcirc_k$-contractible in a closed cover $\G{}$, then there
  exists a closed $(k+1)$-regular interior $\N{n}$-cover isomorphic
  to~$\F{}$ and refining $\st_{\st^3_\F{} \G{}} \F{}$. Moreover, given
  a $Z(\F{})$-set~$Z$, we may require that there exists a neighborhood
  of $Z$ on which the constructed cover is equal to~$\F{}$.
\end{lemma}
\begin{proof}
  Let $\F{} = \{ F_i \}_{i \in I}$ be a $k$-regular ($k < n$) closed
  interior $\N{n}$-cover of a space $X$ that is
  $\circledcirc_k$-contractible in a closed cover $\G{}$. Let
  $\mathcal{J} = \{ J \subset I \colon J \neq \emptyset \text{ and }
  F_J \neq \emptyset \}$.  By proposition~\ref{pro:countable homotopy
    groups}, for every $J \in \mathcal{J}$ the $k$-dimensional
  homotopy group $\pi_k(F_J)$ is countable. Therefore for every $J \in
  \mathcal{J}$ there exists a collection of maps $\{ \psi^J_m \colon
  S^k \to F_J \}_{m \in \mathbb{N}}$ that generates $\pi_k(F_J)$. By
  the assumption that $\F{}$ is $\circledcirc_k$-contractible, we may
  additionaly require that $\psi^J_m$'s are
  $\circledcirc_\F{}$-contractible in $\G{}$, i.e.  that for every $J
  \in \mathcal{J}$ and every $m \in \mathbb{N}$ there exists an
  element $G^J_m$ of~$\G{}$ and a $\circledcirc_\F{}$-contraction
  $\Psi^J_m \colon B^{k+1} \to G^J_m$ of $\psi^J_m$. By
  remark~\ref{rem:interior cover by annuli}, we may assume that the
  following condition is satisfied.

  \begin{enumerate}
  \item There exists an integer $l^J_m$ such that the collection $\{
    \bigcup \{ A^j_{l^J_m} \colon 1 \leq j \leq l^J_m,
    \Psi^J_m(A^j_{l^J_m}) \subset F_i \} \}_{i \in I}$ is an interior
    cover of $B^{k+1}$.
  \end{enumerate}

  By the assumptions, $\F{}$ is a countable star-finite interior
  cover, so it is locally finite. Hence every element of $\mathcal{J}$
  is a finite subset of~$I$. In particular, $\mathcal{J}$ is
  countable.  We endow $\mathcal{J}$ with the discrete topology. Let
  $\Psi \colon B^{k+1} \times \mathcal{J} \times \mathbb{N} \to X$ be
  a map defined by the formula $\Psi(x, J, m) = \Psi^J_m(x)$, i.e. let
  $\Psi$ be a disjoint union of $\Psi^J_m$'s.  The domain of $\Psi$ is
  a $k$-dimensional locally finite countable simplicial complex.  Let
  $\U{}$ be an open cover of $X$ obtained via lemma~\ref{lem:close
    approximation} applied to~$\F{}$. By
  theorem~\ref{thm:approximation within a cover}, there exists a
  $\U{}$-approximation of $\Psi$ within~$\F{}$ by a closed embedding
  $\varPhi$ with image disjoint from~$Z$. Let $\varPhi^J_m \colon
  B^{k+1} \to X$ be a map defined by the formula $\varPhi^J_m(x) =
  \varPhi(x,J,m)$, i.e. an approximation of $\Psi^J_m$ obtained as a
  restriction of $\Phi$ to an appropriate subset of its domain.
  Observe that $\varPhi^J_m$'s satisfy the following conditions.

  \begin{enumerate}\addtocounter{enumi}{1}
  \item $\varphi^J_m = {\varPhi^J_m}_{|S^k}$ maps $S^k$ into $F_J$ and
    $\{ \varphi^J_m \}_{m \in \mathbb{N}}$ generates $\pi_k(F_J)$.
  \item There exists an integer $l^J_m$ such that the collection $\{
    \bigcup \{ A^j_{l^J_m} \colon 1 \leq j \leq l^J_m,
    \varPhi^J_m(A^j_{l^J_m}) \subset F_i \} \}_{i \in I}$ is an
    interior cover of $B^{k+1}$.
  \item Each $\varPhi^J_m$ is an embedding, the collection $\{ \im
    \varPhi^J_m \}_{J \in \mathcal{J}, m \in \mathbb{N}}$ is a discrete
    collection in $X$ and $\im \varPhi^J_m \subset \st_\F{} G^J_m$.
  \end{enumerate}

  Condition (2) is satisfied because $\varphi^J_m$ is a
  $\U{}$-approximation within~$\F{}$ of $\psi^J_m$; by the choice of
  $\U{}$, every two $\U{}$-close maps into $F_J$ are homotopic and by
  the choice of $\psi^J_m$'s, $\pi_k(F_J)$ is generated by $\{
  \psi^J_m \}_{m \in \mathbb{N}}$.  Condition (3) is satisfied because
  $\varphi^J_m$ is an approximation within $\F{}$ of $\psi^J_m$, so
  the collection defined in (1) refines the collection defined in (3).
  Condition (4) is satisfied because $\varPhi$ is a closed embedding
  into $X$ and is $\F{}$-close to $\Psi$.

  Let $L = \bigcup_{J \in \mathcal{J}, m \in \mathbb{N}} \im
  \varPhi^J_m$.  It is a locally finite simplicial complex embedded as
  a closed subset of $X$. By the construction it is disjoint from $Z$.
  For each $i \in I$, let
  \[ 
  L_i = \bigcup \{ \varPhi^J_m(A^j_{l^J_m}) \colon J \in \mathcal{J}, m \in
  \mathbb{N}, 1 \leq j \leq l^J_m, \varPhi^J_m(A^j_{l^J_m}) \subset
  F_i \}.
  \]
  By (3), $\mathcal{L} = \{ L_i \}_{i \in I}$ is an interior cover of
  $L$ by subcomplexes. By theorem~\ref{thm:embedded core}, we may
  enlarge the pair $\langle L, \mathcal{L} \rangle$ to an $n$-core
  $\langle K_0, \mathcal{K}_0 \rangle$ of $\F{}$, with $\mathcal{K}_0 =
  \{ K^0_i \}_{i \in I}$ and $K_0$ disjoint from $Z$. For each $J \in
  \mathcal{J}$ and $m \in \mathbb{N}$ fix a triangulation of a copy of
  $B^{k+1}$ such that $A^j_{l^J_m}$ is a subcomplex of $B^{k+1}$ for
  each $1 \leq j \leq l^J_m$. To avoid undue proliferation of
  notation, each copy is denoted by the same symbol $B^{k+1}$.  Let
  $p^J_m \colon C(B^{k+1}) \to B^{k+1}$ be a projection along $[0,1]$.
  By theorem~\ref{thm:approximation within a cover}, we may
  approximate each $\widetilde \varPhi^J_m$ within $\F{}$ $\rel
  B^{k+1} \times \{ 0 \}$ by a closed embedding $\widehat \varPhi^J_m$
  with image disjoint from~$Z$. We may require that images of
  $\widehat \varPhi^J_m$'s form a discrete collection in $X$ (by
  simultaneously approximating the disjoint union of $\varPhi^J_m$'s,
  like we did when approximating $\Psi^J_m$'s).  Since the
  approximation is within $\F{}$, $\im \widehat \varPhi^J_m \subset
  \st_\F{} G^J_m$ and $\widehat \varPhi^J_m({\varPhi^J_m}^{-1}(K^0_i))
  \subset F_i$ for each $i \in I$.
  
  Let $K = K_0 \cup \bigcup_{J \in \mathcal{J}, m \in \mathbb{N}} \im
  \widehat \varPhi^J_m$ and $\mathcal{K} = \{ K_i \}_{i \in I}$, with
  $K_i = K^0_i \cup \bigcup \widehat
  \varPhi^J_m({\varPhi^J_m}^{-1}(K^0_i))$. By the construction,
  $\langle K, \mathcal{K} \rangle$ is an $n$-core of~$\F{}$.  By
  theorem~\ref{thm:construction of an exact core}, there exists a
  deformation of $\F{}$ to a closed star-finite interior $\N{n}$-cover
  ~$\F{0}$ of~$X$ such that $\langle K, \mathcal{K} \rangle$ is an
  exact $n$-core of~$\F{0}$. Since $Z$ is disjoint from $K$, by
  letting $L = K$ when applying theorem~\ref{thm:construction of an
    exact core}, we may require that $\F{0}$ is equal to $\F{}$ on a
  neighborhood of $Z$. We may also require that the deformation is so
  small, that $\F{0} \prec \st_\G{} \F{}$. Observe that $\langle K_0,
  \mathcal{K}_0 \rangle$ is an exact $n$-core of $\F{0}$ as well.
  
  Order $\mathcal{J} \times \mathbb{N}$ into a sequence $\{ (J_q, m_q)
  \}_{q \in \mathbb{N}}$. For each $q \in \mathbb{N}$ let $C_q = \im
  \widehat \varPhi^{J_q}_{m_q}$ and let $D_q = \widehat
  \varPhi^{J_q}_{m_q}(B^{k+1} \times \{ 1 \})$ be the upper base of
  $C_q$.  By the construction of $\widehat \varPhi^J_m$'s, there
  exists a discrete collection $\{ U_q \}_{q \in \mathbb{N}}$ of open
  subsets of $X$ such that $D_q \subset U_q$ and $U_q \cap (K_0 \cup
  \bigcup_{r \neq q} C_r) = \emptyset$ for each $q \in \mathbb{N}$. By
  the construction, $D_q \subset \st_\F{} G^{J_q}_{m_q}$, hence we may
  require (and we do) that $U_q \subset \st^2_\F{} G^{J_q}_{m_q}$ (we
  use the assumption that $\F{}$ is an interior cover that
  refines~$\G{}$). By proposition~\ref{pro:compact is ZF} and by
  corollary~\ref{cor:approximation within a cover}, there exists an
  approximation of the identity on $U_q$ within $\F{0}$ by a closed
  embedding $h_q \colon U_q \to U_q$ with image disjoint from $D_q$.
  By lemma~\ref{lem:local adjustment} we may require (and we do) that
  $id_{X \setminus U_q } \cup h_q$ is a closed embedding into $X$.
  Let $V_q = U_q \setminus h_q(U_q)$. By the construction, $V_q$ is an
  open neighborhood of $D_q$.
  
  We shall recursively construct a sequence of closed interior
  $\N{n}$-covers $\F{q} = \{ F^q_i \}_{i \in I}$ of $X$, $q \geq 1$,
  that are isomorphic to $\F{0}$ and satisfy the following conditions.

  \begin{enumerate}\addtocounter{enumi}{4}
  \item There exists an exact $n$-core $K_q$ of $\F{q}$ such that
    $K_{q-1} \cup C_q$ is a subcomplex of $K_q$ and $K_q \setminus
    (K_{q-1} \cup C_q) \subset V_q$.
  \item $\F{q}$ is equal to $\F{q-1}$ on $C_q \cup (X \setminus V_q)$.
  \item The inclusion $K_{q-1} \cap F^{q-1}_J \subset K_q \cap F^q_J$
    induces epimorphisms on homotopy groups of dimensions less than or
    equal to $k$.
  \item $\varphi^{J_q}_{m_q}$ is a trivial hole in $\F{q}$.
  \end{enumerate}

  Before we will do the construction, we will show that under these
  conditions the limit $\F{\infty} = \lim_{q \to \infty} \F{q}$ (see
  definition~\ref{def:definition of limit}) is a cover that satisfies
  assertion of the lemma. By the construction, $\F{}$ is isomorphic to
  $\F{0}$ and we required that $\F{0}$ is isomorphic to $\F{q}$ for
  each $q \geq 1$. By (6) and by the assumption that the collection
  $\{ V_q \}$ is discrete, the sequence $\F{q}$ stabilizes locally, so
  if elements of $\F{\infty}$ intersect, then the corresponding
  elements of $\F{q}$ intersect for large enough $q$. By (5),
  $\F{\infty}$ is equal to $\F{0}$ on $K_0$ and by the construction,
  $K_0$ is an exact core of $\F{0}$, hence for each $J \in
  \mathcal{J}$ the intersection $K_0 \cap F^0_J = K_0 \cap F^\infty_J$
  is non-empty. Therefore if elements of $\F{0}$ intersect, then the
  corresponding elements of $\F{\infty}$ intersect. Hence $\F{\infty}$
  is isomorphic to $\F{}$. As we noted, the sequence $\F{q}$
  stabilizes locally, hence by theorem~\ref{thm:open subsets of
    n-spaces}, $\F{\infty}$ is a closed interior $\N{n}$-cover of $X$.
  Let $\V{} = \{ V_q \}_{q \geq 1}$. By the construction, $\F{0} \prec
  \st_\F{} \F{}$ and $\F{\infty} \prec \st_\V{} \F{0}$. We assumed
  that $\V{} \prec \st^2_\F{} \G{}$, hence $\F{\infty} \prec
  \st_{\st^2_\F{} \G{}} \st_\F{} \F{}$. By a direct computation,
  $\F{\infty} \prec \st_{\st^3_\F{} \G{}} \F{}$. What is left to do is
  to show that $\F{\infty}$ is a $(k+1)$-regular cover.

  Let $K_\infty = \bigcup_{q \geq 1} K_q$. Let $f_q \colon X \to X$ be
  a closed embedding defined as the identity on $X \setminus
  \bigcup_{r > q} U_r$ and by the formula $f_q(x) = h_r(x)$ for each
  $x \in U_r$ with $r > q$. Let $\mathcal{H}_q = \{ H^q_i \}_{i \in
    I}$ and let $H^q_i = f_q(F^q_i)$ for each $i \in I$. Observe that
  $\mathcal{H}_q$ is a closed interior $\N{n}$-cover of $X \setminus
  \bigcup_{r > q} V_r$ and $K_q$ is its exact $n$-core. Also,
  $\mathcal{H}_{k}$ is equal to $\mathcal{H}_{k-1}$ on $X \setminus
  U_{k}$ and $\mathcal{H}_{k-1}$ is a shrinking of $\mathcal{H}_{k}$
  on $U_{k} \setminus V_{k}$.  By lemma~\ref{lem:passing to the
    limit}, $K_\infty$ is an exact $n$-core of $\lim_{q \to \infty}
  \mathcal{H}_{k}$.  Observe that $\lim_{q \to \infty} \mathcal{H}_q =
  \F{\infty}$, so $K_\infty$ is also an exact $n$-core of
  $\F{\infty}$. Observe that if $J \in \mathcal{J}$, then $K_\infty
  \cap F^\infty_J = \bigcup_{q \geq 1} K_q \cap F^q_J$. By the
  assumptions, $K_0 \cap F^0_J$ has trivial homotopy groups up to
  dimension~$k$, so by (8), $K_q \cap F^q_J$ has trivial homotopy
  groups up to dimension~$k$. Hence $K_\infty \cap F^\infty_J$ has
  trivial homotopy groups up to dimension~$k$, so $\F{\infty}$ is
  $k$-regular. To see that it is $(k+1)$-regular consider a map
  $\varphi \colon S^k \to F^\infty_J$. Since $K_\infty$ is an exact
  $n$-core of $\F{\infty}$, $\varphi$ is homotopic in $F^\infty_J$
  with a map $\varphi' \colon S^k \to F^\infty_J \cap K_\infty$. Since
  the image of $\varphi'$ is compact, it lies in $K_q \cap F^\infty_J$
  for some $q \geq 1$. It follows from (6) and from the assumption
  that $U_q$ is disjoint from $K_0 \cup \bigcup_{r \neq q} C_q$ that
  $K_q \cap F^\infty_J = K_q \cap F^q_J$. By (7), the inclusion $K_0
  \cap F^0_J \subset K_q \cap F^q_J$ induces an epimorphism on
  $k$-dimensional homotopy groups. Hence $\varphi'$ is homotopic, in
  $K_q \cap F^q_J$, with a finite product of $\varphi^J_m$'s, for some
  $m \in \mathbb{N}$. By (8), $\varphi^J_m$ is contractible in $K_r
  \cap F^r_J \subset K_\infty \cap F^\infty_J$ for $r$ such that $J =
  J_r$ and $m = m_r$. Hence $\F{\infty}$ is $(k+1)$-regular.

  {%
\begin{figure}[ht]
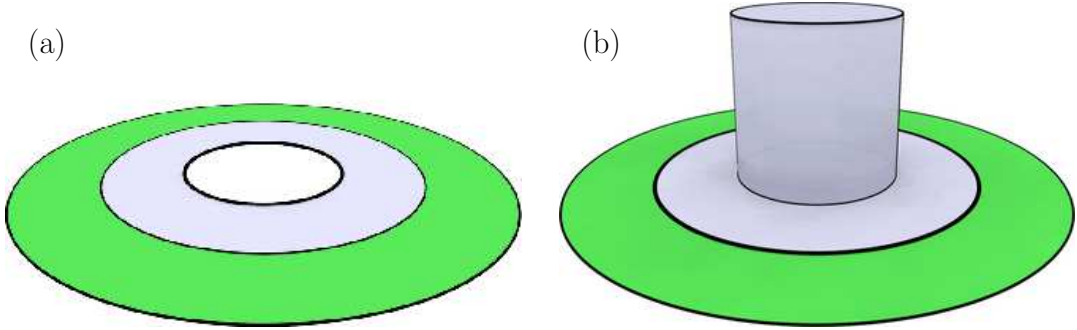

\begin{center}
$\begin{array}{c@{\hspace{5mm}}c}
\resizebox{0.45\textwidth}{!}{\input{psfigures/ocontractible1.pstex_t}} &
\resizebox{0.45\textwidth}{!}{\input{psfigures/ocontractible2.pstex_t}}
\end{array}$
\end{center}
\caption{{After patching of
    $\xi_1$, $\xi_2$ becomes a small hole in the middle
    set.}}
\label{fig:patching of o-contractible hole}
\end{figure}
}

  \emph{The construction of $\F{q}$}. Assume that we already
  constructed $\F{q-1}$. Let $\mathcal{F}^0_{q-1} = \F{q-1}$ and let
  $K_{q-1}^0 = K_{q-1} \cup C_q$.  By the construction, $K_{q-1}^0$ is
  an exact $n$-core of $\mathcal{F}^0_{q-1}$. Let $l =
  l^{J_{q-1}}_{m_{q-1}}$. We shall construct a finite sequence
  $\mathcal{F}^1_{q-1}, \mathcal{F}^2_{q-1}, \ldots,
  \mathcal{F}^l_{q-1}$ of closed interior $\N{n}$-covers of $X$ that
  are isomorphic to $\F{q-1}$ and satisfy the the following conditions
  for each $0 < m \leq l$. Below we let $\mathcal{F}^m_{q-1} = \{
  F^{m, q-1}_i \}_{i \in I}$, $F^{m,q-1}_J = \bigcap_{i \in J}
  F^{m,q-1}_i$ for each $J \in \mathcal{J}$ and $\xi_m =
  \varPhi^{J_q}_{m_q}(xm/l)$.

  \begin{enumerate}\addtocounter{enumi}{8}
  \item There exists an exact $n$-core $K^m_{q-1}$ of
    $\mathcal{F}^m_{q-1}$ such that $K^{m-1}_{q-1}$ is a subcomplex of
    $K^m_{q-1}$ and $K^m_{q-1} \setminus K^{m-1}_{q-1} \subset V_q$.
  \item $\mathcal{F}^m_{q-1}$ is equal to $\mathcal{F}^{m-1}_{q-1}$ on
    $K^{m-1}_{q-1} \cup (X \setminus V_q)$.
  \item The inclusion $K^{m-1}_{q-1} \cap F^{m-1, q-1}_J \subset
    K^m_{q-1} \cap F^{m,q-1}_J$ induces epimorphisms on homotopy
    groups of dimensions less than or equal to $k$.
  \item $\xi_m$ is a trivial hole in $\mathcal{F}^m_{q-1}$.
  \end{enumerate}
  
  Obviously, $\mathcal{F}^l_{q-1}$ is a cover $\F{q}$ that we are
  looking for.  Assume that we already contructed
  $\mathcal{F}^{m-1}_{q-1}$ and $K^{m-1}_{q-1}$.  By the construction,
  the collection $\{ A^j_l \}_{1 \leq j \leq l}$ of subsets of $D_q$
  refines $\mathcal{F}^{m-1}_{q-1}$.  Hence for some $J \in
  \mathcal{J}$, $\xi_m$ is homotopic with $\xi_{m-1}$ in
  $F^{m-1,q-1}_J$.  By (12), $\xi_{m-1}$ is a trivial hole in
  $\mathcal{F}^{m-1}_{q-1}$, hence $\xi_m$ is a small hole in
  $\mathcal{F}^{m-1}_{q-1}$. Let $B_m$ be the ``obvious'' contraction
  of $\xi_m$ (cf. figure~\ref{fig:patching of o-contractible hole}).
  By lemma~\ref{lem:patching of a small hole}, there exists a closed
  interior $\N{n}$-cover $\mathcal{F}^m_{q-1}$ obtained from
  $\mathcal{F}^{m-1}_{q-1}$ by patching of $B_m$ in $V_q$. We let
  $K^m_{q-1}$ to be an exact $n$-core of $\mathcal{F}^m_{q-1}$ that
  satisfies conditions of the definition of a cover obtained by
  patching. A direct verification shows that conditions (9)-(12) are
  satisfied. We are done.
\end{proof}

\begin{remark}
  \label{rem:ok contractible}
  If $\G{} = \st^q \E{}$, then the cover constructed in
  lemma~\ref{lem:ok contractible} refines $\st^{21q + 10} \E{}$. It is
  beause $\F{}$ refines $\G{}$, the constructed cover refines
  $\st_{\st^3_\F{} \G{}} \F{}$ and
  \[
  \st_{\st^3_\F{} \G{}} \F{} \prec \st_{\st^3_\G{} \G{}} \st^3_{\G{}}
  \G{} = \st \st^3 \G{} = \st^{10} \G{},
  \]
  the last equality by lemma~\ref{lem:star of a star}. By the same
  lemma, $\st^{10} \st^q \E{} = \st^{21q + 10} \E{}$.
\end{remark}

As an application we show that every $n$-semiregular cover is
$\circledcirc_0$-contractible, proving theorem~\ref{thm:pump up the
  regularity} for $k = 0$.

\section{Proof of theorem~\lowercase{\ref{thm:pump up the regularity}} for $k = 0$}

Let $\F{} = \{ F_i \}_{i \in I}$ be an $n$-semiregular ($n > 0$)
closed star-finite interior $\N{n}$-cover $n$-contractible in a closed
cover $\mathcal{E}$. By the assumption that $\F{}$ is $n$-semiregular,
there exists an anticanonical map~$\lambda$ of~$\F{}$. By the
definition, every approximation within~$\F{}$ of~$\lambda$ is an
anticanonical map of $\F{}$. Hence by theorem~\ref{thm:approximation
  within a cover}, we may assume that $\lambda$ is a closed embedding.
Fix a non-empty subset $J$ of $I$.  Let $\sigma_J$ be a barycenter of
a simplex of the nerve of $\F{}$ that is spanned by vertices
corresponding to elements of $\F{}$ that are indexed by $J$. We let
$c_0 = \lambda(\sigma_J)$ be the base point of $F_J$. By the
assumption, $\F{}$ is an $\N{n}$-cover and $n > 0$, so $F_J$ is
separable and locally connected, hence there are only countably many
connected components of $F_J$.  Let $c_1, c_2, \ldots$ be a sequence
of points of $F_J$ such that every connected component of $F_J$ has an
element in the sequence. A set of maps $\varphi_k \colon S^0 \to F_J$
that map~$1$ to $c_0$ and $-1$ to $c_k$ generates $\pi_0(F_J)$.  We
will show that for each $k$, $\varphi_k$ is
$\circledcirc_\F{}$-contractible in an element of $\st \E{}$.

By the assumption that $\F{}$ is $n$-contractible in $\E{}$, there
exists a path $\gamma \colon [0,1] \to X$ that connects $c_k$ with
$c_0$ inside an element of~$\E{}$, i.e. such that $\gamma(0) = c_k$,
$\gamma(1) = c_0$ and $\im \gamma$ is a subset of some $E \in \E{}$.
Since $\F{}$ is an interior cover, there exists a positive integer $l$
such that $\{ \gamma([(m-1)/l,m/l]) \}_{0 < m \leq l}$ refines $\F{}$.
For each $0 < m \leq l$, let $j_m$ be an element of $I$ such that
$\gamma([(m-1)/l,m/l]) \subset F_{j_m}$. Without a loss of generality
we may assume that $\gamma$ is constant on small neighborhoods of~$0$
and~$1$, so we may assume that $j_1$ and $j_l$ are elements of $J$.
Consider a function $[(m-1)/l, m/l] \mapsto \lambda(\bst v(F_{j_m}))$
defined on a collection $\{ [(m-1)/l,m/l] \}_{0 < m \leq l}$. It is a
carrier and a constant map $\beta \colon \{ 0, 1 \} \to X$ that maps
$0$ and $1$ to $c_0$ is carried by it. By lemma~\ref{lem:regular
  stars}, the cover of the nerve of $\F{}$ by barycentric stars of its
vertices is regular for the class of metric spaces, hence by the
carrier theorem, $\beta$ can be extended over $[0,1]$ to a map $\beta
\colon [0,1] \to X$ such that $\beta([(m-1)/l,m/l]) \subset
\lambda(\bst v(F_{j_m})) \subset F_{j_m}$ for each $0 < m \leq l$. Let
$\varPhi \colon [-1, 1] \to X$ be a map defined by the formula
\[
  \varPhi(x) = \left\{ \begin{array}{ll}
      \gamma(x + 1) & x \in [-1,0] \\
      \beta(x) & x \in [0,1]
    \end{array} \right..
\]
This map is a $\circledcirc_\F{}$-contraction of $\varphi$ directly
from the definition. The image of $\varPhi$ is contained in $\st_\F{}
\im \gamma$, $\im \gamma$ is contained in an element of $\E{}$ and
$\F{}$ refines $\E{}$, hence the image of $\varPhi$ is contained in an
element of $\st \E{}$.

We proved that $\F{}$ is $\circledcirc_0$-contractible in $\st \E{}$.
By lemma~\ref{lem:ok contractible} and by remark~\ref{rem:ok
  contractible}, there exists a $1$-regular $n$-semiregular closed
interior $\N{n}$-cover of $X$ that is isomorphic to $\F{}$, that
refines $\st^{31} \E{}$ and that is equal to $\F{}$ on some open
neighborhood of $Z$. We are done.

\section{$\Cup$-Contractibility}

\begin{definition*}
  Let $\F{}$ and $\G{}$ be covers of a space $X$. Let $A$ be a
  subspace of $X$. Let $\varphi$ be a map of a pair $(B^k, S^{k-1})$
  into a pair $(X, A)$. We say that $\varphi$ is
  \df{map!uf-contractible@$\Cup_{\F{}}$-contractible}
  {$\Cup_{\F{}}$-contractible to $A$ in $\G{}$} if it admits a
  homotopy $\varPhi \colon (B^k, S^{k-1}) \times [0, 1] \to (X, A)$
  with a constant map and such that the image of $\varPhi$ lies in an
  element of $\G{}$ and for some integer $l$, the collection $\{ B^k
  \times [(m-1)/l, m/l] \}_{1 \leq m \leq l}$ refines
  $\varPhi^{-1}(\F{})$.  We shall call $\varPhi$ a
  \df{uF-contraction@$\Cup_{\F{}}$-contraction}{$\Cup_{\F{}}$-contraction
    of $\varphi$ to $A$ in $\G{}$}.
\end{definition*}

\begin{definition*}
  Let $0 < k < n$ and let $A$ be a subspace of a space $X$. Let $\F{}
  = \{ F_i \}_{i \in I}$ be a cover of $X$. Assume that each non-empty
  $F_J$ has a basepoint $x_J$ in the intersection $F_J \cap A$. We say
  that $\F{}$ is
  \df{cover!uk-contractible@$\Cup_k$-contractible}{$\Cup_k$-contractible
    to $A$ in a cover $\G{}$} if for each non-empty subset $J$ of $I$
  such that $F_J$ is non-empty there exists a set of maps $\{
  \varphi^J_m \colon (B^k, S^{k-1}) \to (F_J, F_J \cap A) \}_{m \geq
    1}$ that are $\Cup_\F{}$-contractible in $\G{}$ and such that
  $i_*^{-1}(\{ \varphi^J_m \}_{m \geq 1})$ generates $\pi_k(F_J)$,
  where $i_* \colon \pi_k(F_J) \to \pi_k(F_J, F_J \cap A)$ is a
  function induced by the inclusion $(F_J, x_J) \subset (F_J, F_J \cap
  A)$.
\end{definition*}

\begin{lemma}
  \label{lem:u-contraction and zero homomorphisms}
  Let $0 < k < n$ and let $A$ be a subspace of a space $X$. Let $\F{}
  = \{ F_i \}_{i \in I}$ be a cover of $X$. As usual, let $F_J =
  \bigcap_{i \in J} F_i$ for each non-empty subset $J$ of $I$. Assume
  that each non-empty $F_J$ has a base-point $x_J$ in the intersection
  $F_J \cap A$.  Let $\G{}$ be a closed cover of $X$. Assume that for
  each non-empty subset $J$ of $I$ the $(k-1)$th homotopy group of $A
  \cap F_J$ is trivial and that the inclusion of $A \cap F_J$ into
  $F_J$ induces a trivial (zero) homomorphism on $k$th homotopy groups.
  If $\F{}$ is $\Cup_k$-contractible to $A$ in $\G{}$, then $\F{}$ is
  $\circledcirc_k$-contractible in $\st \G{}$.
\end{lemma}
\begin{proof}
  Fix $m \geq 1$ and fix a non-empty $J \subset I$ such that $F_J \neq
  \emptyset$. Let $\{ \varphi^J_m \}_{m \geq 1}$ be a set of
  generators of $\pi_k(F_J)$ such that compositions $\varphi^J_m \circ
  \xi$ are $\Cup_\F{}$-contractible to $A$ in $\G{}$. Pick $\varphi =
  \varphi^J_m$. We will show that $\varphi$ is
  $\circledcirc_\F{}$-contractible in $\st \G{}$, thus proving the
  assertion directly from the definition of
  $\circledcirc_k$-contractibility.  Let $B^k_+$ be the north
  hemisphere of $S^k$, let $B^k_-$ be the south hemisphere of $S^k$
  and let $S^{k-1} = B^k_+ \cap B^k_-$ be the equator of $S^k$.
  Without a loss of generality, we may let $B^k_+$ denote the domain
  of $\varphi$. Let $\varPhi$ be a $\Cup_{\F{}}$-contraction to $A$ of
  $\varphi \circ \xi$ in some element $G$ of $\G{}$. It is a homotopy
  of $\varphi$ with a constant map.  For $0 \leq t \leq 1$ we let
  $\varPhi_t \colon B^k_+ \to X$ denote $\varPhi$ restricted to a time
  $t$. By the definition of $\Cup_{\F{}}$-contraction to $A$, for each
  $t$, $\varPhi_t(S^{k-1}) \subset A$ and there exists a positive
  integer $l$ and a sequence of elements $j_1, j_2, \ldots, j_l \in I$
  such that
  \begin{enumerate}
  \item[(*)] $\varPhi_t \subset F_{j_q}$ for each $0 < q \leq l$ and
    each $t$ such that $(q-1)/l \leq t \leq q/l$.
  \end{enumerate}
  We will extend $\varPhi$ over $B^k_- \times [0,1]$ to a
  $\circledcirc_\F{}$-contraction of $\varphi$.

  By the definition, $\varPhi_0$ restricted to $S^{k-1} \times \{0 \}$
  is constant. We extend $\varPhi_0$ onto $B^k_- \times \{ 0 \}$ to be
  constant as well, with value determined by the value of $\varPhi$ on
  $S^{k-1} \times \{ 0 \}$.  Again by the definition, $\varPhi_1$
  restricted to $S^{k-1} \times \{1\}$ is constant and we extend
  $\varPhi_1$ onto $B^k_- \times \{ 1 \}$ to be a constant map.  For
  each integer $0 < q < l$, $\varPhi_{q/l}$ maps $S^{k-1} \times \{
  q/l \}$ into $F_{j_{q-1}} \cap F_{j_q} \cap A$. By the assumptions,
  this restriction is contractible in $F_{j_{q-1}} \cap F_{j_q} \cap
  A$. We extend $\varPhi$ over $B^k_- \times \{ q/l \}$ to be any such
  contraction.  For each $0 < q \leq l$, the (extended) map
  $\varPhi_{|S^{k-1} \times [(q-1)/l,q/l] \cup B^k_- \times \{
    (q-1)/l, q/l \}}$ is a map from a $k$-dimensional sphere into $A
  \cap F_{j_q}$. By the assumptions, every such map extends over
  $B^k_- \times [ (q-1)/l, q/l ]$ to a map into $F_{j_q}$. We define
  an extension of $\varPhi$ onto $B^k_- \times [(q-1)/l,q/l]$ to be
  any such extension.

  From the construction, the extended $\varPhi$ is a
  $\circledcirc_\F{}$-contraction of $\varPhi_0$ and $\varPhi_0$ is
  homotopic with $\varphi$ in $F_J$. Its image lies in $\st_\F{}
  \G{}$. Since $\F{}$ refines $\G{}$, it is a
  $\circledcirc_\F{}$-contraction in $\st \G{}$.

\end{proof}

\begin{example}
  Let $\F{} = \{ F_i \}_{i \in I}$ be a cover of a space $X$. Assume
  that $\F{}$ admits an anticanonical map $\lambda$ that is a closed
  embedding and such that $\lambda(\bst v(F_i)) = F_i \cap \im
  \lambda$ for each $i \in I$.  By lemma~\ref{lem:u-contraction and
    zero homomorphisms}, if $\F{}$ is $\Cup_k$-contractible to $\im
  \lambda$ in $\G{}$, then $\F{}$ is $\circledcirc_k$-contractible in
  $\st \G{}$.
\end{example}

\section{Patching of large holes} 

\begin{definition}\label{def:patch}
  Let $\F{} = \{ F_i \}_{i \in I}$ be a cover of a space $X$ and let
  $F_J = \bigcap_{i \in J} F_i$ for each non-empty subset $J$ of $I$.
  Let $\tau$ be the set of simplices of a triangulated $k$-dimensional
  sphere $S^k$. Let $\varphi \colon S^k \to F_J$ be a $k$-hole in $F_J$.
  Let
  \[
    C = \left(S^k \times [-1, 1] \times \{ 0 \}) \cup ( (S^k)^{(k-2)}
    \times [-1, 1] \times [0, 1]\right),
  \]
  where $(S^k)^{(k-2)}$ denotes the $(k-2)$-dimensional skeleton of
  $S^k$ (empty if $k = 1$). We say that a map $\varPhi \colon C \to X$
  is a \df{patch}{patch for $\varphi$} if the following conditions are
  satisfied.

  \begin{enumerate}
    \item $\varPhi$ restricted to $S^k = S^k \times \{ -1 \} \times \{
    0 \}$ is equal to $\varphi$,
    \item $\varPhi$ restricted to $S^k \times \{ 1 \} \times \{ 0 \}$
    is a trivial hole in $F_J$,
    \item there is an integer $l$ such that the collection 
      \[\{ \varPhi((\delta \times B_{m, l}) \cap \dom \varPhi) \colon \delta \in \tau, 1
      \leq m \leq l \}\] refines $\F{}$, where $B_{m, l} = \{ x \in
      [-1,1] \times [0,1] \colon |x|_\infty \in [\frac{m-1}{l},
      \frac{m}{l}] \}$ (see figure~\ref{fig:bml}).

      {%
\begin{figure}[ht]
\begin{center}
\resizebox{0.60\textwidth}{!}{\input{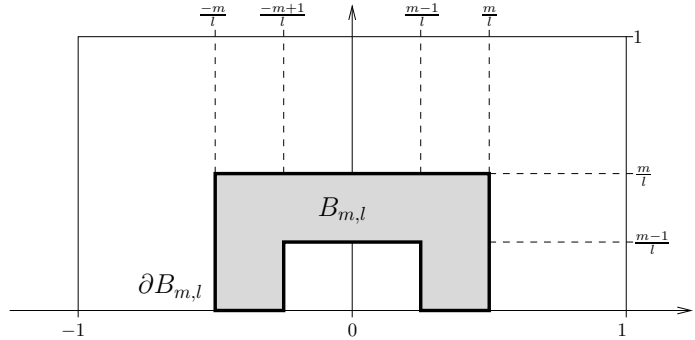}}
\end{center}
\caption{{$B_{m,l} = \{ x \in [-1,1] \times [0,1] \colon
        |x|_\infty \in [\frac{m-1}{l}, \frac{m}{l}] \}$.}}
\label{fig:bml}
\end{figure}
}
  \end{enumerate}
\end{definition}

\begin{remark}\label{rem:interior cover of a patch}
  Let $\F{} = \{ F_i \}_{i \in I}$. If $\varPhi$ is a patch for a
  $k$-hole $\varphi$ if $\F{}$, then by the definition there exists a
  triangulation $\tau$ of $S^k$ and an integer $l$ such that the
  collection
  \[
  \left\{ \bigcup \left\{ \delta \times B_{m,l} \colon \delta \in
      \tau, 1 \leq m \leq l, \varPhi(\delta \times B_{m,l}) \subset
      F_i \right\} \right\}_{i \in I}
  \]
  is a cover of the domain of~$\varPhi$ by subcomplexes.  We can
  always replace $\varPhi$ by another patch such that this collection
  (possibly for a finer triangulation~$\tau$ and a larger integer~$l$)
  is an interior cover of the domain of $\varPhi$. The new patch can
  be obtained by a reparametrization of $\varPhi$. In particular, we
  may assume that its image is equal to the image of $\varPhi$.
\end{remark}

We'll prove below two lemmas. The first of them states that every hole
in $\F{}$ admits a patch the size of which is controlled, while the
next lemma enables one to use this patch to patch the hole. Some new
holes may be added in the process, but they will be handled by
lemma~\ref{lem:u-contraction and zero homomorphisms}.

\begin{lemma}
  \label{lem:the existence of patches}
  If $\F{}$ is a star-finite $n$-semiregular interior cover that is
  $n$-contractible in a cover $\E{}$ and $\varphi$ is a $k$-hole in
  $\F{}$ for some $k < n$, then $\varphi$ has a patch whose image is
  contained in an element of $\st \E{}$, for some sufficiently fine
  triangulation of $S^k$.
\end{lemma}
\begin{proof}
  We denote the underlying space by $X$. Let $\F{} = \{ F_i \}_{i \in
    I}$ and let $\varphi \colon S^k \to F_J$ be a $k$-hole in $F_J$
  for some non-empty $J \subset I$. Identify $S^k$ with $S^k \times \{
  -1 \} \times \{ 0 \}$ and let $\varPhi \colon S^k \times \{ -1 \}
  \times \{ 0 \} \to X$ to be equal to $\varphi$. We'll extend
  $\varPhi$ consecutively onto $S^k \times [-1, 0] \times \{ 0 \}$,
  $S^k \times [0,1] \times \{ 0 \}$ and $(S^k)^{(k-2)} \times [-1, 1]
  \times [0, 1]$ to obtain a patch for $\varphi$. We shall fix a
  triangulation of $S^k$ after the first extension. In the first step
  we will use the assumption that $\F{}$ is $n$-contractible in
  $\E{}$; in the second step that it is $n$-semiregular; in the third
  step that it is $k$-regular.

  First step.  Let $C_1 = S^k \times [-1, 0] \times \{ 0 \}$. Let
  $\lambda$ be an anticanonical map of $\F{}$. For each element $F$ of
  $\F{}$ we let $v(F)$ be the vertex in the nerve of $\F{}$ that
  corresponds to $F$. Let~$\sigma_J$ be the barycenter of a simplex of
  $N(\F{})$ spanned by the vertices $v(F_j)$, $j \in J$. We define
  $\varPhi$ on $S^k \times \{ 0 \} \times \{ 0 \}$ to be a constant
  map with value $\lambda(\sigma_J)$.  Since $\F{}$ is
  $n$-contractible in $\E{}$ and $\varPhi$ (so far defined) has image
  contained in $F_J$ and a domain whose dimension is less than $n$,
  $\varPhi$ extends over $C_1$ to a map with image contained in some
  element of~$\E{}$. We define $\varPhi$ on $C_1$ to be any such
  extension.
  
  Second step. Let $C_2 = S^k \times [0, 1] \times \{ 0 \}$.  We first
  fix a triangulation $\tau$ of $S^k$ and an integer $l$ such that for
  every $1 \leq m \leq l$ and every $\delta \in \tau$ the set
  $\varPhi(\delta \times [-m/l, (-m+1)/l]) \times \{ 0 \}$ is a subset
  of $F_{i_{\delta, m}}$ for some $i_{\delta, m} \in I$. Such a pair
  exists because $\F{}$ is an interior cover and $C_1$ is compact.
  Without a loss of generality we may assume that $i_{\delta, 1} \in
  J$ and $i_{\delta, l} \in J$ for each $\delta \in \tau$. A function
  that assigns $\bst v(F_{i_{\delta, m}})$ to a subset $\delta \times
  [(m-1)/l, m/l] \times \{ 0 \}$ of $C_2$ is a carrier. A constant map
  on $S^k \times \{ 0, 1 \} \times \{ 0 \}$ with value $\sigma_J$ is
  carried by it.  By lemma~\ref{lem:regular stars}, the cover of the
  nerve of $\F{}$ by barycentric stars of vertices is $n$-regular.
  Therefore by the carrier theorem, the constant map extends over
  $C_2$ to a map that maps $\delta \times [(m-1)/l, m/l]) \times \{ 0
  \}$ into $\bst v(F_{i_{\delta, m}})$.  We define $\varPhi$ on $C_2$
  to be a composition with $\lambda$ of any such extension.
  
  Third step. Let $C_3 = (S^k)^{(k-2)} \times [-1, 1] \times [0, 1]$.
  By the construction, $\varPhi((C_1 \cup C_2) \cap (\delta \times
  B_{m,l})) \subset F_{i_{\delta, m}}$ for each $\delta \in \tau$ and
  $1 \leq m \leq l$. Therefore the map $\delta \times B_{m,l} \mapsto
  F_{i_{\delta, m}}$ is a carrier.  By assumptions $\F{}$ is a
  $k$-regular cover and by the definition $C_3$ is $k$-dimensional.
  Hence by the carrier theorem, $\varPhi$ extends over $C_3$ to a map
  carried by this carrier. We define $\varPhi$ on $C_3$ to be any such
  extension.
  
  By the construction $\varPhi$ is a patch for $\varphi$. Observe that
  $\varPhi(C_1)$ is contained in some element of $\E{}$. Therefore by
  condition (3) of the definition of a patch, since $\F{}$ refines $\E{}$, the
  image of $\varPhi$ must be contained in an element of $\st \E{}$.
\end{proof}

\begin{lemma}
  \label{lem:patching of a large hole}
  Let $0 < k < n$. Let $\F{} = \{ F_i \}_{i \in I}$ be a closed
  countable $k$-regular star-finite interior $\N{n}$-cover of a space
  $X$. Let $J$ be a nonempty subset of $I$ and let~$\varphi$ be a
  closed embedding of $S^k$ into $F_J = \bigcap_{i \in I} F_i$. Let
  $\varPhi$ be a closed embedding and a patch for $\varphi$. Let $L$
  be an exact $n$-core of $\F{}$ such that $\im \varPhi$ is a
  subcomplex of $L$.  Let $A$ be a subcomplex of $L$ such that $\im
  \varPhi$ is a subcomplex of $A$. Let $\mathcal{H}$ be a cover of~$L$
  such that $\F{}/L$ is $\Cup_k$-contractible to $A$ in $\mathcal{H}$.
  For every neighborhood $U$ of $\im \varPhi$ there exists a cover
  $\G{}$ of~$X$ and an exact $n$-core $K$ of $\G{}$ such that the
  following conditions are satisfied.
  \begin{enumerate}
  \item $\G{}$ is a $k$-regular closed interior $\N{n}$-cover of $X$
    that is isomorphic to $\F{}$ and is equal to $\F{}$ on $L \cup (X
    \setminus U)$,
  \item $L$ is a subcomplex of $K$ and $K \setminus L \subset U$,
  \item $\G{} / K$ is $\Cup_k$-contractible to $A$ in $\mathcal{H} \cup
    \{ U \cap K \}$,
  \item the element $\varphi$ of $\pi_k(L \cap F_J) = \pi_k(L \cap
    G_J)$ becomes $0$ in $\pi_k(K \cap G_J)$.
  \end{enumerate}
\end{lemma}

In the next section we apply lemma~\ref{lem:patching of a large hole}
recursively as to kill all generators of $\pi_k(L \cap F_J)$, $J
\subset I$, to be able to use lemma~\ref{lem:u-contraction and zero
  homomorphisms} and reduce the proof of theorem~\ref{thm:pump up the
  regularity} to the previously solved case of
$\circledcirc_k$-contractible covers.

\begin{proof}
  Let $\tau$ be a set of simplices of $S^k$ and let $l$ be a positive
  integer such that conditions of definition~\ref{def:patch} are
  satisfied for $\varPhi$. Order $(k-1)$-dimensional simplices from
  $\tau$ into a sequence $\delta_0, \delta_1, \ldots, \delta_q$. Let
  $C_{0, 0}$ denote the domain of $\varPhi$. For each $0 \leq p \leq
  q$ and for each $1 \leq m \leq l$, let $C_{p, m}$ be a simplicial
  complex defined recursively by the identity $C_{p, m} = C_{p, m-1}
  \cup (\delta_p \times B_{m,l})$ with $C_{p,0}$ identified with
  $C_{p-1,l}$ for $p > 0$. Note that 
  \[ 
  C_{q, l} = (S^k \times [-1, 1] \times \{ 0 \}) \cup \left(
    (S^k)^{(k-1)} \times [-1, 1] \times [0, 1] \right)
  \]
  and compare $C_{q,l}$ with complex $C$ from
  definition~\ref{def:patch}. Let $D_{p, m} = (\delta_p \times
  B_{m,l}) \cap C_{p,m-1}$ for each $0 \leq p \leq q$ and each $0 < m
  \leq l$.

  {%
\begin{figure}[ht]
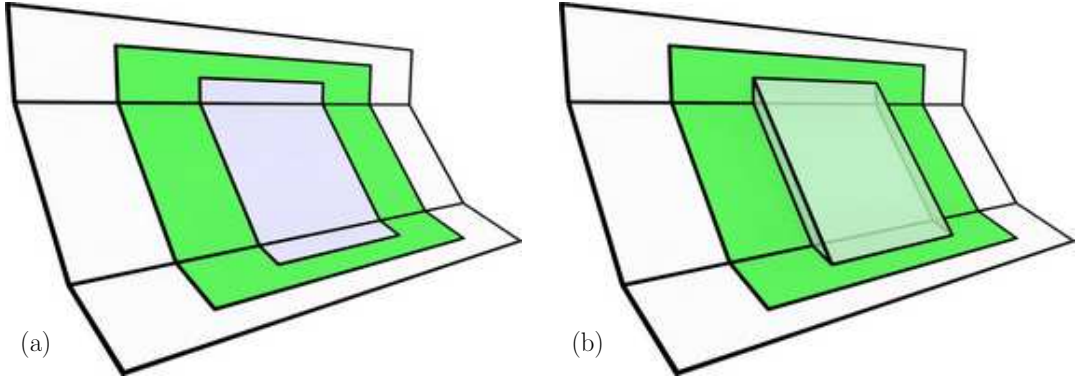

\begin{center}
$\begin{array}{c@{\hspace{5mm}}c}
\resizebox{0.45\textwidth}{!}{\input{psfigures/patching1.pstex_t}} &
\resizebox{0.45\textwidth}{!}{\input{psfigures/patching2.pstex_t}}
\end{array}$
\end{center}
\caption{{Patching of
    $(k-1)$-holes.}}
\label{fig:patching}
\end{figure}
}

  Since $\varPhi$ is an embedding, we may identify $C_{0, 0}$ with the
  image of $\varPhi$ and treat it as a subcomplex of $L$. Below we
  shall recursively embed $C_{p,m}$'s into $X$ and we will use same
  symbols to denote images of these embeddings. Let $\F{0,0} = \F{}$
  and $L_{0,0} = L$. We shall recursively construct a sequence of
  covers $\F{p,m} = \{ F^{p,m}_i \}_{i \in I}$ and their exact
  $n$-cores $L_{p,m}$ in such a way that the following conditions are
  satisfied ($L_{p,0}$ is identified with $L_{p-1, l}$ for each $p >
  0$).

  \begin{enumerate}[(i)]
  \item $\F{p,m}$ is a $k$-regular closed interior $\N{n}$-cover of
    $X$ that is isomorphic to $\F{}$ and is equal to $\F{}$ on $L \cup
    (X \setminus U)$,
  \item $C_{p,m}$ is embedded into $X$ in such a way that $L_{p,m} =
    L_{p,m-1} \cup C_{p,m}$ and $C_{p,m}$ is a subset of $U$,
  \item $\F{p,m} / L_{p,m}$ is $\Cup_k$-contractible to $A$ in
    $\mathcal{H} \cup \{ U \cap L_{p,m} \}$,
  \item if $D_{p,m} \subset F^{p,m-1}_i$, then $\delta_p \times
    B_{m,l} \subset F^{p,m}_i$; if $\partial D_{p,m} \subset
    F^{p,m-1}_i$, then $B_{p,m+1} \subset F^{p,m}_i$.
  \end{enumerate}

  Identify $\F{p,0}$ with $\F{p-1,l}$ for each $p > 0$ and assume that
  we already contructed $\F{p,m-1}$. By conditions (i) and (iv),
  conditions of lemma~\ref{lem:patching of a small hole} are satisfied
  with $\F{} = \F{p,m-1}$, $K = L_{p,m-1}$ and $B = B_{p,m}$. Hence
  there exists a closed interior $\N{n}$-cover $\F{p, m}$ with an
  exact $n$-core $L_{p, m}$ that is obtained from $\F{p, m-1}$ by
  patching of $B_{p, m}$ in $U$. We may identify $C_{p,m} \setminus
  C_{p,m-1}$ with $L_{p,m} \setminus L_{p,m-1}$, as these sets are
  homeomorphic by the definitions. Conditions (i), (ii), (iii) and
  (iv) are satisfied directly from the definitions.

  Order $k$-dimensional simplices of $\tau$ into a sequence $\Delta_0,
  \Delta_1, \ldots, \Delta_Q$. Let $E_{0,0} = C_{q, l}$. For each $0
  \leq p \leq Q$ and for each $1 \leq m \leq l$, let $E_{p,m}$ be a
  simplicial complex defined recursively by the identity $E_{p,m} =
  E_{p,m-1} \cup (\Delta_p \times B_{m,l})$ if $k < n - 1$ or by the
  identity $E_{p,m} = E_{p,m-1} \cup (\Delta_p \times \partial
  B_{m,l})$ if $k = n - 1$. In either case we identify $E_{p,0}$ with
  $E_{p-1,l}$ for each $p > 0$. Let $F_{p,m} = (\Delta_p \times
    B_{m,l}) \cap E_{p,m-1}$ for each $0 \leq p \leq Q$ and each $0 <
    m \leq l$.
 
  {%
\begin{figure}[ht]
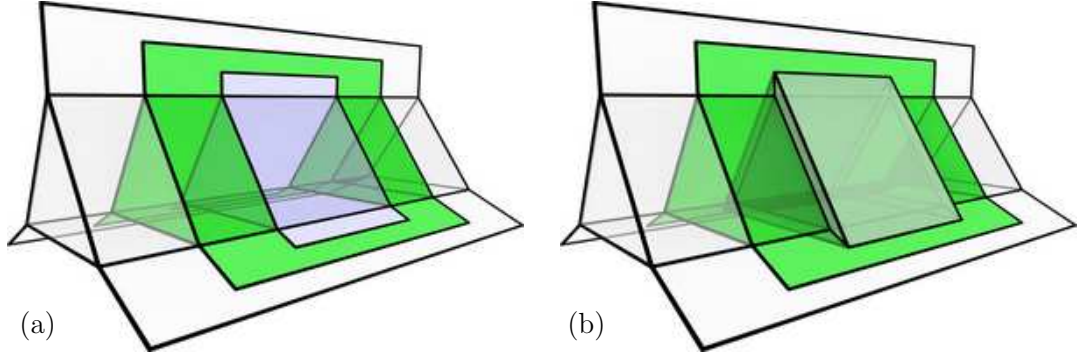

\begin{center}
$\begin{array}{c@{\hspace{5mm}}c}
\resizebox{0.45\textwidth}{!}{\input{psfigures/patching3.pstex_t}} &
\resizebox{0.45\textwidth}{!}{\input{psfigures/patching4.pstex_t}}
\end{array}$
\end{center}
\caption{{Patching of
    $k$-holes.}}
\label{fig:patching of k}
\end{figure}
}

  By the definition, $E_{0, 0}$ is a subset of $X$. Below we shall
  recursively embed $E_{p,m}$'s into $X$ and we will use same symbols
  to denote images of these embeddings. Let $\G{0,0} = \F{q,l}$ and
  let $K_{0,0} = L_{q,l}$. We shall recursively construct a sequence
  of covers $\G{p,m} = \{ G^{p,m}_i \}_{i \in I}$ and their exact
  $n$-cores $K_{p,m}$ in such a way that the following conditions are
  satisfied ($K_{p,0}$ is identified with $K_{p-1,l}$ for each $p >
  0$).

  \begin{enumerate}[(I)]
  \item $\G{p,m}$ is a $k$-regular closed interior $\N{n}$-cover of
    $X$ that is isomorphic to $\F{}$ and is equal to $\F{}$ on $L \cup
    (X \setminus U)$,
  \item $E_{p,m}$ is embedded into $X$ in such a way that $K_{p,m} =
    K_{p,m-1} \cup E_{p,m}$ and $E_{p,m}$ is a subset of $U$.
  \item $\G{p,m} / K_{p,m}$ is $\Cup_k$-contractible to $A$ in
    $\mathcal{H} \cup \{ U \cap K_{p,m} \}$,
  \item if $F_{p,m} \subset G^{p,m-1}_i$, then $\Delta_p \times
    B_{m,l} \subset G^{p,m}_i$ if $k < n - 1$ and $\Delta_p \times
    \partial E_{p,m} \subset G^{p,m}_i$ if $k = n - 1$; if $\partial
    F_{p,m} \subset G^{p,m-1}_i$, then $F_{p,m+1} \subset G^{p,m}_i$.
  \end{enumerate}

  Again, we obtain $\G{p,m}$ from $\G{p,m-1}$ by patching of $E_{p,m}$
  in $U$.  All four conditions are verified easily.

  Let $\G{} = \G{Q, l}$ and $K = K_{Q,l}$. Observe that if $\im
  \varphi \subset F_i$ for $i$ in $I$, then by (iv) and (IV), for
  every simplex $\delta \in \tau$, $\delta \times \{ x \in [-1,1]
  \times [0,1] \colon |x|_\infty = 1 \} \subset G^{Q, l}_i$. Hence
  $\varphi$ is contractible in $G^{Q,l}_i$. By (I), (II) and (III),
  $\G{}$ satisfies conditions (1), (2) and (3) of the lemma. We are
  done. 
\end{proof}

\section{Proof of theorem~\lowercase{\ref{thm:pump up the regularity}} for $k > 0$}

Let $\F{} = \{ F_i \}_{i \in I}$ be a $k$-regular $n$-semiregular ($0
< k < n$) closed countable star-finite interior $\N{n}$-cover of a
space $X$. Assume that $\F{}$ is $n$-contractible in a closed cover
$\E{}$ of $X$. Let $Z$ be a $Z(\F{})$-set. We shall construct a closed
$(k+1)$-regular $n$-semiregular closed countable star-finite interior
$\N{n}$-cover of $X$ that refines $\st^{472} \E{}$, that is isomorphic to
$\F{}$ and that is equal to $\F{}$ on a neighborhood of $Z$.

Let $\mathcal{J}$ be a set of those non-empty subsets $J$ of $I$, for
which the intersection $F_J = \bigcap_{i \in I} F_i$ is non-empty.
The set $\mathcal{J}$ is countable, because $\F{}$ is countable and
locally finite. By proposition~\ref{pro:countable homotopy groups},
the set of homotopy classes of maps from $S^k$ into~$F_J$ is countable
for each $J$ in $\mathcal{J}$. Hence there exists a sequence $J_m$ of
elements of $\mathcal{J}$ and a sequence of maps $\varphi_m \colon S^k
\to F_{J_m}$ such that for each $J \in \mathcal{J}$, the set $\{
\varphi_m \}_{m \geq 1, J_m = J}$ contains elements from all homotopy
classes of maps from $S^k$ into~$F_J$. By lemma~\ref{lem:the existence
  of patches}, for each $\varphi_m$ there exists a patch $\varPhi_m$
whose image is contained in an element of $\st \E{}$. We assume that
the assertion of remark~\ref{rem:interior cover of a patch} is true
for $\varPhi_m$, i.e. that the collection defined in its statement is
an interior cover of $\varPhi_m$. Let~$\mathcal{O}$ be an open cover
of $X$ obtained via lemma~\ref{lem:close approximation} applied to
$\F{}$.  Let $\varPhi$ be a disjoint union of $\varPhi_m$'s. By
theorem~\ref{thm:complete complexes}, the domain of $\varPhi$ is
Polish. By theorem~\ref{thm:approximation within a cover}, there
exists a $\mathcal{O}$-approximation $\Psi$ of $\varPhi$ within~$\F{}$
by a closed embedding. By proposition~\ref{pro:compact is ZF} and by
corollary~\ref{cor:approximation within a cover}, we may assume that
the image of $\Psi$ is disjoint from $Z$. By a slight abuse of
notation, we may let $\Psi_m$ denote the restriction of $\Psi$ to the
domain of $\varPhi_m$ and let $\psi_m$ denote the restriction of
$\Psi$ to the domain of $\varphi_m$.  By the construction, the
collection of images of $\Psi_m$'s is discrete in $X$, $\psi_m$ maps
$S^k$ into $F_{J_m}$, for each $J$ in $\mathcal{J}$ the set $\{ \psi_m
\}_{m \geq 1, J_m = J}$ contains elements from all homotopy classes of
maps from $S^k$ into $F_J$, $\Psi_m$ is a patch for $\psi_m$ and the
image of $\Psi_m$ is contained in $\st_\F{} \E{}$.

Let $\tau_m$ be a triangulation of $S^k$ and let $l_m$ be an integer
such that $\varPsi_m$ satisfies conditions of the definition of a
patch for $\tau = \tau_m$ and $l = l_m$. Let $M = \im \Psi$ and let
\[
\mathcal{M} = \left\{ M_i = \bigcup \left\{ \varPhi(\delta \times
    B_{m,l_m}) \colon \delta \in \tau_m, 1 \leq m \leq l_m,
    \varPhi(\delta \times B_{m,l_m}) \subset F_i \right\} \right\}_{i
  \in I}.
\]
By the assumption we made on $\varPhi_m$'s, $\mathcal{M}$ is an
interior cover of $M$. By the definition, it refines $\F{}$. Hence by
theorem~\ref{thm:embedded core}, there exists an $n$-core $\langle
L_0, \mathcal{L}_0 \rangle$ of $\F{}$ such that~$M$ is a subcomplex of
$L_0$ and such that $\mathcal{L}_0$ is equal to $\mathcal{M}$ on $M$.
Let $ \widetilde L$ be a disjoint union $L_0 \sqcup \bigsqcup_{m \geq
  1} (\im \varPhi_m \times [0, 1])$, divided by an equivalence
relation that for each $m \geq 1$ identifies $\im \varPhi_m \subset
L_0$ with $\im \varPhi_m \times \{ 0 \} \subset \im \varPhi_m \times
[0, 1]$.  Let $\widetilde p$ be a restriction to the $n$-dimensional
skeleton of $\widetilde L$ of a natural projection from $\widetilde L$
onto $L_0$.  By theorem~\ref{thm:approximation within a cover}, there
exists an approximation of $\widetilde p$ within~$\F{}$ $\rel L_0$ by
a closed embedding $p$. By proposition~\ref{pro:compact is ZF} and by
corollary~\ref{cor:approximation within a cover}, we may assume that
the image of $p$ is disjoint from $Z$.  We let $\mathcal{L} = \{
p(\widetilde p^{-1}(L_0^i)) \}_{i \in I}$. We let $L = \im p$. We let
$L_m = L_{m-1} \cup p(\widetilde p^{-1}(\im \Psi_m))$.  By the
construction, for each $m \geq 0$, $\langle L_m, \mathcal{L} / L_m
\rangle$ is an $n$-core of $\F{}$.  Let $\widetilde \Xi_m$ be the
unique map onto $\im \varPhi_m \times \{ 1 \}$ such that $p \circ
\Xi_m = \varPhi_m$ and let $\widetilde \xi_m$ be the unique map onto
$\im\varphi_m \times \{ 1 \}$ such that $p \circ \xi_m = \varphi_m$.
Let $\Xi_m = p \circ \widetilde \Xi_m$ and let $\xi_m = p \circ
\widetilde \xi_m$. By the construction, $\Xi_m$ is an approximation of
$\Psi_m$ within $\F{}$.  Hence $\Xi_m$ is a patch for $\xi_m$. The
image of $\Psi_m$ lies in an element of $\st_\F{} \E{}$ and so does
the image of $\Xi_m$. By the construction, $\xi_m$ is homotopic with
$\psi_m$ in $F_{J_m}$.  Let $\U{} = \{ U_m \}_{m \geq 1}$ be a
discrete collection of open subsets of $X$ such that $\im \Xi_m
\subset U_m$, $U_m \cap Z = \emptyset$, $L_{m-1} \cap U_m = \emptyset$
and $(L_{m+1} \setminus L_m) \cap U_m = \emptyset$ for each $m \geq
1$. Since the collection $\{ \im \Xi_m \}_{m \geq 1}$ refines
$\st_\F{} \E{}$, we may assume that $\U{}$ refines $\st^2_\F{} \E{}$,
because $\F{}$ is an interior cover (which implies that for each set
$A$, $\st_\F{} A$ is a neighborhood of $A$).

{%
\begin{figure}[ht]
\begin{center}
\resizebox{0.60\textwidth}{!}{\input{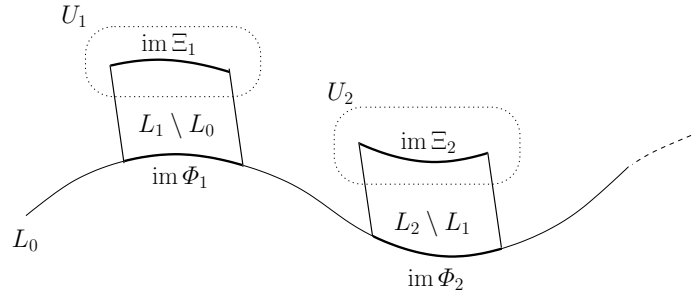}}
\end{center}
\caption{{Some of the sets constructed in the proof of
  theorem~\ref{thm:pump up the regularity}.}}
\label{fig:patches in a core}
\end{figure}
}

By theorem~\ref{thm:construction of an exact core}, there exists a
deformation of $\F{}$ to a closed interior $\N{n}$-cover $\F{-1} = \{
F^{-1}_i \}_{i \in I}$ of $X$ with an exact $n$-core $\langle L,
\mathcal{L} \rangle$. We may assume that $\F{-1}$ is equal to $\F{}$
on a neighborhood of $Z$ and that $\F{-1}$ refines $\st \F{}$. By the
construction, each $\Xi_m$ is a patch for $\xi_m$ and for each $J \in
\mathcal{J}$, the set $\{ \xi_m \}_{m \geq 1, J_m = J}$ contains
elements from all homotopy classes of maps from $S^k$ into~$F_J$.

By lemma~\ref{lem:local adjustment}, for each $m \geq 1$ there exists
an open cover $\U{m}$ of $U_m$ such that every map from $U_m$ into
$U_m$ that is $\U{m}$-close to the identity on $U_m$ extends
continuously by the identity on $X \setminus U_m$. By
theorem~\ref{thm:open subsets of n-spaces}, each $U_m$ is an
$\N{n}$-space and by remark~\ref{rem:closed locally compact is ZF},
$\im \Xi_m$ is a $Z(\F{} / U_m)$-set in $U_m$. By
corollary~\ref{cor:approximation within a cover}, there exists a
$\U{m}$-approximation $h_m \colon U_m \to U_m$ within $\F{-1}$ of the
identity on $U_m$ by a closed embedding with image disjoint from $\im
\Xi_m$.  We let $V_m = U_m \setminus \im h_m$.  Note that $V_m$ is an
open neighborhood of $\im \Xi_m$. For each $m \geq 0$, let $H_m \colon
X \to X \setminus \bigcup_{l > m} V_l$ be a homeomorphism defined by
the formula
\[
  H_m(x) = \left\{ 
    \begin{array}{ll}
      h_l(x) & x \in U_l \text{ and } l > m \\
      x & x \in X \setminus \bigcup_{l > m} U_l
    \end{array}
  \right..
\]
Let $\F{0} = \{ F^0_i = H_0(F^{-1}_i) \}_{i \in I}$. Observe that
$L_0$ is an exact $n$-core of $\F{0}$. We let $K_0 = L_0$. By the
construction, $H_m$ is an approximation within $\F{-1}$ of the
identity on $X$, hence $\F{0}$ refines $\F{-1}$. Hence $\F{0}$ refines
$\st \F{}$.

We shall recursively construct a sequence of closed interior
$\N{n}$-covers $\F{m} = \{ F^m_i \}_{i \in I}$ of $X \setminus
\bigcup_{l > m} V_l$ such that
\begin{enumerate}[(i)]
\item $\F{m}$ is $k$-regular and is isomorphic to $\F{m-1}$,
\item $\xi_m$ is a trivial hole in $F^m_{J_m}$,
\item $\F{m}$ is $\Cup_k$-contractible to $L_m$ in $\st_{\{ U_l \}_{l
      \leq m}} \F{-1}$,
\item $\F{m}$ is equal to $\F{m-1}$ on $X\setminus U_m$ and is its
  swelling on $U_m \setminus V_m$,
\item $\F{m}$ has an exact $n$-core $K_m$ such that $L_m$ is a
  subcomplex of $K_m$ and $K_m \setminus K_{m-1} \subset V_m$.
\end{enumerate}
Let $m > 0$ and assume that we already constructed $\F{m-1}$ and
$K_{m-1}$. Let $\F{m-1/2} = \{ F^{m-1/2}_i \}_{i \in I}$ with
\[
F^{m-1/2}_i = \left\{ \begin{array}{ll}
h^{-1}_m(x) & x \in U_m \\
x & x \in X \setminus U_m
\end{array}\right..
\]
The collection $\F{m-1/2}$ is a closed interior $\N{n}$-cover of $X
\setminus \bigcup_{l > m} V_l$. Let $K_{m-1/2} = K_{m-1} \cup L_m$. By
the construction, $K_{m-1/2}$ is an exact $n$-core of $\F{m-1/2}$.
Obviously, $\F{m-1/2}$ is $k$-regular, is equal to $\F{m-1}$ on $X
\setminus U_m$ and is its swelling on $U_m \setminus V_m$. By (iii),
$\F{m-1}$ is $\Cup_k$-contractible to $L_{m-1}$ in $\st_{\{ U_l \}_{l
    < m}} \F{-1}$. Since $K_{m-1}$ is an exact $n$-core of $\F{m-1}$,
this is equivalent to $\Cup_k$-contractibility of $\F{m-1} / K_{m-1}$
to $L_{m-1}$. Since $K_{m-1}$ is an exact $n$-core of $\F{m-1/2}$ as
well, $\F{m-1/2}$ is $\Cup_k$-contractible to $L_{m-1}$ (and obviously
to $L_m$ as well) in $\st_{\{ U_l \}_{l < m}} \F{-1}$. Hence
conditions of lemma~\ref{lem:patching of a large hole} are satsified
with $\F{} = \F{m-1/2}$, $X = X \setminus \bigcup_{l > m} V_l$, $J =
J_m$, $\varphi = \xi_m$, $\varPhi = \Xi_m$, $L = K_{m-1/2}$, $L_0 =
L_m$, $\mathcal{H} = \st_{\{ U_l \}_{l < m}} \F{-1}$ and $U = U_m$.
Hence there exists a closed interior $\N{n}$-cover $\F{m}$ of $X
\setminus \bigcup_{l > m} V_l$ and an exact $n$-core $K_m$ of $\F{m}$
that satisfy conditions (i)-(v). The construction of the sequence
$\F{m}$ is done.

Take $\F{\infty} = \lim_{m \to \infty} \F{m}$. By the construction
$\U{}$ is a discrete collection in $X$, hence $\F{\infty}$ is a closed
interior $\N{n}$-cover of $X$.  We shall verify that
\begin{enumerate}\addtocounter{enumi}{-1}
\item $\F{\infty}$ is $k$-regular and is isomorphic to $\F{}$,
\item for each $J$ in $\mathcal{J}$, the inclusion of $L \cap
  F^\infty_J$ into $F^\infty_J$ induces a trivial (zero) homomorphism
  on $k$th homotopy groups,
\item $\F{\infty}$ is $\Cup_k$-contractible to $L$ in $\st_\U{}
  \F{-1}$.
\end{enumerate}
If (0), (1) and (2) are satisfied, then by
lemma~\ref{lem:u-contraction and zero homomorphisms}, $\F{\infty}$ is
$\circledcirc_k$-contractible in $\st \st_\U{} \F{0}$. Observe that
$\U{}$ refines $\st^2_\F{} \E{}$, $\F{-1}$ refines $\st \F{}$ and $\st
\F{}$ refines $\st_\F{} \E{}$, hence $\st \st_\U{} \F{-1}$ refines
$\st \st_{\st^2_\F{} \E{}} \st^2_\F{} \E{} = \st \st \st^2_\F{} \E{}$,
which refines $\st \st \st^2 \E{} = \st \st^7 \E{} = \st^{22} \E{}$,
the last equalities by lemma~\ref{lem:star of a star}.  Then by
lemma~\ref{lem:ok contractible} and by remark~\ref{rem:ok
  contractible}, there exists a $(k+1)$-regular closed interior
$\N{n}$-cover that is isomorphic to $\F{\infty}$, that refines $\st
\st^{21\cdot 22 + 10} \E{} = \st^{472} \E{}$ and that is equal to
$\F{\infty}$ on a neighborhood of $Z$. By the construction,
$\F{\infty}$ is isomorphic to $\F{}$ and equal to $\F{}$ on a
neighborhood of $Z$.  Hence verification of (0), (1) and (2) will
finish the proof.

By lemma~\ref{lem:passing to the limit}, $K_\infty = \bigcup_{m \geq
  0} K_m$ is an exact $n$-core of $\F{\infty}$. By the construction,
$K_\infty \cap F^\infty_J = \bigcup K_m \cap F^m_J$. Since $K_m$ is an
exact $n$-core of $\F{m}$ and $\F{m}$ is $k$-regular, each $K_m \cap
F^m_J$ is $l$-connected for each $l < k$. Therefore $\bigcup_{m \geq
  0} K_m \cap F^m_J$ is $l$-connected for each $l < k$. Since
$K_\infty$ is an exact $n$-core of $\F{\infty}$, $\F{\infty}$ is
$k$-regular. Condition (0) is satisfied.

Let $J$ in $\mathcal{J}$. By the construction, the inclusion $L_0 \cap
F^\infty_J \subset L \cap F^\infty_J$ is an $n$-homotopy equivalence.
Hence to prove (1), it suffices to check that the inclusion of $L_0
\cap F^\infty_J$ into $F^\infty_J$ induces a trivial homomorphism on
$k$th homotopy groups. Let $\varphi$ be a map from $S^k$ into $L_0
\cap F^\infty_J$. By the construction, $\varphi$ is homotopic in $L_0
\cap F^\infty_J$ with $\xi_m$ for some $m$ such that $J = J_m$. By
(ii), $\xi_m$ is nullhomotopic in $K_m \cap F^m_{J}$. Hence it is
nullhomotopic in $F^\infty_J$, as $K_m \cap F^m_J \subset F^\infty_J$.

Every element of $\pi_k(F^\infty_J)$ is homotopic in $F^\infty_J$ with
a map into $F^\infty_J \cap K_m$ for some $m$. Hence
$\pi_k(F^\infty_J)$ is generated by maps into $F^m_J \cap K_m$. Such
maps, in turn, are $\Cup_k$-contractible to $L$ in $\st_\U{} \F{-1}$,
by (iii). Hence (2) is satisfied. We are done!


\bibliography{../references}

\newpage
\addcontentsline{toc}{chapter}{Index}
\printindex

\end{document}